# A singular integration by parts formula for the exponential Euclidean QFT on the plane

BY FRANCESCO C. DE VECCHI[abc], MASSIMILIANO GUBINELLI[dbe], MATTIA TURRA[bf]

*a*. Department of Mathematics, University of Pavia
Via Adolfo Ferrata 5, 27100 Pavia, Italy

*b*. Institute for Applied Mathematics &
Hausdorff Center for Mathematics, University of Bonn
Endenicher Allee 60, 53115 Bonn, Germany

*d*. Mathematical Institute, University of Oxford,
Woodstock Road, OX2 6GG Oxford, United Kingdom

*c*. *Email:* `francescocarlo.devecchi@unipv.it`
*e*. *Email:* `gubinelli@maths.ox.ac.uk`
*f*. *Email:* `mattia.turra@iam.uni-bonn.de`

**Abstract**

We give a novel characterization of the Euclidean quantum field theory with exponential interaction $v$ on $\mathbb{R}^2$ through a renormalized integration by parts (IbP) formula, or otherwise said via an Euclidean Dyson–Schwinger equation for expected values of observables. In order to obtain the well-posedness of the singular IbP problem, we import some ideas used to analyse singular SPDEs and we require the measure to "look like" the Gaussian free field (GFF) in the sense that a suitable Wasserstein distance from the GFF is finite. This guarantees the existence of a nice coupling with the GFF which allows to control the renormalized IbP formula.

**A.M.S. subject classification:** 60H17, 81S20, 81T08, 28C20

# Table of contents

# 1 Introduction

One of the first steps in constructive quantum field theory (QFT) is to build a family of distributions satisfying Wightman axioms [GJ87, SW64], which can thus be interpreted as the Wightman functions of a field theory with unique ground state, invariant with respect to the Poincaré group. Wightman axioms, however, do not identify uniquely a particular QFT. Instead different QFTs are expected to satisfy different *Dyson–Schwinger equations*: a systems of partial differential equations (PDEs) relating Wightman functions and encoding the local and hyperbolic equations of motions for the quantum fields [dM97, Sch51].

An important progress in the analysis of Wightman QFTs was the introduction of Schwinger functions, namely the analytic continuations on imaginary time of Wightman functions, which are described by a set of axioms introduced by Osterwalder and Schrader (see [GJ87, OS73, OS75, Sim74]). As observed by Nelson (see [Nel73a, Nel73b, Nel73c]), in many cases (such as the scalar bosonic QFT) Schwinger functions are the moments of a probability measure $\nu$ on Schwarz distributions (see [GJ87, Sim74] for systematic applications of this approach). In particular, the Dyson–Schwinger equations translate in an integration by parts (IbP) formula for the measure $\nu$ [ALZ06, DG74, FR77, GH21] and it becomes natural to investigate the problem of existence and uniqueness of probability measures satisfying prescribed integration by parts formulas.

The characterization of a measure through some IbP formula is a classical subject in stochastic analysis which has different formulations, such as existence and uniqueness of a measure with given logarithmic gradient [Bog10] or the unique closability of a minimally defined pre-Dirichlet form [ABR22, ARZ93, BR95] (and the references therein). The application of logarithmic gradients and integration by parts formulas to quantum field theory was already proposed by Kirillov in the case of sine–Gordon models, see [Kir94a, Kir94b, Kir95a, Kir95b] where the problem is considered without renormalization. In the aforementioned works by Kirillov, the author exploits also a Lyapunov functions technique to show existence of solutions to the integration by parts formulas; a generalization of this technique is applied to the related problem of non-singular (i.e. with no need for renormalization) Fokker–Planck–Kolmogorov equations by Bogachev and Röckner [BR99, BR01]. The problem of uniqueness of solutions to integration by parts formula or, similarly, to an infinite-dimensional Fokker–Planck–Kolmogorov equation, is solved only in some particular cases, see for instance the books [Bog10, BKRS15], and the works by Bogachev, Da Prato and Röckner [BDR09, BDR11] and by Röckner, Zhu and Zhu [RZZ14], where a dissipative non-regular drift (without renormalization) is considered. Let us also mention the studies about uniqueness of solutions to Fokker–Planck–Kolmogorov equations or of invariant measures of the $P(\varphi)_2$ stochastic quantization equation on the two-dimensional torus [DD20, KT22, RZZ17a, RZZ17b, TW18]. Although in this case the problem of existence and uniqueness is solved, it is not clear whether the techniques used in the aforementioned papers can be extended to the models on the non-compact space $\mathbb{R}^2$ or in dimension greater than two. It is worth to mention that the study of uniqueness in the framework of Dirichlet forms for the exponential model has been discussed by Albeverio, Kawabi and Röckner in [AKR12] for the one-dimensional non-singular case, and by the same authors together with Mihalache in [AKMR21] for the two-dimensional setting on the torus. See also [AMR15] for a review of the existing literature on the Dirichlet form approach to the dynamical problem. Let us mention that our concern here is different from that in Dirichlet form theory since here the measure is an unknown in the problem and not given a priori.

The key open problem that we address in this paper is to provide a suitable setting in which existence and uniqueness of measures satisfying given some *singular* IbP formula, that is one involving a renormalization procedure in its definition, which is the usual situation in constructive Euclidean QFT. Instead of attempting a general framework, we concentrate in a particular case where we can establish a reasonable well-posedness theory for the singular IbP problem: we test our ideas on the EQFT with exponential interaction and positive mass $m > 0$, or Høegh–Krohn model [AH74] on the whole space $\mathbb{R}^2$. The exponential interaction in the case of mass $m = 0$ [AGH79, AJPS97, RV14] is a classical example of conformal field theory [Mus10, Sch08] and it finds important applications in Liouville quantum gravity [DS11, KRV20]. As far as we know, this contribution of ours is the first which manages to address this question for an EQFT requiring renormalization and in the infinite volume limit.

Let us give a more detailed description of the problem that we consider here. Let $\mathcal{S}(\mathbb{R}^2)$ be the space of Schwartz functions and denote its dual, that is the space of tempered distributions, by $\mathcal{S}'(\mathbb{R}^2)$. We fix a Banach space (or a topological vector space) $E \subset \mathcal{S}'(\mathbb{R}^2)$ and we consider a family $\mathfrak{F}$ of functions $F: E \to \mathbb{R}$ which are Fréchet differentiable. In particular, for any $\varphi \in E$, we can consider the derivative $D_f F(\varphi)$ of $F$ in the direction $f \in \mathcal{S}(\mathbb{R}^2)$; the map $f \mapsto D_f F(\varphi)$ is linear and bounded in $E$, and thus in $\mathcal{S}(\mathbb{R}^2)$. This means that, since $\mathcal{S}'(\mathbb{R}^2)$ is the topological dual of $\mathcal{S}(\mathbb{R}^2)$, there exists a unique $\nabla_\varphi F(\varphi) \in \mathcal{S}'(\mathbb{R}^2)$ such that $D_f F(\varphi) = \langle \nabla_\varphi F(\varphi), f \rangle$, where $\langle \cdot, \cdot \rangle = \langle \cdot, \cdot \rangle_{\mathcal{S}', \mathcal{S}}$ is the duality between $\mathcal{S}'(\mathbb{R}^2)$ and $\mathcal{S}(\mathbb{R}^2)$.

Let us denote by $\mathcal{M}$ a family of probability measures on $E$, and let $B: E \to \mathcal{S}'(\mathbb{R}^2)$ a given map. Then a generic IbP problem, or Euclidean Schwinger–Dyson equation for a measure $\nu \in \mathcal{M}$, has the general form

$$\int_E \langle \nabla_\varphi F - FB, f \rangle \, d\nu = 0, \qquad \text{for any } F \in \mathrm{Cyl}_E^b.$$

We are interested in *local* functionals $B$ which have the form

$$B(\varphi)(x) = p(\varphi(x)), \qquad x \in \mathbb{R}^2,$$

for some smooth function $p: \mathbb{R} \to \mathbb{R}$. This locality of the IbP formula is peculiar of EQFT where locality (or reflection positivity, or domain Markov property) is structurally linked to the finite speed of propagation of signals in the Minkowski theory. Unfortunately, such kind of functions $B$ are seldom well-defined on the set $E$ on which we could hope that any solution $\nu$ would be supported. Typically, this support looks very much like the support of the GFF and therefore non-linear local functionals are not automatically well-defined and need to be approached via an ultraviolet regularization and subsequent renormalization. In this sense, we talk about a *singular* IbP formula adopting the term from the recent literature in singular stochastic PDEs.

Given this motivation we consider a sequence of maps $(B_\varepsilon)_{\varepsilon > 0}$ such that, for every $\varepsilon > 0$, we have $B_\varepsilon: E \to \mathcal{S}'(\mathbb{R}^2)$ and for which we recover locality in the limit as $\varepsilon \to 0$. They are typically of the form

$$B(\varphi)(x) = p_\varepsilon((g_\varepsilon * \varphi)(x)), \qquad x \in \mathbb{R}^2$$

where $(g_\varepsilon)_{\varepsilon \geqslant 0}$ is some sequence of local smoothing kernels for which $(g_\varepsilon * \varphi) \to \varphi$ in $\mathcal{S}'(\mathbb{R}^d)$ and $(p_\varepsilon: \mathbb{R} \to \mathbb{R})_\varepsilon$ is a sequence of smooth function chosen to deliver the expected renormalization, typical of EQFT in two and three dimensions. The main problem we analyze in this paper can therefore be formulated abstractly as follows:

**Problem A.** *We say that a measure $v \in \mathcal{M}$ satisfies the integration by parts formula with respect to $(B_\varepsilon)_{\varepsilon > 0}$ and $\mathcal{M}$ if, for any $f \in \mathcal{S}(\mathbb{R}^2)$, we have*

$$\int_E \langle \nabla_\varphi F(\varphi), f \rangle \, v(\mathrm{d}\varphi) = \lim_{\varepsilon \to 0} \int_E F(\varphi) \langle B_\varepsilon(\varphi), f \rangle v(\mathrm{d}\varphi), \qquad \text{for any } F \in \mathrm{Cyl}_E^b, \tag{1.1}$$

*where $\mathrm{Cyl}_E^b$ is the set of smooth and bounded cylinder functions (cf. the Notation section at the end of the present introduction).*

Let us remark that the problem of existence and uniqueness of solutions to Problem A strongly depends on the subset $\mathcal{M}$ of the space $\mathcal{P}(E)$ of (Radon) probability measures on $E$. One of the main problems is that, if we consider $\mathcal{M} = \mathcal{P}(E)$, i.e. we consider a generic Radon probability measure on $E$, then it is not clear if $(B_\varepsilon)_\varepsilon$ admits a limit in probability as $\varepsilon \to 0$ and whether the limit depends on the measure $v$ or not. For example, if we take $B_\varepsilon$ to be the drift of $\Phi_2^4$ measure in $\mathbb{R}^2$, namely

$$B_\varepsilon(\varphi) = (g_\varepsilon * \varphi)^4 - 6 c_\varepsilon (g_\varepsilon * \varphi)^2 + 3 c_\varepsilon^2$$

with $c_\varepsilon = \|(-\Delta + m^2)^{-1/2} g_\varepsilon\|^2_{L^2(\mathbb{R}^2)}$, then it is known that $B_\varepsilon$ converges to the unique limit $:\varphi^4:$ (where $:\cdot:$ is the Wick product of Gaussian random fields, see Chapter 1 in [Sim74]) when the measure $v$ is absolutely continuous with respect to the Gaussian free field with mass $m > 0$. On the other hand, if $v$ is supported on the space of smooth functions, such a limit does not exist. It is useful then to consider a class of measures $v$ for which it is possible to make sense (almost surely) of the limit $\lim_{\varepsilon \to 0} B_\varepsilon$.

The set $\mathcal{M}$ can therefore be neither too large nor too small and we find useful to focus on measures which "look like" the Gaussian free field with mass $m > 0$. We formalize this idea introducing an appropriate Wasserstein metric which encodes this proximity without requiring absolute continuity, which in the full $\mathcal{S}'(\mathbb{R}^2)$ is anyway never an appropriate condition. Let $B$ be a convex cone in $E$ equipped with a norm $\|\cdot\|_B$ that is stronger than the one of $E$, then we define the function $d_B: \mathcal{P}(E) \times \mathcal{P}(E) \to [0, +\infty]$ as

$$d_B(v, v') = \inf_{\pi \in \Pi(v, v')} \int_{E \times E} \|x - y\|_B \chi_B(x - y) \, \pi(\mathrm{d}x, \mathrm{d}y), \qquad v, v' \in \mathcal{P}(E),$$

where $\Pi(v, v') \subset \mathcal{P}(E^2)$ is the set of couplings between $v$ and $v'$, and

$$\chi_B(x) := \begin{cases} 1, & \text{if } x \in B \\ +\infty & \text{if } x \notin B \end{cases}.$$

Let us suppose that $E$ contains the support of the measure $v^{\text{free}}$ of the Gaussian free field with mass $m > 0$ and define

$$\mathcal{M}_B := \{ v \in \mathcal{P}(E) \mid d_B(v, v^{\text{free}}) < +\infty \}. \tag{1.2}$$

The space $\mathcal{M}_B$ strongly depends on the convex cone $B$ together with its norm. For example, if we consider $B = H^1(\mathbb{R}^2)$ we have

$$\mathcal{M}_{H^1(\mathbb{R}^2)} = \{ v \in \mathcal{P}(E) \mid v \text{ is absolutely continuous with respect to } v^{\text{free}} \}.$$

Another simple case is when $B = E$ equipped with its natural norm, this gives

$$\mathcal{M}_E = \mathcal{W}^1(E),$$

namely the Wasserstein space of measures on $E$ (see Chapter 6 in [Vil09]).

The class of measures $\mathcal{M}_B$ encodes the existence of sufficiently regular couplings between our target measures and the GFF. This mirrors the situation in singular SPDEs and in other recent development in EQFT where the existence of such couplings has been a key technical aspect to develop a suitable stochastic analysis of singular dynamics and EQFTs. In the fundamental work [DD03], Da Prato and Debussche indeed introduced a notion of solution to stochastic quantization equations in two dimensions as the sum of the Gaussian free field with mass $m$ and a solution to a (random) non-linear PDE. In other words, they look at the solution as a perturbation of the solution to the linear stochastic heat equation. More recently, Barashkov and Gubinelli [BG20, BG21a, BG21b] developed a coupling approach based on an optimal control problem, Bauerschmidt and Hofstetter [BH20] used a similar coupling to study the pathwise properties of the two dimensional sine–Gordon EQFT on a torus and Shenfeld and Mikulincer [She22, MS21] linked these ideas with current developments in the theory of optimal transport and functional inequalities. Similar ideas are fundamental also in the context of singular stochastic PDEs, see regularity structures [CCHS22a, CCHS22b, Hai14, HS16], paracontrolled calculus [AK20, AK21, GH21, GIP15], and renormalization group theory [Duc21, Kup16].

As we already mentioned, we would like to provide complete well-posedness results within this framework of singular IbP formulas and, for this reason, we focus on the specific case of the exponential interaction, namely we take $E = B_X + B_Y$, where $B_X = C_\ell^{-\delta}(\mathbb{R}^2)$, i.e. the (weighted) Besov–Hölder space with negative regularity $-\delta$ (see Appendix A), and $B_Y$ is taken to be

$$B_Y \subset B_{p,p,\ell}^{s(\alpha)-\delta}(\mathbb{R}^2),$$

where $s(\alpha) > 0$ satisfies some conditions depending on the parameter $\alpha$ (see Definition 2.6 below) and $B_{p,p,\ell}^{s(\alpha)-\delta}$ is a weighted Besov space (see Appendix A). Moreover, we set

$$B_\varepsilon(\varphi) := (-\Delta + m^2)\varphi + \alpha f_\varepsilon e^{\alpha(g_\varepsilon * \varphi) - \frac{\alpha^2}{2} c_\varepsilon}, \tag{1.3}$$

where $\alpha, m \in \mathbb{R}_+$, $f_\varepsilon$ is a smooth, spatial cut-off function such that $f_\varepsilon \to 1$, $g_\varepsilon = \varepsilon^{-2} g(\varepsilon^{-1} \cdot)$ is a regular mollifier, and

$$c_\varepsilon := \int_{\mathbb{R}^2} g_\varepsilon(z)(-\Delta + m^2)^{-1} g_\varepsilon(z)\, dz \tag{1.4}$$

is a renormalization constant diverging logarithmically to $+\infty$ as $\varepsilon \to 0$. Finally, we consider the space of measures $\mathcal{M}$ in Problem A to be $\mathcal{M}_{B_Y}$ (see equation (1.2)), that is an intermediate regime between the case $\mathcal{M}_{H^1}$ (i.e. the space of measures that are absolutely continuous with respect to the Gaussian free field) and $\mathcal{M}_E$ (which coincides with the Wasserstein space $\mathcal{W}_E^1$).

**Remark 1.1.** In the following, we shall consider the number $\gamma_{\max} \approx 0.55$ given by the maximum taken over all $r > 1$, satisfying Definition 2.6, of

$$\frac{2(r-1)^2}{r((r-1)^2 + 1)}. \tag{1.5}$$

Moreover, we let $\tilde{\gamma}_{\max} := 3 - 2\sqrt{2} \approx 0.172 < \gamma_{\max}$.

In the present setting, it is possible to obtain the following results.

**Theorem 1.2.** *Suppose that $\alpha^2 < 4\pi\tilde{\gamma}_{\max}$. Consider $\mathcal{M}_{B_Y}$ (see Definition 2.2) with $E = B_X + B_Y$, where*

$$B_X = C_\ell^{-\delta}(\mathbb{R}^2) \quad \text{and} \quad B_Y = B_Y^{\leqslant} := B_{p,p,\ell}^{s(\alpha)-\delta}(\mathbb{R}^2) \cap \{f : \mathbb{R}^2 \to \mathbb{R}, f \leq 0\}.$$

*Then there exists a unique solution to Problem A with respect to $(B_\varepsilon)_\varepsilon$ (given by equation (1.3)) and the space of measures $\mathcal{M}_{B_Y}$.*

It is also possible to obtain an existence result for the whole regime $\alpha^2 < 8\pi$.

**Theorem 1.3.** *Suppose that $\alpha^2 < 8\pi$ and assume that the same hypotheses on the spaces $B_X$, $B_Y$, $E$, $\mathcal{M}_{B_Y}$ and the drift $(B_\varepsilon)_\varepsilon$ as in Theorem 1.2 hold. Then there exists a solution to Problem A with respect to $(B_\varepsilon)_\varepsilon$ and $\mathcal{M}_{B_Y}$.*

Theorem 1.3 is obtained by building suitable Lyapunov functions independent of $\varepsilon > 0$, similarly to what was done by Kirillov for the not-renormalized equation (see [Kir94a, Kir94b, Kir95a, Kir95b]), and reducing the infinite-dimensional problem to the existence of a symmetric invariant measure for a finite-dimensional differential operator.

An important consequence of Theorem 1.2 and Theorem 1.3 is a differential characterization of the exponential measure (see Theorem 2.13 for details).

Let us mention that in the present paper we also get a uniqueness result for a slightly more restrictive formulation of Problem A (cf. Problem B in Section 2.1) in the regime $\alpha^2 < 4\pi\gamma_{\max}$. The possible measures solving this latter formulation of the IbP problem contain the invariant measure of the stochastic quantization equation with exponential interaction.

## Plan of the paper

Let us present here the structure of the paper. In Section 2, we discuss the IbP Problem A in the general setting, establishing some equivalent formulations that will be useful to address Problem A itself. We also consider the more restrictive Problem B and Problem B-sym, their relation with Problem A and some properties of the solutions to such problems, such as the negativity of the coupling between the Gaussian free field and the quantum field measure with exponential interaction. Section 3 is devoted to the study of uniqueness of solutions to Problem A showing the uniqueness result stated in Theorem 1.2. Since the proof relies on some properties of the solution to the resolvent equation associated to the drift $B_\varepsilon$, most parts of the section is dedicated to the study of such an object. The existence of solutions to Problem A is proved in Section 4. The proof is based on an approximation method and it involves Lyapunov functions. In Appendices A and B we recall the definitions and properties of weighted Besov spaces and Wick exponential, respectively, which are used throughout the paper. Appendices C, D, and E are concerned with some technical, analytical results on (S)PDEs exploited in the paper.

## Notation

We fix here some notation adopted throughout the paper.

We write $a \lesssim b$ or $b \gtrsim a$ if there exists a constant $C > 0$, independent of the variables under consideration, such that $a \leq C b$, and $a \simeq b$ if both $a \lesssim b$ and $b \lesssim a$. If the aforementioned constant $C$ depends on a variable, say $C = C(x)$, then we use the notation $a \lesssim_x b$, and similarly for $\gtrsim$.

The space $L(A, B)$ is the space of linear and bounded functionals from the Banach space $A$ into the Banach space $B$. $L(A, B)$ is equipped with its natural operator norm $\|\cdot\|_{L(A,B)}$. We also write $L(A) = L(A, A)$.

The space of Schwartz functions on $\mathbb{R}^d$ is denoted by $\mathcal{S}(\mathbb{R}^d)$ and its dual, that is the space of tempered distributions, is denoted by $\mathcal{S}'(\mathbb{R}^d)$. Moreover, we let $C_b^m(\mathbb{R}^d)$ be the space of $m$-times differentiable functions on $\mathbb{R}^d$ with continuous and bounded derivatives. We also use the notation $B_{p,q,\ell}^r(\mathbb{R}^d)$ to denote the weighted Besov space on $\mathbb{R}^d$ (see Appendix A for more results on Besov spaces). For $k, \ell > 0$, we define the weight $\rho_\ell^k(x) = (1 + k|x|^2)^{-\ell/2}$, $x \in \mathbb{R}^d$, and $\rho_\ell := \rho_\ell^1$.

Let $K$ be a convex subset of $\mathcal{S}'(\mathbb{R}^d)$, we define the set $\mathrm{Cyl}_K$ of *cylinder functions* on $K$ to be the set of functions $F: K \to \mathbb{R}$ such that there exist a function $\tilde{F} \in C^2(\mathbb{R}^n)$, with at most linear growth at infinity and all bounded derivatives, and $u_1, \ldots, u_n \in \mathcal{S}(\mathbb{R}^d)$, such that $F(\kappa) = \tilde{F}(\langle u_1, \kappa \rangle, \ldots, \langle u_n, \kappa \rangle)$, for every $\kappa \in K$. We say that $F \in \mathrm{Cyl}_K$ has compact support in Fourier variables if the Fourier transforms $\hat{u}_j$ of $u_j$, $j = 1, \ldots, n$, have compact supports on $\mathbb{R}^d$. Moreover, we denote by $\mathrm{Cyl}_K^b$ the subset of $\mathrm{Cyl}_K$ such that, sticking with the previous notation, the function $\tilde{F}$ is also bounded.

We adopt the notation $P_t = e^{-(-\Delta + m^2)t}$ for the heat kernel with mass $m$ (cf. Appendix A).

## Acknowledgements

The authors thank Sergio Albeverio for the kind suggestions and remarks on various early drafts of the paper. This work has been partially funded by the German Research Foundation (DFG) under Germany's Excellence Strategy - GZ 2047/1, project-id 390685813. The first author is partially funded by the Istituto Nazionale di Alta Matematica "Francesco Severi" (INdAM), Gruppo Nazionale per l'Analisi Matematica, la Probabilità e le loro Applicazioni (GNAMPA): "Analisi armonica e stocastica in problemi di quantizzazione e integrazione funzionale".

# 2 Formulation of the problem

## 2.1 Reformulation of Problem A

In this section, we want to reformulate Problem A when the space of measures $\mathcal{M}$ is the space $\mathcal{M}_B$ introduced in (1.2). For this reason, we consider $E = B_X + B_Y$, where $B_X$ is some space supporting the measure of the Gaussian free field with mass $m > 0$, and $B_Y = B$. We consider also the natural projection $P^X: B_X \times B_Y \to B_X$ and the map $P^{X+Y}: B_X \times B_Y \to E$ such that $(X, Y) \mapsto X + Y$.

Let us reformulate Problem A by means of a second order operator. Consider $B_\varepsilon$ to be regular enough and let

$$\mathfrak{L}_\varepsilon F := \frac{1}{2} \operatorname{tr}_{L^2(\mathbb{R}^2)}(\nabla_\varphi^2 F) - \langle B_\varepsilon, \nabla_\varphi F \rangle, \qquad F \in \mathrm{Cyl}_E.$$

**Problem A'.** *We say that a measure $\nu \in \mathcal{M}$ satisfies the* symmetric Fokker–Planck–Kolmogorov *equation related to $B_\varepsilon$ if*

$$\lim_{\varepsilon \to 0} \int [(\mathfrak{L}_\varepsilon F)\, G - F\, (\mathfrak{L}_\varepsilon G)]\, \mathrm{d}\nu = 0, \qquad \text{for any } F \in \mathrm{Cyl}_E^b, G \in \mathrm{Cyl}_E. \tag{2.1}$$

**Remark 2.1.** A consequence of equation (2.1) is that

$$\lim_{\varepsilon \to 0} \int \mathfrak{L}_\varepsilon F \mathrm{d}\nu = 0, \qquad \text{for any } F \in \mathrm{Cyl}_E. \tag{2.2}$$

That is what is usually called *Fokker-Planck-Kolmogorov equation*.

We now want to lift the problem from the space $E$ to the space $B_X \times B_Y$. In order to do so, we introduce the following notion:

**Definition 2.2.** *The subset of measures $\mathcal{M}$ satisfies the* coupling hypotheses *if, for any $\nu \in \mathcal{M}$, there exists a Radon measure $\mu$ on $B_X \times B_Y$ with the following properties:*

  i. $P_*^X \mu = \mu^{\mathrm{free}}$, *where $\mu^{\mathrm{free}}$ is the law of the Gaussian free field on $B_X$,*

  ii. $P_*^{X+Y} \mu = \nu$,

  iii. $\int \|Y\|_{B_Y} \mu(\mathrm{d}X, \mathrm{d}Y) < +\infty$,

*we call $\mu$ a coupling of $\nu$ with the free field. We denote by $\mathcal{M}_{B_X \times B_Y}$ the set of Radon measures on $B_X \times B_Y$ satisfying condition i. and iii.*

**Remark 2.3.** With the same notation as in Definition 2.2, it is clear that $\mathcal{M}_{B_Y}$ defined as in (1.2) where $B = B_Y$, coincides with the space $P_*^{X+Y} \mathcal{M}_{B_X \times B_Y}$. Indeed, $\nu \in \mathcal{M}_{B_Y}$ if and only if there exists a coupling $\pi(\mathrm{d}x, \mathrm{d}z)$ between $\nu(\mathrm{d}z)$ and the free field $\mu^{\mathrm{free}}(\mathrm{d}x)$ such that the difference $x - z$ is supported on $B_Y$. The coupling $\pi$ is related with a measure $\mu$ in $\mathcal{M}_{B_X \times B_Y}$ via the transformation $y = x - z$.

We consider now an operator $\mathscr{L}_\varepsilon$ on the space of regular functions on $B_X \times B_Y$ of the form

$$\mathscr{L}_\varepsilon \Phi(X, Y) := \frac{1}{2} \mathrm{tr}(\nabla_X^2 \Phi) - \langle (-\Delta + m^2)X, \nabla_X \Phi \rangle - \langle B_\varepsilon(X + Y) - (-\Delta + m^2)Y, \nabla_Y \Phi \rangle.$$

**Problem A''.** *We say that a measure $\nu \in \mathcal{M}_{B_Y}$ satisfies the* symmetric Fokker–Planck–Kolmogorov *equation related to $B_\varepsilon$ if*

$$\lim_{\varepsilon \to 0} \int \left[ \mathscr{L}_\varepsilon(F \circ P^{X+Y}) \, G \circ P^{X+Y} - F \circ P^{X+Y} \mathscr{L}_\varepsilon(G \circ P^{X+Y}) \right] \mu(\mathrm{d}X, \mathrm{d}Y) = 0, \quad \text{for any } F \in \mathrm{Cyl}_E^b, G \in \mathrm{Cyl}_E, \tag{2.3}$$

*where $\mu$ is a coupling of $\nu$ with the free field.*

**Remark 2.4.** As in the case of Problem A' (see also Remark 2.1), equation (2.3) implies

$$\lim_{\varepsilon \to 0} \int \mathscr{L}_\varepsilon(F \circ P^{X+Y}) \, \mu(\mathrm{d}X, \mathrm{d}Y) = 0, \qquad \text{for any } F \in \mathrm{Cyl}_E, \tag{2.4}$$

We often say that $\mu$ solves Problem A'' if there is a measure $\nu$ solving the FPK equation with $\mu$ as coupling of $\nu$ with the free field.

Problem A'' is based on a formulation of integration by parts formula related to stationary solutions to the Fokker–Planck–Kolmogorov equation, making the operator $\mathscr{L}_\varepsilon$ symmetric on its domain (see Problem A' in Section 2.1 or [ABR99, BR95] for more details on the relations between the formulations).

It is convenient to introduce new equations related to Problem A'' described above, where the argument of $\mathscr{L}_\varepsilon$ in equation (2.4) is not necessarily of the form $F \circ P^{X+Y}$.

**Problem B.** *We say that a Radon measure $\mu$ on $B_X \times B_Y$ satisfies the* Fokker–Planck–Kolmogorov *equation related to $B_\varepsilon$, if*

$$\lim_{\varepsilon \to 0} \int \mathscr{L}_\varepsilon \Phi \, d\mu = 0, \qquad \text{for any } \Phi \in \mathrm{Cyl}_{B_X \times B_Y}. \tag{2.5}$$

**Problem B-sym.** *We say that $\mu$ satisfies the* symmetric Fokker–Planck–Kolmogorov equation *related to $B_\varepsilon$ if, furthermore,*

$$\lim_{\varepsilon \to 0} \int \left[ \mathscr{L}_\varepsilon(F \circ P^{X+Y}) \, G \circ P^{X+Y} - F \circ P^{X+Y} \, \mathscr{L}_\varepsilon(G \circ P^{X+Y}) \right] d\mu = 0, \tag{2.6}$$

*for any $F, G: B_X + B_Y \to \mathbb{R}$ such that $F \circ P^{X+Y} \in \mathrm{Cyl}_E^b$, $G \circ P^{X+Y} \in \mathrm{Cyl}_E$.*

We have then the following result.

**Theorem 2.5.** *Suppose that $\sup_{\varepsilon > 0} \int |\langle B_\varepsilon(\varphi), \nabla_\varphi F \rangle| \, \nu(d\varphi) < +\infty$, for any $F \in \mathrm{Cyl}_E$ and $\nu \in \mathscr{M}$. Then the following statements hold:*

  i. *Problem A is equivalent to Problem A'.*

  ii. *Assume further that $\mathscr{M} = \mathscr{M}_{B_Y}$, then Problem A, Problem A' and Problem A'' are equivalent.*

  iii. *Assume the same hypothesis as in point ii.. Then, a solution to Problem B-sym is a solution to Problem A''.*

*Summarizing, if all the hypotheses of points i., ii., and iii., then Problem A is equivalent to Problem A'', and Problem B implies Problems A and A'', that is*

$$\text{Problem A} \iff \text{Problem A'} \iff \text{Problem A''} \impliedby \text{Problem B-sym}$$

**Proof.** Let us give the proof point by point.

**Proof of point i.** Let us prove that a solution to Problem A' is also a solution to Problem A. Let $\nu \in \mathscr{M}$ solve Problem A', take $F \in \mathrm{Cyl}_E^b$, and consider $G(\cdot) = \langle \cdot, f \rangle$, for $f \in \mathscr{S}(\mathbb{R}^2)$. Then, by equation (2.2), we have

$$\lim_{\varepsilon \to 0} \int \mathfrak{L}_\varepsilon(F \cdot G) \, d\nu = 0,$$

where

$$\begin{aligned}\mathfrak{L}_\varepsilon(F \cdot G) &= \mathfrak{L}_\varepsilon(F) \, G + F \mathfrak{L}_\varepsilon(G) + 2\langle \nabla_\varphi F, \nabla_\varphi G \rangle \\ &= \mathfrak{L}_\varepsilon(F) \, G + F \mathfrak{L}_\varepsilon(G) + 2\langle \nabla_\varphi F, f \rangle.\end{aligned}$$

Thus, we have

$$\begin{aligned}
0 &= \lim_{\varepsilon\to 0}\int \mathfrak{L}_\varepsilon(F\cdot G)\mathrm{d}\nu = \lim_{\varepsilon\to 0}\int (\mathfrak{L}_\varepsilon(F)\,G + F\mathfrak{L}_\varepsilon(G) + 2\langle \nabla_\varphi F, f\rangle)\mathrm{d}\nu \\
&= \lim_{\varepsilon\to 0}\int (2F\mathfrak{L}_\varepsilon(G) + 2\langle \nabla_\varphi F, f\rangle)\mathrm{d}\nu = \lim_{\varepsilon\to 0}\int (2F\langle B_\varepsilon, f\rangle + 2\langle \nabla_\varphi F, f\rangle)\mathrm{d}\nu.
\end{aligned}$$

We now show that any solution to Problem A is a solution to Problem A'. Let $\mu\in\mathcal{M}$ solve Problem A. Consider an orthonormal basis $(e_n)_n$ of $L^2(\mu)$ and take $F = \langle \nabla G, e_n\rangle$, for some $G\in\mathrm{Cyl}_E$. We have then, by equation (1.1),

$$\int \nabla^2 G(e_n, e_n)\mathrm{d}\nu = \int \langle \nabla F, e_n\rangle\mathrm{d}\nu = \int \langle \nabla G, e_n\rangle\langle B_\varepsilon, e_n\rangle\mathrm{d}\nu.$$

Taking the sum over $n\in\mathbb{N}$ on both sides, noticing that by the properties of cylinder functions we can exchange the integral with the sum, and exploiting Parseval's identity we have

$$\int \mathrm{tr}_{L^2(\mathbb{R}^2)}(\nabla^2 G)\mathrm{d}\nu = \int \langle \nabla G, B_\varepsilon\rangle\mathrm{d}\nu,$$

that is $\int \mathfrak{L}_\varepsilon G \mathrm{d}\nu = 0$. Now, consider again an orthonormal basis $(e_n)_n$ of $L^2(\mu)$, but take $F = G\langle \nabla H, e_n\rangle$, for some $G, H\in\mathrm{Cyl}_E$. Let us note that

$$\langle \nabla F, e_n\rangle = \langle \nabla (G\langle \nabla H, e_n\rangle), e_n\rangle = \langle \nabla G, e_n\rangle\langle \nabla H, e_n\rangle + G\nabla^2 H(e_n, e_n).$$

On the other hand, equation (1.1) gives

$$\int \langle \nabla (G\langle \nabla H, e_n\rangle), e_n\rangle\mathrm{d}\nu = \int G\langle \nabla H, e_n\rangle\langle B_\varepsilon, e_n\rangle\mathrm{d}\nu.$$

This yields

$$\int G\langle \nabla H, e_n\rangle\langle B_\varepsilon, e_n\rangle\mathrm{d}\nu = \int (\langle \nabla G, e_n\rangle\langle \nabla H, e_n\rangle + G\nabla^2 H(e_n, e_n))\mathrm{d}\nu,$$

and taking the sum over $n\in\mathbb{N}$ and using again the properties of cylinder functions, we get

$$\int G\langle \nabla H, B_\varepsilon\rangle\mathrm{d}\nu = \int (\langle \nabla G, \nabla H\rangle + G\mathrm{tr}_{L^2(\mathbb{R}^2)}(\nabla^2 H))\mathrm{d}\nu,$$

which implies

$$0 = \int (\langle \nabla G, \nabla H\rangle + G\mathrm{tr}(\nabla^2 H) - G\langle \nabla H, B_\varepsilon\rangle)\mathrm{d}\nu = \int (\langle \nabla G, \nabla H\rangle + G\mathfrak{L}_\varepsilon H)\mathrm{d}\nu,$$

namely

$$\int G\mathfrak{L}_\varepsilon(H)\mathrm{d}\nu = -\int \langle \nabla G, \nabla H\rangle\mathrm{d}\nu.$$

Doing the same computation exchanging the roles of $G$ and $H$ we also get $\int H\mathfrak{L}_\varepsilon(G)\nu(\mathrm{d}\varphi) = -\int \langle \nabla G, \nabla H\rangle \nu(\mathrm{d}\varphi)$, that gives

$$\int G\mathfrak{L}_\varepsilon(H)\mathrm{d}\nu - \int H\mathfrak{L}_\varepsilon(G)\mathrm{d}\nu = 0,$$

and concludes the proof of the first point of the theorem.

**Proof of point ii.** Assume further that $\mathcal{M} = \mathcal{M}_{B_Y}$. Notice that $\mathscr{L}_\varepsilon(F \circ P^{X+Y}) = \mathfrak{L}_\varepsilon(F) \circ P^{X+Y}$, and since $P_*^{X+Y}\mu = \nu$ we have that

$$\int \mathscr{L}_\varepsilon(F \circ P^{X+Y})\,\mathrm{d}\mu = \int \mathfrak{L}_\varepsilon(F) \circ P_*^{X+Y}\,\mathrm{d}\mu = \int \mathfrak{L}_\varepsilon F\,\mathrm{d}\nu,$$

which gives the equivalence between equations (2.4) and (2.2) and in the same way we have the equivalence between equations (2.3) and (2.1).

**Proof of point iii.** Let $\mu$ be a solution to Problem B. We only have to show that $\mu \in \mathcal{M}_{B_X \times B_Y}$, that is, the only condition that we need to verify is that $P_*^X \mu = \mu^{\text{free}}$.

Since the terms involving derivatives with respect to $X$ of $\mathscr{L}_\varepsilon$ do not depend on $\varepsilon$, in this proof we write $\mathscr{L}$ for the operator $\mathscr{L}_\varepsilon$. Furthermore, for a solution $\mu$ to (2.5)–(2.6) we have that, if $\Phi$ does not depend on $Y$,

$$0 = \lim_{\varepsilon \to 0} \int \mathscr{L}_\varepsilon \Phi(X)\mu(\mathrm{d}X, \mathrm{d}Y) = \int \mathscr{L}_\varepsilon \Phi(X)\mu(\mathrm{d}X, \mathrm{d}Y) = \int \mathscr{L}\Phi(X)\mu(\mathrm{d}X, \mathrm{d}Y).$$

*Step 1.* We start with the proof that the free field measure with mass $m > 0$ solves equations (2.5)–(2.3).

Consider, for $X_0 \in B_X$,

$$\mu_t \sim X_t = P_t X_0 + \int_0^t P_{t-s} \xi_s\,\mathrm{d}s,$$

and let $\mu_\infty$ be the limit

$$\mu_\infty(\mathrm{d}X) = \lim_{T \to +\infty} \frac{1}{T} \int_0^T \mu_t(\mathrm{d}X)\,\mathrm{d}t.$$

We want to show that $\int \mathscr{L}\Phi\,\mathrm{d}\mu_\infty = 0$, for any $\Phi$ smooth cylindric function. By Itô's formula (see Theorem 4.32 in [DZ14]), taking expectation and dividing by $T > 0$, we have that

$$\begin{aligned}
\frac{1}{T}\mathbb{E}[\Phi(X_T) - \Phi(X_0)] &= \frac{1}{T}\int_0^T \mathbb{E}[\mathscr{L}\Phi(X_s)]\,\mathrm{d}s \\
&= \frac{1}{T}\int_0^T \int \mathscr{L}\Phi(X)\mu_s(\mathrm{d}X)\,\mathrm{d}s \\
&= \int \mathscr{L}\Phi(X)\Big(\frac{1}{T}\int_0^T \mu_s(\cdot)\,\mathrm{d}s\Big)(\mathrm{d}X),
\end{aligned}$$

and letting $T \to +\infty$, recalling that $\Phi$ is bounded, we get

$$0 = \int \mathscr{L}\Phi(X)\Big(\lim_{T \to +\infty} \frac{1}{T}\int_0^T \mu_s(\cdot)\,\mathrm{d}s\Big)(\mathrm{d}X) = \int \mathscr{L}\Phi(X)\mu_\infty(\mathrm{d}X).$$

*Step 2.* We now show uniqueness of the measure $\mu$. Consider

$$G^\lambda(X_0) = \mathbb{E}_\xi \int_0^{+\infty} e^{-\lambda s} F(X_s)\,\mathrm{d}s,$$

where $F \in \mathrm{Cyl}_{B_X}$ of the form

$$F(X) = \tilde{F}(\langle u_1, X\rangle, \ldots, \langle u_n, X\rangle).$$

*Step 2.a.* Let us prove that $\mathscr{L}_\varepsilon G^\lambda$ is well-defined, recall

$$\mathscr{L}_\varepsilon G^\lambda(X_0) = \mathscr{L} G^\lambda(X_0) = \frac{1}{2}\operatorname{tr}(\nabla^2_{X_0} G^\lambda) - \langle(-\Delta + m^2)X_0, \nabla_{X_0} G^\lambda\rangle,$$

where

$$\nabla_{X_0} G^\lambda = \mathbb{E}_\xi \int_0^{+\infty} e^{-\lambda s}\langle\nabla_{X_0} F(X_s), \nabla_{X_0} X_s\rangle\, ds = \mathbb{E}_\xi \int_0^{+\infty} e^{-\lambda s}\sum_{k=1}^n \partial_k \tilde{F}(\langle u_1, X\rangle, \ldots, \langle u_n, X\rangle)\, P_t u_k\, ds,$$

and

$$\begin{aligned}\nabla^2_{X_0} G^\lambda &= \mathbb{E}\int_0^{+\infty} e^{-\lambda s}\nabla_{X_0}[\langle\nabla_{X_0} F(X_s), \nabla_{X_0} X_s\rangle]\, ds \\ &= \mathbb{E}_\xi \int_0^{+\infty} e^{-\lambda s}\sum_{k,j=1}^n \partial^2_{k,j}\tilde{F}(\langle u_1, X\rangle, \ldots, \langle u_n, X\rangle)\langle P_t u_k, P_t u_j\rangle\, ds,\end{aligned}$$

since

$$\nabla_{X_0} X_t(h) = P_t h, \qquad \nabla^2_{X_0} X_t(h) = 0.$$

This gives

$$\begin{aligned}\mathscr{L} G^\lambda(X_0) &= \frac{1}{2}\operatorname{tr}_{L^2(\mathbb{R}^2)}\Big(\mathbb{E}_\xi \int_0^{+\infty} e^{-\lambda s}\sum_{k,j=1}^n \partial^2_{k,j}\tilde{F}(\langle u_1, X\rangle, \ldots, \langle u_n, X\rangle)\langle P_t u_k, P_t u_j\rangle\, ds\Big) \\ &\quad - \mathbb{E}_\xi \int_0^{+\infty} e^{-\lambda s}\sum_{k=1}^n \partial_k \tilde{F}(\langle u_1, X\rangle, \ldots, \langle u_n, X\rangle)\langle(-\Delta + m^2)X_0, P_t u_k\rangle\, ds,\end{aligned}$$

and thus $\mathscr{L} G^\lambda$ is well-defined.

*Step 2.b.* We now show that

$$\int (\lambda - \mathscr{L}) G^\lambda\, d\mu = \int F\, d\mu. \tag{2.7}$$

By Itô formula (see Theorem 4.32 in [DZ14]), we have

$$e^{-\lambda t} G^\lambda(X_t) - G^\lambda(X_0) = \int_0^t e^{-\lambda s}\mathscr{L}_\varepsilon G^\lambda(X_s)\, ds + \int_0^t e^{-\lambda s}\nabla_{X_0} G^\lambda(X_s)\, dX_s - \lambda\int_0^t e^{-\lambda s} G^\lambda(X_s)\, ds,$$

and taking the expectation we get

$$e^{-\lambda t}\mathbb{E}[G^\lambda(X_t)] - \mathbb{E}[G^\lambda(X_0)] = \mathbb{E}\int_0^t e^{-\lambda s}\mathscr{L}_\varepsilon G^\lambda(X_s)\, ds - \lambda\mathbb{E}\int_0^t e^{-\lambda s} G^\lambda(X_s)\, ds. \tag{2.8}$$

On the other hand,

$$e^{-\lambda t} G^\lambda(X_t) - G^\lambda(X_0) = \int_0^{+\infty} e^{-\lambda(t+s)}\mathbb{E}_\xi[F(X_{t+s})]\, ds - \int_0^{+\infty} e^{-\lambda s}\mathbb{E}_\xi[F(X_s)]\, ds,$$

and taking the expectation

$$e^{-\lambda t}\mathbb{E}[G^\lambda(X_t)] - \mathbb{E}[G^\lambda(X_0)] = -\mathbb{E}\int_0^t e^{-\lambda s}\mathbb{E}_\xi[F(X_s)]\, ds. \tag{2.9}$$

Comparing equations (2.8) and (2.9), we get

$$\lambda \mathbb{E}\int_0^t e^{-\lambda s} G^\lambda(X_s)\,ds - \mathbb{E}\int_0^t e^{-\lambda s} \mathscr{L}_\varepsilon G^\lambda(X_s)\,ds = \mathbb{E}\int_0^t e^{-\lambda s}\mathbb{E}_\xi[F(X_s)]\,ds,$$

dividing by $t>0$ and letting $t\to 0$, we have the result.

*Step 2.c.* Take $Z\sim\mu$ and $Z'\sim\mu_\infty$, and call $\psi^Z$ and $\psi^{Z'}$ the solutions to equation

$$(\partial_t - \Delta + m^2)\psi = \xi,$$

with initial conditions $\psi(0)=Z$ and $\psi(0)=Z'$, respectively. We show that $\mu\sim\mu_\infty$.

Consider the product measure $\mu\otimes\mu_\infty$ between $\mu$ and $\mu_\infty$. Since we are working in separable spaces, any Borel measure is a Radon measure (see, e.g. Theorem 9 in Chapter 2, Section 3 in [Sch73]), and therefore we have that, for every $\varepsilon>0$, there exists $r>0$ such that

$$\mu\otimes\mu_\infty(\|Z\|>r, \|Z'\|>r) < \varepsilon. \tag{2.10}$$

Since

$$\begin{aligned}|\tilde{F}(\langle u_1,\psi^Z\rangle,\ldots,\langle u_n,\psi^Z\rangle) - \tilde{F}(\langle u_1,\psi^{Z'}\rangle,\ldots,\langle u_n,\psi^{Z'}\rangle)| &\leq \|\nabla\tilde{F}\|_{L^\infty}\sum_{k=1}^n |\langle u_k,\psi^Z - \psi^{Z'}\rangle| \\ &\leq \|\nabla\tilde{F}\|_{L^\infty}\sum_{k=1}^n \|u_k\|\|P_t(Z-Z')\| \\ &\leq C\lambda\int_0^{+\infty} e^{-\lambda t}e^{-m^2 t}(\|Z\|+\|Z'\|)\,dt \\ &\leq C'\frac{\lambda}{\lambda+m^2}(\|Z\|+\|Z'\|),\end{aligned}$$

then we have, recalling $\int \lambda G_\lambda\,d\mu_\infty = \int F\,d\mu_\infty$,

$$\int \lambda G_\lambda\,d\mu - \int F\,d\mu_\infty = \mathbb{E}[\lambda G_\lambda(Z) - \lambda G_\lambda(Z')] \leq C'\frac{\lambda}{\lambda+m^2}\mathbb{E}(\|Z\|+\|Z'\|).$$

The right-hand side of the previous inequality converges to zero as $\lambda\to+\infty$ by property (2.10), and thus $\int F\,d\mu_\infty = \int F\,d\mu$ since, by equation (2.7) and being $\mu$ a solution to Problem B, we have $\lambda\int G_\lambda\,d\mu = \int F\,d\mu$. This yields $\mu\sim\mu_\infty$ and uniqueness of the measure. □

## 2.2 A priori deductions for exponential interaction

Let $\alpha, m\in\mathbb{R}_+$, and $\gamma := \alpha^2/(4\pi)$.

**Definition 2.6.** *We define*

$$B_X := C_\ell^{-\delta}(\mathbb{R}^2),$$

*for some $\delta>0$ small enough, and with $\ell>0$. Let us recall that $\ell$ denotes the presence of a weight, in fact $C_\ell^\sigma(\mathbb{R}^2)$ and $B_{p,q,\ell}^\sigma(\mathbb{R}^\sigma)$ are weighted Besov spaces (see Appendix A). Also, we let $B_Y$ to be a space of Besov functions with positive regularity, in particular, we choose*

$$B_Y \subset B_{p,p,\ell}^{s-\delta}(\mathbb{R}^2),$$

*where $1 < p < +\infty$ and $s > 0$ satisfies the following condition depending on $\gamma$:*

$$0 < s < \gamma + 2 - \sqrt{8\gamma},$$

*arising when dealing with the problem of existence of solutions to the Fokker–Planck–Kolmogorov equation (cf. Theorem 4.4). We suppose further that there exists $r > 1$ such that*

$$\frac{1}{p} + \frac{1}{r} \leq 1, \qquad s - \gamma(r-1) - 2\delta > 0, \qquad \gamma r < 2.$$

*Notice that such a condition is always satisfied for some $s > 0$, $p > 1$ and $r > 1$, whenever $\gamma < 2$. Let us also recall the definition of the space*

$$B_Y^\leq := B_{p,p,\ell}^{s-\delta}(\mathbb{R}^2) \cap \{f : \mathbb{R}^2 \to \mathbb{R}, f \leq 0\} \tag{2.11}$$

*featured in Theorem 1.2.*

Let us introduce the functional $\mathscr{G}_\varepsilon : B_X \times B_Y \to \mathscr{S}'(\mathbb{R}^2)$

$$\mathscr{G}_\varepsilon(X, Y) := \alpha f_\varepsilon : e^{\alpha(g_\varepsilon * X)} : e^{\alpha(g_\varepsilon * Y)}. \tag{2.12}$$

where the term $:e^{\alpha(g_\varepsilon * X)}:$ is defined as

$$:e^{\alpha(g_\varepsilon * X)}: = e^{\alpha(g_\varepsilon * X) - \frac{\alpha^2}{2} c_\varepsilon}, \tag{2.13}$$

with the constant $c_\varepsilon$ introduced in equation (1.4). Here, $f_\varepsilon$ is a smooth cut-off function with compact support and with derivatives uniformly bounded in $\varepsilon$ such that $f_\varepsilon \to 1$ uniformly with respect to a polynomially weighted norm, and we let $g_\varepsilon = \varepsilon^{-2} g(\varepsilon^{-1} \cdot)$, where $g$ is a positive, smooth, compactly supported function on $\mathbb{R}^2$ with Lebesgue integral equal to 1, and such that there exists its convolutional square root $\tilde{g}_\varepsilon$, i.e. $g_\varepsilon = \tilde{g}_\varepsilon * \tilde{g}_\varepsilon$, which is also positive, smooth, and compactly supported function with $\int \tilde{g}_\varepsilon(x) \, dx = 1$. We also assume that $g_\varepsilon(x) = g_\varepsilon(-x)$, for all $x \in \mathbb{R}^2$. Notice that this last property implies that $\hat{g}_\varepsilon$ takes real values.

Consider an operator $\mathscr{L}_\varepsilon$ of the form

$$\mathscr{L}_\varepsilon \Phi(X, Y) = \frac{1}{2} \operatorname{tr}(\nabla_X^2 \Phi) - \langle (-\Delta + m^2) X, \nabla_X \Phi \rangle - \langle (-\Delta + m^2) Y + \mathscr{G}_\varepsilon(X, Y), \nabla_Y \Phi \rangle, \tag{2.14}$$

which is well-defined for $\Phi : B_X \times B_Y \to \mathbb{R}$ in a suitable class $\mathscr{F}$ of regular functions to be specified below (see Definition 3.1). For the moment, we only require the set $\operatorname{Cyl}_{B_X \times B_Y}$ of smooth cylindrical functions to be contained in the family $\mathscr{F}$.

**Remark 2.7.** Our choice of $B_X$ and $B_Y$ implies that the space $B_X$ contain the Gaussian free field with mass $m$ and that

$$\begin{aligned}
(-\Delta + m^2) X &\in C_\ell^{-2-\delta}(\mathbb{R}^2), \\
(-\Delta + m^2) Y &\in B_{p,p,\ell}^{s-2-\delta}(\mathbb{R}^2), \\
:\exp(\alpha X): &\in B_{r,r,\ell'}^{-\gamma(r-1)-\delta}(\mathbb{R}^2), \qquad r\ell' > 2, \\
\exp(\alpha Y) &\in L^\infty(\mathbb{R}^2) \cap B_Y, \qquad \text{if } Y \leq 0,
\end{aligned}$$

where all the parameters are as in Definition 2.6 (see [ADG21, HKK21, HKK22] and Appendix B). For simplicity, we adopt the notation

$$B_{\exp}^{\tilde{s},\tilde{\ell}} = B_{\tilde{s},\tilde{s},\tilde{\ell}}^{-\gamma(\tilde{s}-1)-\delta}(\mathbb{R}^2),$$

and we omit the parameters whenever $\tilde{s} = r$ and $\tilde{\ell} = \ell'$, i.e. $B_{\exp} = B_{\exp}^{\tilde{s},\tilde{\ell}}$.

Let us introduce also the notation

$$\mathscr{G} = \mathscr{G}(X,Y) := \alpha \!:\! e^{\alpha X} \!:\! e^{\alpha Y},$$

where the measurable function $:\exp(\alpha X):$ is defined on a subset of $B_X$ of full measure (with respect to the free field measure with mass $m$) as the limit of $:\exp(\alpha g_\varepsilon * X):$ as $\varepsilon \to 0$ (see Proposition B.1 for details). The operator $\mathscr{L}_\varepsilon$ defined in equation (2.14) is an approximation for the operator

$$\mathscr{L}\Phi := \frac{1}{2}\mathrm{tr}_{L^2(\mathbb{R}^2)}(\nabla_X^2 \Phi) - \langle (-\Delta + m^2)X, \nabla_X \Phi \rangle - \langle (-\Delta + m^2)Y + \mathscr{G}(X,Y), \nabla_Y \Phi \rangle, \tag{2.15}$$

for $\Phi \in \mathrm{Cyl}_{B_X \times B_Y}$.

**Theorem 2.8.** *Let $B_X$ and $B_Y$ be as in Definition 2.6 with the additional condition $B_Y = B_Y^{\leqslant}$ and consider a Radon measure $\mu \in \mathcal{M}_{B_X \times B_Y}$. Then, for any $\Phi, \Psi \in \mathrm{Cyl}_{B_X \times B_Y}$, we have*

$$\lim_{\varepsilon \to 0} \int (\mathscr{L}_\varepsilon \Phi) \Psi \, \mathrm{d}\mu = \int (\mathscr{L}\Phi) \Psi \, \mathrm{d}\mu. \tag{2.16}$$

**Proof.** We show the proof only for the case $\Psi = 1$, the general case can be deduced via Lebesgue's dominated convergence theorem. Let us prove that $\mathscr{L}\Phi \in L^1(\mu)$, for $\Phi \in \mathrm{Cyl}_{B_X \times B_Y}$. By definition, we have that $\mathrm{tr}_{L^2(\mathbb{R}^2)}(\nabla_X^2 \Phi)$ is bounded, and therefore in $L^1(\mu)$. Moreover,

$$|\langle (-\Delta + m^2)X, \nabla_X \Phi \rangle| = |\langle X, (-\Delta + m^2)\nabla_X \Phi \rangle| \lesssim \|X\|_{B_X} \|\nabla_X \Phi\|_{B_{1,1,-\ell'}^{2-\delta}},$$

and since $\|X\|_{B_X} \in L^1(\mu)$, $P_*^X \mu$ being the Gaussian free field, and $\|\nabla_X \Phi\|_{B_{1,1,-\ell'}^{2-\delta}} \in L^\infty(B_X \times B_Y, \mathbb{R})$, then also the term $\langle (-\Delta + m^2)X, \nabla_X \Phi \rangle$ is in $L^1(\mu)$. Now, the term $\langle (-\Delta + m^2)Y, \nabla_Y \Phi \rangle$ can be handled similarly by using the hypothesis $\|Y\|_{B_Y} \in L^1(\mu)$.

We are left to consider the term $\langle \mathscr{G}, \nabla_Y \Phi \rangle = \langle \alpha : e^{\alpha X}: e^{\alpha Y}, \nabla_Y \Phi \rangle$. First of all, we note that the product $:\exp(\alpha X): \exp(\alpha Y)$ is well-defined. Indeed, since $P_*^X \mu$ is the Gaussian free field we can exploit Proposition B.1 to get $:\exp(\alpha X): \in B_{\exp}^{r,\ell'}$ (cf. Remark 2.7). Furthermore, since $Y \leqslant 0$ and $Y \in B_{p,p,\ell}^{s-\delta}(\mathbb{R}^2)$, by composition of a smooth bounded function with a Besov function we have $\exp(\alpha Y) \in B_{p,p,\ell}^{s-\delta}(\mathbb{R}^2)$. Thus, by Theorem A.3, provided

$$\frac{1}{r'} := \frac{1}{r} + \frac{1}{p} < 1, \qquad s - \frac{\alpha^2}{4\pi}(r-1) - 2\delta > 0,$$

we have that $:\exp(\alpha X): \exp(\alpha Y) \in B_{\exp}^{r',\ell+\ell'}$, and the product is well-defined and continuous. On the other hand, since we have also $\exp(\alpha Y) \in L^\infty(\mathbb{R}^2)$, then Proposition A.4 implies that $:\exp(\alpha X): \exp(\alpha Y) \in B_{\exp}^{r,\ell}$ and

$$\|:e^{\alpha X}: e^{\alpha Y}\|_{B_{\exp}^{r,\ell}} \lesssim \|:e^{\alpha X}:\|_{B_{\exp}^{r,\ell}}.$$

Thus,
$$|\langle \mathcal{G}, \nabla_Y \Phi \rangle| \lesssim \left\| :e^{\alpha X}: \right\|_{B^{r,\ell}_{\exp}} \left\| \nabla_Y \Phi \right\|_{(B^{r,\ell}_{\exp})^*}.$$

By hypothesis we have $\left\| \nabla_Y \Phi \right\|_{(B^{r,\ell}_{\exp})^*} \in L^\infty(B_X \times B_Y, \mathbb{R})$ and, since $r\alpha^2 < 8\pi$ by our assumptions in Definition 2.6, we also have that $\left\| :\exp(\alpha X): \right\|_{B^{r,\ell}_{\exp}} \in L^r(\mu)$. Hence,
$$|\langle \mathcal{G}, \nabla_Y \Phi \rangle| \in L^r(\mu).$$

We are left to show that equation (2.16) holds. This is equivalent to showing
$$\lim_{\varepsilon \to 0} \int (\mathcal{L}_\varepsilon \Phi - \mathcal{L}\Phi) \Psi \, d\mu = \lim_{\varepsilon \to 0} \int \langle \mathcal{G}_\varepsilon - \mathcal{G}, \nabla_Y \Phi \rangle \Psi \, d\mu = 0.$$

By the same reasoning as in the previous part of the proof and by Proposition B.1, we have that $\langle \mathcal{G}_\varepsilon - \mathcal{G}, \nabla_Y \Phi \rangle \in L^r(\mu)$, uniformly with respect to $0 < \varepsilon < 1$. In order to show that the previous limit holds, it is sufficient to prove that
$$\langle \mathcal{G}_\varepsilon - \mathcal{G}, \nabla_Y \Phi \rangle \to 0,$$

in probability as $\varepsilon \to 0$. On the other hand, we know that
$$\left\| \mathcal{G}_\varepsilon - \mathcal{G} \right\|_{B^{r,\ell'}_{\exp}} < +\infty,$$

uniformly with respect to $\varepsilon$, and by Proposition B.1 we have that $\mathcal{G}_\varepsilon \to \mathcal{G}$, weakly in $\mathcal{S}'(\mathbb{R}^2)$. Therefore, we have that $\mathcal{G}_\varepsilon \to \mathcal{G}$ weakly in $B^{r,\ell'}_{\exp}$, and in probability with respect to $\mu$. By definition of weak convergence and since $\nabla_Y \Phi \in (B^{r,\ell'}_{\exp})^*$ by hypothesis (see Definition 3.1), we have that $\langle \mathcal{G}_\varepsilon - \mathcal{G}, \nabla_Y \Phi \rangle \to 0$, in probability. Finally, by property iv. in Definition 3.1 and the fact that
$$\left\| \mathcal{G}_\varepsilon - \mathcal{G} \right\|_{B^{r,\ell}_{\exp}} \in L^r(\mu)$$

uniformly in $\varepsilon > 0$ by Proposition B.1, we have that $\langle \mathcal{G}_\varepsilon - \mathcal{G}, \nabla_Y \Phi \rangle \in L^p(\mu)$ uniformly in $\varepsilon > 0$, for some $p > 1$, and thus the result follows from Lebesgue's dominated convergence theorem. □

**Remark 2.9.** By Theorem 2.8 it is possible to take the limit $\varepsilon \to 0$ inside the integration with respect to $\mu$ in equations (2.4)–(2.3) and (2.5)–(2.6) so that the mentioned equations are equivalent to the one where the limit disappears and the operator $\mathcal{L}_\varepsilon$ is replaced by $\mathcal{L}$. Another consequence of the previous theorem is that, if we consider $\tilde{\mathcal{L}}_\varepsilon$ to be another approximation of $\mathcal{L}$ such that
$$\lim_{\varepsilon \to 0} \int (\tilde{\mathcal{L}}_\varepsilon \Phi) \Psi \, d\mu = \int (\mathcal{L}\Phi) \Psi \, d\mu, \quad \text{for any } \Phi, \Psi \in \mathrm{Cyl}_{B_X \times B_Y}, \tag{2.17}$$

then Problems A″ and B with the operator $\tilde{\mathcal{L}}_\varepsilon$ are equivalent to Problems A″ and B with the operator $\mathcal{L}_\varepsilon$. This means that the formulations given in Problems A″ and B do not depend on the precise form of the approximating operator $\mathcal{L}_\varepsilon$ but only on its limit $\mathcal{L}$.

We now want to prove a result that justifies our restriction in taking $B_Y = B_Y^{\leqslant}$ in the results above and in the rest of the paper. Indeed, if we focus on Problem B, the solutions are always supported on the space of negative functions on the $Y$ component. In the light of point *iii.* in Theorem 2.5, the next result implies that any solution $\mu$ to Problem B belongs to the spaces $\mathscr{M}_{B_X \times B_Y}$.

**Theorem 2.10.** *Let $\mu$ be a solution to Problem B, then we have*

$$\operatorname{supp} P_*^Y(\mu) \subset \{Y \leqslant 0\}.$$

**Proof.** Suppose that $Y \in B_Y = B_{p,p,\ell}^{s-\delta}(\mathbb{R}^2)$ (see Definition 2.6). Define, for $r, k > 0$,

$$\rho_r^k(x) = (1 + k|x|^2)^{-r/2}, \qquad x \in \mathbb{R}^2,$$

and

$$f(x) = x \vee 0, \qquad I(x) = \int_0^x f(y)\, dy, \qquad x \in \mathbb{R}.$$

We consider, for $\eta > 0$, the convolution $g_\eta * Y$ which, by definition of $B_Y$, belongs to $B_{2,2,\ell'}^{s'}$ for any $\ell', s' \geqslant 1$ such that $p\ell' > 2$. Take

$$F(Y) = \arctan(\langle I(g_\eta * Y), \rho_r^k \rangle),$$

so that, by integration by parts,

$$\begin{aligned}
\mathscr{L}_\varepsilon F(Y) &= \frac{1}{1 + \langle I(g_\eta * Y), \rho_r^k \rangle^2} \langle f(g_\eta * Y), \rho_r^k g_\eta * (-(-\Delta + m^2) Y - \mathscr{G}_\varepsilon(X, Y)) \rangle \\
&= -\frac{1}{1 + \langle I(g_\eta * Y), \rho_r^k \rangle^2} \langle f(g_\eta * Y), \rho_r^k((-\Delta + m^2)(g_\eta * Y) + g_\eta * \mathscr{G}_\varepsilon(X, Y)) \rangle \\
&= -\frac{1}{1 + \langle I(g_\eta * Y), \rho_r^k \rangle^2} \big[ \langle f(g_\eta * Y), \rho_r^k m^2 (g_\eta * Y) \rangle + \langle f(g_\eta * Y), \rho_r^k(-\Delta)(g_\eta * Y) \rangle \\
&\qquad\qquad + \langle f(g_\eta * Y), \rho_r^k(g_\eta * \mathscr{G}_\varepsilon(X, Y)) \rangle \big] \\
&= -\frac{1}{1 + \langle I(g_\eta * Y), \rho_r^k \rangle^2} \big[ \langle f(g_\eta * Y), \rho_r^k m^2 (g_\eta * Y) \rangle + \langle f'(g_\eta * Y)(g_\eta * \nabla Y), \rho_r^k(g_\eta * \nabla Y) \rangle \\
&\qquad\qquad + \langle f(g_\eta * Y), \nabla \rho_r^k(g_\eta * \nabla Y) \rangle + \langle f(g_\eta * Y), \rho_r^k(g_\eta * \mathscr{G}_\varepsilon(X, Y)) \rangle \big].
\end{aligned}$$

Consider the term $\langle f(g_\eta * Y), \nabla \rho_r^k(g_\eta * \nabla Y) \rangle$, and multiply and divide by $\rho_r^k$ to get

$$\langle f(g_\eta * Y), \nabla \rho_r^k(g_\eta * \nabla Y) \rangle = \langle \rho_r^k f(g_\eta * Y), \frac{\nabla \rho_r^k}{\rho_r^k}(g_\eta * \nabla Y) \rangle.$$

By Young inequality, we have

$$\begin{aligned}
\mathscr{L}_\varepsilon F(Y) &\leqslant -\frac{1}{1 + \langle I(g_\eta * Y), \rho_\ell^k \rangle^2} \Big[ \langle f(g_\eta * Y), \rho_r^k m^2(g_\eta * Y) \rangle + \langle f'(g_\eta * Y)(g_\eta * \nabla Y), \rho_r^k(g_\eta * \nabla Y) \rangle \\
&\quad - \frac{1}{2} \langle \rho_r^k f^2(g_\eta * Y), \left(\frac{\nabla \rho_r^k}{\rho_r^k}\right)^2 \rangle - \frac{1}{2} \langle \rho_r^k, |g_\eta * \nabla Y|^2 \mathbb{I}_{\{g_\eta * Y \geqslant 0\}} \rangle \Big],
\end{aligned}$$

where we also used that the term $\langle f(g_\eta * Y), \rho_r^k(g_\eta * \mathcal{G}_\varepsilon(X, Y))\rangle$ is positive. We want to show that $\nabla \rho_r^k / \rho_r^k$ is bounded. Taking $y = \sqrt{k} x$, we have, for $k > 0$ small enough,

$$\frac{\nabla \rho_r^k}{\rho_r^k} = \frac{2kx_1}{1 + k|x|^2} = \frac{2\sqrt{k}\, y_1}{1 + |y|^2} \leq 2\sqrt{k} \sup_{y \in \mathbb{R}^2} \frac{y_1}{1 + |y|^2} \leq C_k,$$

which gives, if $C_k/2 - m^2 < -\zeta$,

$$\mathcal{L}_\varepsilon F(Y) < -\frac{1}{1 + \langle I(g_\eta * Y), \rho_\ell^k\rangle^2} \left[ \zeta \langle f^2(g_\eta * Y), \rho_r^k\rangle + \frac{1}{2}\langle \rho_r^k, |g_\eta * \nabla Y|^2 \mathbb{I}_{\{g_\eta * Y \geq 0\}}\rangle \right],$$

since $f(g_\eta * Y) = g_\eta * Y$ on $\mathrm{supp}(f)$.

Since the right-hand side does not depend on $\varepsilon$, and the left-hand side converges to 0 when integrated with respect to $\mu$, we have

$$\int \frac{\zeta}{1 + \langle I(g_\eta * Y), \rho_\ell^k\rangle^2} \left[ \langle f^2(g_\eta * Y), \rho_r^k\rangle + \frac{1}{2}\langle \rho_r^k, |g_\eta * \nabla Y|^2 \mathbb{I}_{\{g_\eta * Y \geq 0\}}\rangle \right] \mathrm{d}\mu \leq 0,$$

and since the terms $2^{-1}\langle \rho_r^k, |g_\eta * \nabla Y|^2 \mathbb{I}_{\{g_\eta * Y \geq 0\}}\rangle$, $\rho_r^k$, $|g_\eta * \nabla Y|^2$, and $\zeta(1 + \langle I(g_\eta * Y), \rho_\ell^k\rangle^2)^{-1}$ are all positive, we have

$$\int \zeta \langle f^2(g_\eta * Y), \rho_r^k\rangle \mathrm{d}\mu \leq 0.$$

By Fatou's Lemma we can take the limit as $\eta \to 0$ and get that

$$\mathbb{I}_{\{Y \geq 0\}} = 0, \qquad \mu\text{-a.s.} \qquad \square$$

**Remark 2.11.** A consequence of Theorem 2.10 is that for the case of exponential interaction, a solution to Problem B with $B_Y = B_{p,p,\ell}^{s-\delta}(\mathbb{R}^2)$ is necessarily a solution to Problem A'' with $B_Y = B_Y^\leq$.

## 2.3 A description of the exponential measure

What has been described in the previous section has the main aim of giving a characterization to a unique measure which we call the exponential QFT. In this section, we want to connect our discussion of the IbP formula with the standard approach to construct EQFT via Gibbsian modifications of the GFF. Assume that $\mu_M^{\mathrm{free}}$ is the measure of the free field with mass $m$ on the torus $\mathbb{T}_M^2$ of size $M$. We define the interaction, for $\varepsilon > 0$,

$$V_{M,\varepsilon}^{\exp,\alpha}(\varphi) := \int_{\mathbb{T}_M^2} f_\varepsilon(x) e^{\alpha(g_\varepsilon * \varphi) - \frac{\alpha^2}{2} c_\varepsilon} \mathrm{d}x,$$

with $f$, $g$, $c_\varepsilon$ and $\alpha$ are as in Section 2.2. We also consider

$$Z_{M,\varepsilon} := \int e^{-V_{M,\varepsilon}^{\exp,\alpha}(\varphi)} \mu_M^{\mathrm{free}}(\mathrm{d}\varphi).$$

**Definition 2.12.** *We say that the measure $v_m^{\exp,\alpha}$ on $\mathcal{S}'(\mathbb{R}^2)$ is the measure related with the Euclidean QFT having action*

$$S(\varphi) = \frac{1}{2}\int_{\mathbb{R}^2}(|\nabla\varphi(x)|^2 + m^2\varphi(x)^2)\,\mathrm{d}x + \int_{\mathbb{R}^2}:e^{\alpha\varphi(x)}:\mathrm{d}x,$$

*if there are two sequences $\varepsilon_n \to 0$ and $M_{n'} \to +\infty$, with $\varepsilon_n > 0$ and $N_n \in \mathbb{N}$ such that*

$$v_m^{\exp,\alpha}(\mathrm{d}\varphi) = \lim_{n'\to+\infty}\lim_{n\to+\infty} Z_{M_{n'},\varepsilon_n}^{-1} e^{-V_{M_{n'},\varepsilon_n}^{\exp,\alpha}(\varphi)}\mu_{M_{n'}}^{\mathrm{free}}(\mathrm{d}\varphi), \tag{2.18}$$

*where the limits are taken in weak sense in the space of probability measures.*

Such a measure $v_m^{\exp,\alpha}$ was first built by Albeverio and Høegh-Krohn in [AH74] (see also [FP77]) using techniques from constructive QFT. More recently, this model was studied in the context of stochastic quantization on the torus or on a compact manifold (see e.g. [AKMR21, BG21a, Gar20, HKK21, HKK22, ORTW20, ORW21]), and on $\mathbb{R}^2$ (see [ADG21]). See also [BD21] for the related $\cosh(\Phi)_2$-model.

A consequence of the results we presented in the previous sections of the paper (see Theorem 1.2, Theorem 1.3) is the following differential characterization of the measure related to the exponential interaction.

**Theorem 2.13.** *The following statements hold:*

   i. *For any $\alpha^2 < 8\pi$, there exists a measure related to the exponential interaction defined by the limit (2.18).*

   ii. *Let $B_X = C_\ell^{-\delta}(\mathbb{R}^2)$ and $B_Y = B_Y^{\leqslant}$. Then, for any $\alpha^2 < 4\pi\tilde{\gamma}_{\max}$, the measure $v_m^{\exp,\alpha}$ is the unique measure in the space $\mathcal{M}_{B_X+B_Y}$ such that, for any $F \in \mathrm{Cyl}_E^b$ and $h \in \mathcal{S}(\mathbb{R}^2)$, we have*

$$\int \langle \nabla_\varphi F, h\rangle\, v(\mathrm{d}\varphi) = \lim_{\varepsilon\to 0}\int F(\varphi)\,\langle(-\Delta+m^2)\varphi + \alpha f_\varepsilon e^{\alpha g_\varepsilon*\varphi-\frac{\alpha^2}{2}c_\varepsilon}, h\rangle\, v(\mathrm{d}\varphi).$$

# 3 Uniqueness of solution

## 3.1 Proof of uniqueness of solutions to Problem A″

In this section, we discuss the uniqueness of solutions to Problems A, A″, and B in the case where $B_X$, $B_Y$, $\mathcal{M} = \mathcal{M}_{B_Y}$, and $E = B_X + B_Y$ are as in Section 2.2 and $B_\varepsilon$ is the drift of the exponential interaction (see equation (1.3)). The method that we adopt to prove uniqueness is the study of the resolvent equation for the operator $\mathcal{L}_\varepsilon$ defined in equation (2.14), namely

$$(\lambda - \mathcal{L}_\varepsilon)G_\varepsilon^\lambda = F, \qquad \lambda \in \mathbb{R}_+,. \tag{3.1}$$

We are interested in classical solutions to the previous equation, namely solutions $G_\varepsilon^\lambda$ that are at least in $C^2(B_X \times B_Y)$ and to which the operator $\mathcal{L}_\varepsilon$ can be applied.

Let us first introduce here a class of functions where the solutions of the resolvent equation belongs.

**Definition 3.1.** *Recall that $\gamma = \alpha^2/(4\pi)$. Let $p, s, \ell, r$ be the same parameters as in Definition 2.6 and let us introduce the space*

$$\hat{B}_X = B_{\tilde{p},\tilde{p},\ell}^{-\delta}(\mathbb{R}^2),$$

*for some $1 < \tilde{p} < +\infty$ large enough. Note that $B_X \subset \hat{B}_X$. Denote by $\mathscr{F}$ the class of bounded, measurable functions $\Phi : \hat{B}_X \times B_Y \to \mathbb{R}$ such that*

  i. $\nabla_Y \Phi \in C^0(\hat{B}_X \times B_Y, B_{l,l,-\ell}^{(2-s)\wedge(\gamma(r-1))+\delta}(\mathbb{R}^2))$, *for any $\delta > 0$, $1 < l < +\infty$.*

  ii. $\nabla_X \Phi \in C^0(\hat{B}_X \times B_Y, B_{1,1,-\ell}^{2-\delta}(\mathbb{R}^2))$, *so that $\langle \nabla_X \Phi, (-\Delta + m^2)X \rangle$ is well-defined.*

  iii. *The operator $\nabla_X^2 \Phi \in C^0(\hat{B}_X \times B_Y, L(\hat{B}_X, B_X))$ can be extended in a unique continuous way to an operator $\nabla_X^2 \Phi \in C^0(\hat{B}_X \times B_Y, L(H_\ell^{-\kappa}(\mathbb{R}^2), H_{-\ell}^{\kappa}(\mathbb{R}^2)))$, where $\kappa > 1$ and $\ell > 1$.*

  iv. *There exists some $f_\Phi(X) \in L^p(\mu^{\text{free}})$, $p \in [1, +\infty)$, such that*

$$\left\| \nabla_Y \Phi(X,Y) \right\|_{B_{l,l,-\ell}^{(2-s)\wedge(\gamma(r-1))+\delta}(\mathbb{R}^2)} \leq f_\Phi(X),$$
$$\left\| \nabla_X \Phi(X,Y) \right\|_{B_{1,1,-\ell}^{2-\delta}(\mathbb{R}^2)} \leq f_\Phi(X),$$
$$\left\| \nabla_X^2 \Phi(X,Y) \right\|_{L(H_\ell^{-\kappa}(\mathbb{R}^2), H_{-\ell}^{\kappa}(\mathbb{R}^2))} \leq f_\Phi(X).$$

*We define the space $\mathfrak{F}$ as the space of functions $F: \tilde{B}_X + B_Y \to \mathbb{R}$ such that*

$$F \circ P^{X+Y} \in \mathscr{F}.$$

**Remark 3.2.** *Under condition iii. of the previous definition, since the immersion $H_\ell^{-\kappa}(\mathbb{R}^2) \hookrightarrow L^2(\mathbb{R}^2)$ is an Hilbert-Schmidt operator, then $\nabla_X^2 \Phi, \rho_{-\ell} \nabla_X^2 \Phi \in L_{\text{loc}}^\infty(\hat{B}_X \times B_Y, \text{TC}(L^2))$, where $\text{TC}(L^2)$ is the space of trace-class operators on $L^2(\mathbb{R}^2)$.*

**Remark 3.3.** *The classes $\mathscr{F}$ and $\mathfrak{F}$ are chosen in such a way that they satisfy the following two important properties: (i) $\mathscr{F}$ and $\mathfrak{F}$ contain $\text{Cyl}_{B_X \times B_Y}$ and $\text{Cyl}_E$, respectively, and (ii) if $B_Y = B_Y^{\leq}$ then, for any Radon measure $\mu \in \mathcal{M}_{B_X \times B_Y}$, we have $\sup_{\varepsilon > 0} \int |\mathcal{L}_\varepsilon \Phi|^\sigma d\mu < +\infty$, for any $\Phi \in \mathscr{F}$ and $\sigma \geq 1$.*

The classes $\mathscr{F}$ and $\mathfrak{F}$ satisfy the following important lemma.

**Lemma 3.4.** *Suppose that $\mu \in \mathcal{M}_{B_X \times B_Y}$ and satisfies equation (2.5) for any $\Phi \in \text{Cyl}_{B_X \times B_Y}$, then it satisfies the same equation for every $\Phi \in \mathscr{F}$. Suppose that $\nu \in \mathcal{M}_{B_Y}$ and satisfies equation (2.1) for any $\Phi \in \text{Cyl}_E$, then it satisfies the same equation for every $\Phi \in \mathfrak{F}$.*

**Proof.** See Appendix D. □

**Proposition 3.5.** *Let $F \in \text{Cyl}_{B_X \times B_Y}$, then there exists a classical solution $G_\varepsilon^\lambda \in C^2(B_X \times B_Y)$ to the resolvent equation (3.1) with the following properties:*

  i. $G_\varepsilon^\lambda \in \mathscr{F}$ *and moreover, if $F = \bar{F} \circ P^{X+Y}$ for some $\bar{F} \in \text{Cyl}_E$, then there exists $\bar{G}_\varepsilon^\lambda \in \mathfrak{F}$ such that we have $G_\varepsilon^\lambda = \bar{G}_\varepsilon^\lambda \circ P^{X+Y}$.*

    *ii.* Suppose that $F$ has compact support in Fourier variables, then there exists $\varepsilon_0 > 0$ such that, for every $\mu_1, \mu_2 \in \mathcal{M}_{B_X \times B_Y}$ and for every $\varsigma \in (0, 1)$, there are two constants $C_{\mu_1, \mu_2, \varsigma} > 0$ and $K > 0$ such that

$$\lambda \left| \int G_\varepsilon^\lambda \, (\mathrm{d}\mu_1 - \mathrm{d}\mu_2) \right| \lesssim \varsigma + \frac{\lambda}{\lambda + K} C_{\mu_1, \mu_2, \varsigma},$$

    where the constant included in the symbol $\lesssim$ does not depend on $\lambda$, $\mu_1$, $\mu_2$, $\varepsilon$ or $\varsigma$.

    *iii.* If $\alpha^2 < 4\pi \tilde{\gamma}_{\max}$, then there exists $q > 1$ such that, for every measurable $\mathcal{K} : B_X \times B_Y \to B_{\exp}^{q, \ell/2}$ and every $\mu \in \mathcal{M}_{B_X \times B_Y}$, we have

$$\int |\langle \nabla_{Y_0} G_\varepsilon^\lambda, \mathcal{K} \rangle | \, \mathrm{d}\mu \lesssim_\lambda \left( \int \|\mathcal{K}\|_{B_{\exp}^{q, \ell/2}}^q \, \mathrm{d}\mu \right)^{1/q},$$

uniformly in $\varepsilon > 0$.

Let $\tilde{\gamma}_{\max}$ be the parameter defined in Remark 1.1. In the light of Proposition 3.5, the proof of which we postpone to Section 3.2, we can then proceed with the proof of uniqueness of solutions to Problem A''.

**Theorem 3.6.** *Let $\alpha^2 < 2\pi \tilde{\gamma}_{\max}$ and take $B_X = C_\ell^{-\delta}(\mathbb{R}^2)$ and $B_Y = B_Y^{\leqslant}$, where the parameters are taken as in Definition 2.6 and the space $B_Y^{\leqslant}$ is defined as in equation (2.11). Then the solution to FPK equation in the sense of Problem A'' is unique.*

**Proof.** Let $\mu_1$ and $\mu_2$ be two solutions of Problem A''. We want to show that, for any $F \in \mathrm{Cyl}_{B_X \times B_Y}$ with compact support in Fourier variables, $\int F \mathrm{d}\mu_1 = \int F \mathrm{d}\mu_2$. This implies $\mu_1 = \mu_2$ since $F \in \mathrm{Cyl}_{B_X \times B_Y}$ with compact support in Fourier variables separates points of the space of Radon measures on $B_X \times B_Y$. Consider the solution $G_\varepsilon^\lambda$ of the resolvent equation (3.1) given by Proposition 3.5 and recall that, by point *ii.* in Proposition 3.5. Since $\mu_1$ and $\mu_2$ are solutions to the FPK equations we have in particular that equation (2.5) holds. Therefore, integrating with respect to $\mu_1$ (respectively, with respect to $\mu_2$) the resolvent equation equation (3.1) and subtracting the integral $\int \mathscr{L} G_\varepsilon^\lambda \, \mathrm{d}\mu_1 = 0$ (respectively, $\int \mathscr{L} G_\varepsilon^\lambda \, \mathrm{d}\mu_2 = 0$) which holds by Lemma 3.4, we get the relation

$$\int F \mathrm{d}\mu_j = \lambda \int G_\varepsilon^\lambda \, \mathrm{d}\mu_j - \int \nabla_{Y_0} G_\varepsilon^\lambda (\mathscr{G} - \mathscr{G}_\varepsilon) \, \mathrm{d}\mu_j, \qquad j = 1, 2,$$

where we omitted the dependences on $(X_0, Y_0)$ for the sake of brevity. Taking the difference of such a relation for $j = 1, 2$ yields

$$\int F \mathrm{d}\mu_1 - \int F \mathrm{d}\mu_2 = \lambda \int G_\varepsilon^\lambda \, \mathrm{d}\mu_1 - \lambda \int G_\varepsilon^\lambda \, \mathrm{d}\mu_2 - \int \nabla_{Y_0} G_\varepsilon^\lambda (\mathscr{G} - \mathscr{G}_\varepsilon) \, \mathrm{d}\mu_1 + \int \nabla_{Y_0} G_\varepsilon^\lambda (\mathscr{G} - \mathscr{G}_\varepsilon) \, \mathrm{d}\mu_2.$$

Now, point *ii.* in Proposition 3.5 gives us, for any $\varsigma > 0$,

$$\left| \int F \mathrm{d}\mu_1 - \int F \mathrm{d}\mu_2 \right| \lesssim \varsigma + (1 + C_\varsigma) \frac{\lambda}{\lambda + K} + \left| \int \nabla_{Y_0} G_\varepsilon^\lambda (\mathscr{G} - \mathscr{G}_\varepsilon) \, \mathrm{d}\mu_1 \right| + \left| \int \nabla_{Y_0} G_\varepsilon^\lambda (\mathscr{G} - \mathscr{G}_\varepsilon) \, \mathrm{d}\mu_2 \right|$$

On the other hand, point *iii.* in Proposition 3.5 implies that

$$\left| \int \nabla_{Y_0} G_\varepsilon^\lambda (\mathscr{G} - \mathscr{G}_\varepsilon) \, \mathrm{d}\mu_j \right| \lesssim_\lambda \left( \int \|\mathscr{G} - \mathscr{G}_\varepsilon\|_{B_{\exp}^{q, \ell/2}}^q \, \mathrm{d}\mu_j \right)^{1/q}.$$

By the proof of Theorem 2.8, Proposition B.1, and using the fact that $P_*^X \mu_j = \mu^{\text{free}}$, we get that

$$\left( \int \| \mathscr{G} - \mathscr{G}_\varepsilon \|_{B_{\exp}^{q,\ell/2}}^q \, \mathrm{d}\mu_j \right)^{1/q} \to 0 \qquad \text{as } \varepsilon \to 0.$$

Thus, if we take the limit $\varepsilon \to 0$ and then $\lambda \to 0$, we get $|\int F \mathrm{d}\mu_1 - \int F \mathrm{d}\mu_2| \lesssim \varsigma$, and since $\varsigma > 0$ is arbitrary, we get $\mu_1 = \mu_2$. □

In the case of Problem B a better result, in the sense that $\gamma_{\max} > \tilde{\gamma}_{\max}$ (cf. Remark 1.1), can be achieved exploiting the properties of the resolvent operator proved in Section 3.3.

**Theorem 3.7.** *Let $\alpha^2 < 2\pi\gamma_{\max}$, $B_X = C_\ell^{-\delta}(\mathbb{R}^2)$ and $B_Y = B_{p,p,\ell}^{s-\delta}(\mathbb{R}^2)$, where the parameters are taken as in Definition 2.6. Then the solution to FPK equation in the sense of Problem B is unique.*

**Proof.** The proof is the same as in the case of Problem A'', but substituting point *iv.* of Proposition 3.5, with Proposition 3.11 below. □

## 3.2 Analysis of the resolvent equation

To solve the resolvent equation as needed in Proposition 3.5, we use a probabilistic representation and look for solutions in the form

$$G_\varepsilon^\lambda(X_0, Y_0) = \mathbb{E}_\xi \left[ \int_0^{+\infty} e^{-\lambda t} F(X_t^\varepsilon, Y_t^\varepsilon) \, \mathrm{d}t \right], \tag{3.2}$$

where $X_t^\varepsilon$ and $Y_t^\varepsilon$ are solutions to the stochastic differential system

$$\partial_t X_t^\varepsilon = -(-\Delta + m^2) X_t^\varepsilon + \xi_t, \tag{3.3}$$
$$\partial_t Y_t^\varepsilon = -(-\Delta + m^2) Y_t^\varepsilon - \mathscr{G}_\varepsilon(X_t^\varepsilon, Y_t^\varepsilon), \tag{3.4}$$

where $\mathscr{G}_\varepsilon$ is defined in equation (2.12), $\xi$ is a Gaussian space-time white noise, and with initial conditions $(X_0, Y_0) \in B_X \times B_Y$. This is a triangular system, where the first equation does not depend on the second variable $Y$. The first equation is a linear SDE with additive noise with stationary solution given by the Ornstein-Uhlenbeck process, whose invariant measure is the Gaussian free field with mass $m > 0$. To prove that $G_\varepsilon^\lambda$ is a classical solution we need some properties on the flow induced by the SPDEs (3.3)–(3.4). We will only state such properties in the present section – postponing their proof to Appendix E – and then proceed by showing that Proposition 3.5 holds.

First, let us write down the equations for the derivatives of the flow. Let us denote by $X$ and $Y$ the solutions to equations (3.3)–(3.4), dropping the dependence on the parameter $\varepsilon > 0$ for simplicity of notation. The derivatives of $X$ solve

$$(\partial_t - \Delta + m^2)\nabla_{X_0} X_t = 0, \qquad \nabla_{X_0} X(0) = \mathrm{id}, \tag{3.5}$$
$$(\partial_t - \Delta + m^2)\nabla_{X_0}^2 X_t = 0, \qquad \nabla_{X_0}^2 X(0) = 0, \tag{3.6}$$
$$(\partial_t - \Delta + m^2)\nabla_{Y_0} X_t = 0, \qquad \nabla_{Y_0} X(0) = 0. \tag{3.7}$$

Whenever $Y_0 \in B_Y^{\leqslant}$,

$$(\partial_t - \Delta + m^2)\nabla_{Y_0}Y_t(Y_0)[h] = -D_Y\mathscr{G}_\varepsilon(X_t, Y_t)[\nabla_{Y_0}Y_t(Y_0)[h]], \qquad \nabla_{Y_0}Y_0 = h, \tag{3.8}$$

for $h \in B_Y \cup B_{\exp}^{r,\ell}$, and

$$\begin{aligned}(\partial_t - \Delta + m^2)\nabla_{X_0}Y_t(Y_0)[h] &= -D_Y\mathscr{G}_\varepsilon(X_t, Y_t)[\nabla_{X_0}Y_t(Y_0)[h]] - D_X\mathscr{G}_\varepsilon(X_t, Y_t)[\nabla_{X_0}X_t[h]],\\ \nabla_{X_0}Y_0 &= 0,\end{aligned} \tag{3.9}$$

for $h \in B_X$, and

$$\begin{aligned}(\partial_t - \Delta + m^2)\nabla_{X_0}^2 Y_t(Y_0)[h, h'] &= -D_Y\mathscr{G}_\varepsilon(X_t, Y_t)[\nabla_{X_0}^2 Y_t(Y_0)[h, h']] \\ &\quad -D_Y^2\mathscr{G}_\varepsilon(X_t, Y_t)[\nabla_{X_0}Y_t(Y_0)[h], \nabla_{X_0}Y_t(Y_0)[h']] \\ &\quad -2D_{X,Y}^2\mathscr{G}_\varepsilon(X_t, Y_t)[\nabla_{X_0}X_t[h], \nabla_{X_0}Y_t(Y_0)[h']] \\ &\quad -D_X^2\mathscr{G}_\varepsilon(X_t, Y_t)[\nabla_{X_0}X_t(Y_0)[h], \nabla_{X_0}X_t(Y_0)[h']], \\ \nabla_{X_0}^2 Y_0 &= 0,\end{aligned} \tag{3.10}$$

for $h, h' \in L^2(\mathbb{R}^2)$. The derivatives $D\mathscr{G}_{\varepsilon,\bar{\varepsilon}}$ have the following expressions

$$D_X\mathscr{G}_\varepsilon(X_t, Y_t)[\varphi] = D_Y\mathscr{G}_\varepsilon(X_t, Y_t)[\varphi] = \alpha^2 f_\varepsilon(:e^{\alpha(X_t * g_\varepsilon)}: e^{\alpha(Y_t * g_\varepsilon)})(g_\varepsilon * \varphi),$$

and

$$\begin{aligned}D_X^2\mathscr{G}_\varepsilon(X_t, Y_t)[\varphi, \psi] = D_Y^2\mathscr{G}_\varepsilon(X_t, Y_t)[\varphi, \psi] &= D_{X,Y}^2\mathscr{G}_\varepsilon(X_t, Y_t)[\varphi, \psi] \\ &= \alpha^3 f_\varepsilon(:e^{\alpha(X_t * g_\varepsilon)}: e^{\alpha(Y_t * g_\varepsilon)})(g_\varepsilon * \varphi)(g_\varepsilon * \psi).\end{aligned}$$

Hereafter, we will use the notation $\gamma = \alpha^2/(4\pi)$, and

$$\tilde{B}_Y = B_Y \cup B_{p \wedge r, p \wedge r, \ell}^{-(2-s) \wedge (\gamma(r-1)) - \delta}(\mathbb{R}^2), \qquad \tilde{B}_X = B_{\infty, \infty, \ell}^{-2+\delta}(\mathbb{R}^2), \tag{3.11}$$

where all the appearing parameters are defined as in Definition 2.6.

**Proposition 3.8.** *For any $\varepsilon > 0$, if $(X_0, Y_0) \in \hat{B}_X \times \{B_Y \cup B_{\exp}^{r,\ell}\}$, then there exists a unique solution $(X, Y)$ to equations (3.3)–(3.4) such that*

$$(X_t, Y_t) \in \hat{B}_X \times \{B_Y \cup B_{\exp}^{r,\ell}\}, \qquad t \in \mathbb{R}_+.$$

*Let $X$ and $Y$ be solutions to equations (3.3)–(3.4).*

i. *For every $\varepsilon > 0$, we have that that the derivatives $\nabla_{X_0}X_t$, $\nabla_{X_0}^2 X_t$, $\nabla_{Y_0}X_t$, $\nabla_{X_0}Y$, $\nabla_{X_0}^2 Y$, and $\nabla_{Y_0}Y$ of $X$ and $Y$ exist and satisfy equations (3.5), (3.6), (3.7), (3.8), (3.9), and (3.10), respectively. Furthermore, they are all continuous functions with respect to $X_0$ and $Y_0$.*

ii. *For every $\delta \in (0, 1)$, $\theta \in (0, 1 - \delta)$, $\ell, \ell' \geqslant 1$, we have the estimates*

$$\begin{aligned}&\|\nabla_{Y_0}Y_t(Y_0)[h]\|_{C_{\ell'}^\theta(\mathbb{R}_+, C_{-\ell}^{2-2\theta-2\delta}(\mathbb{R}^2)) \oplus L^\infty(\mathbb{R}_+, \tilde{B}_Y)} \\ &\lesssim_{g_\varepsilon} \mathfrak{P}_2(\|f_\varepsilon : e^{\alpha(g_\varepsilon * X_t)} : e^{\alpha P_t(g_\varepsilon * h)}\|_{L_{t'}^\infty(\mathbb{R}_+, L_{-\ell}^\infty(\mathbb{R}^2))}, \|h\|_{\tilde{B}_Y}), \quad h \in \tilde{B}_Y,\end{aligned} \tag{3.12}$$

and

$$\|\nabla_{X_0} Y_t(Y_0)[h]\|_{C_\ell^\beta(\mathbb{R}_+, C_{-\ell}^{2-2\theta-2\delta}(\mathbb{R}^2))}$$
$$\lesssim_{g_\varepsilon} \tilde{\mathfrak{P}}_2(\|f_\varepsilon : e^{\alpha(g_\varepsilon * X_t)} : e^{\alpha P_t(g_\varepsilon * h)}\|_{L_\ell^\infty(\mathbb{R}_+, L_{-\ell}^\infty(\mathbb{R}^2))}), \quad h \in \tilde{B}_X, \quad (3.13)$$

where $\mathfrak{P}_2$ and $\tilde{\mathfrak{P}}_2$ are two second degree polynomials.

iii. For every $\ell, \kappa \geq 0$, we have that there exist $\beta, \delta > 0$ such that

$$\|\nabla_{X_0}^2 Y_t(Y_0)\|_{L(H_\ell^{-\kappa}, H_{-\ell}^\kappa)} \lesssim \left( \int_{\mathbb{R}^2} \alpha^2 f_\varepsilon(z') : e^{\alpha(g_\varepsilon * X_t)(z')} : e^{\alpha(g_\varepsilon * Y_t)(z')} e^{(\delta - m^2)t} (1 + |z'|^\beta) \, dz' \right)^2. \quad (3.14)$$

The remainder of the present subsection is devoted to the proof of Proposition 3.5.

**Proof of Proposition** 3.5. Let us start by proving that $G_\varepsilon^\lambda$ is a classical solution of the resolvent equation. We exploit Itô's formula appearing in Theorem 4.32 in [DZ14]. Notice that it can be applied to $G_\varepsilon^\lambda$ since, by point i. in Proposition 3.8, $G_\varepsilon^\lambda$ is a $C^2$-function with trace-class second derivative (continuity of the second derivative can be deduced with similar techniques as the ones adopted in Proposition E.7 and point i. in Proposition 3.5).

We use the notation $X_t^{X_0}$ to denote a process $X$ at time $t$, with starting point $X_0$ at time $t = 0$. Recalling the definition of $G_\varepsilon^\lambda$ in equation (3.2), we have

$$G_\varepsilon^\lambda(X_t^{X_0}, Y_t^{Y_0}) = \int_0^{+\infty} e^{-\lambda s} \mathbb{E}_\xi[F(X_s^{X_t}, Y_s^{Y_t}) | \mathcal{H}_t] \, ds,$$

where $(\mathcal{H}_t)_{t>0}$ is the filtration generated by the white noise at time $t$ and the initial conditions $X_0$ and $Y_0$. By Markovianity of the process, we get

$$G_\varepsilon^\lambda(X_t^{X_0}, Y_t^{Y_0}) = \int_0^{+\infty} e^{-\lambda s} \mathbb{E}_\xi[F(X_{t+s}^{X_0}, Y_{t+s}^{Y_0}) | \mathcal{H}_t] \, ds.$$

Since from now on $X$ and $Y$ always start at $X_0$ and $Y_0$, respectively, when $t = 0$, we drop the dependence on the initial conditions. We have, by Itô's formula,

$$e^{-\lambda t} G_\varepsilon^\lambda(X_t, Y_t) - G_\varepsilon^\lambda(X_0, Y_0)$$
$$= \int_0^t e^{-\lambda s} \mathcal{L}_\varepsilon(G_\varepsilon^\lambda(X_s, Y_s)) \, ds + \int_0^t e^{-\lambda s} \nabla_{X_0} G_\varepsilon^\lambda(X_s, Y_s) \cdot dX_s - \lambda \int_0^t e^{-\lambda s} G_\varepsilon^\lambda(X_s, Y_s) \, ds.$$

On the other side of the equation, we have

$$\int_0^{+\infty} e^{-\lambda(t+s)} \mathbb{E}_\xi[F(X_{t+s}, Y_{t+s}) | \mathcal{H}_t] \, ds - \int_0^{+\infty} e^{-\lambda s} \mathbb{E}_\xi[F(X_s, Y_s) | \mathcal{H}_t] \, ds.$$

Notice that

$$\mathbb{E} \int_0^{+\infty} e^{-\lambda(t+s)} \mathbb{E}[F(X_{t+s}, Y_{t+s}) | \mathcal{H}_t] \, ds = \mathbb{E} \int_t^{+\infty} e^{-\lambda s} \mathbb{E}[F(X_s, Y_s) | \mathcal{H}_t] \, ds,$$

and therefore

$$\mathbb{E} \int_t^{+\infty} e^{-\lambda s} \mathbb{E}[F(X_s, Y_s) | \mathcal{H}_t] \, ds - \mathbb{E} \int_0^{+\infty} e^{-\lambda s} \mathbb{E}[F(X_s, Y_s) | \mathcal{H}_t] \, ds = -\mathbb{E} \int_0^t e^{-\lambda s} \mathbb{E}[F(X_s, Y_s) | \mathcal{H}_t] \, ds.$$

Dividing by $t$ and letting $t \to 0$, we get

$$-\frac{1}{t}\mathbb{E}\int_0^t e^{-\lambda s}\mathbb{E}[F(X_s, Y_s)|\mathscr{H}_t]\,\mathrm{d}s \to -F(X_0, Y_0).$$

Taking expectation also on the first side, we have the desired relation.

**Proof of point i.** When $F \in \mathrm{Cyl}^b_{B_X \times B_Y}$, the derivative of $G^\lambda_\varepsilon(X_0, Y_0)$ with respect to $Y_0$ is given by

$$\nabla_{Y_0} G^\lambda_\varepsilon(X_0, Y_0)(h) = \mathbb{E}_\xi\left[\int_0^{+\infty} e^{-\lambda t}\left(\sum_{k=1}^N \partial_k \tilde{F}\langle \nabla_{Y_0} Y_t(h), v_k\rangle\right)\mathrm{d}t\right].$$

Therefore, we have

$$|\nabla_{Y_0} G^\lambda_\varepsilon(X_0, Y_0)(h)| \lesssim \sup_{k=1,\ldots,N} \mathbb{E}_\xi\left[\int_0^{+\infty} e^{-\lambda t}\|\partial_k \tilde{F}\|_{L^\infty}|\langle \nabla_{Y_0} Y_t(h), v_k\rangle|\,\mathrm{d}t\right].$$

Applying the estimate (3.12), we get

$$\begin{aligned}|\langle \nabla_{Y_0} Y_t(h), v_k\rangle| &\lesssim \rho_{\ell'}(t)\|\nabla_{Y_0} Y_t(Y_0)(h)\|_{C^\theta_\ell(\mathbb{R}_+, C^{2-2\theta-2\delta}_{-\ell}(\mathbb{R}^2))\oplus L^\infty(\mathbb{R}_+, \tilde{B}_Y)}\|v_k\|_{\tilde{B}^*_Y}\\ &\lesssim \rho_{\ell'}(t)\, P_2(\|f_\varepsilon : e^{\alpha(g_\varepsilon * X_t)} : e^{\alpha P_t(g_\varepsilon * h)}\|_{L^\infty_{\ell'}(\mathbb{R}_+, L^\infty_{-\ell}(\mathbb{R}^2))}, \|h\|_{\tilde{B}_Y})\|v_k\|_{\tilde{B}^*_Y}.\end{aligned}$$

Therefore,

$$\begin{aligned}|\nabla_{Y_0} G^\lambda_\varepsilon(X_0, Y_0)(h)| &\lesssim \sup_k (\|\partial_k \tilde{F}\|_{L^\infty}\|v_k\|_{\tilde{B}^*_Y})\int_0^{+\infty} e^{-\lambda t}\rho_{\ell'}(t)\\ &\quad\times\mathbb{E}_\xi[P_2(\|f_\varepsilon : e^{\alpha(g_\varepsilon * X_t)} : e^{\alpha P_t(g_\varepsilon * h)}\|_{L^\infty_{\ell'}(\mathbb{R}_+, L^\infty_{-\ell}(\mathbb{R}^2))}, \|h\|_{\tilde{B}_Y})]\,\mathrm{d}t.\end{aligned}$$

We have that

$$\begin{aligned}&\mathbb{E}_\xi[P_2(\|f_\varepsilon : e^{\alpha(g_\varepsilon * X_t)} : e^{\alpha P_t(g_\varepsilon * h)}\|_{L^\infty_{\ell'}(\mathbb{R}_+, L^\infty_{-\ell}(\mathbb{R}^2))}, \|h\|_{\tilde{B}_Y})]\\ &\lesssim R_{2,\varepsilon}(\|f_\varepsilon e^{\alpha|g_\varepsilon * X_0|}\|_{L^\infty}, e^{\alpha^2 \sup_{t\geq 0}\mathbb{E}[(g_\varepsilon * X_t)^2]}, e^{\alpha|h|_{\tilde{B}_Y}}),\end{aligned}$$

where $R_{2,\varepsilon}$ is a suitable second degree polynomial. Taking the supremum over $h \in \tilde{B}_Y$ with $|h| \leq 1$, recalling that $\tilde{B}^*_Y \subset B^{(2-s)\wedge(\gamma(r-1))+\delta}_{l,l,-\ell}(\mathbb{R}^2)$, for some $l \in (1,+\infty)$, and taking $f_{G^\lambda_\varepsilon}(X_0)$ as in Definition 3.1 to be proportional to $\|f_\varepsilon \exp(\alpha|g_\varepsilon * X_0|)\|^2_{L^\infty} + 1$, we get that $G^\lambda_\varepsilon$ satisfies the first condition of Definition 3.1. For the remaining conditions of Definition 3.1, similar arguments with the application of inequalities (3.13) and (3.14) allow us to conclude that $G^\lambda_\varepsilon \in \mathscr{F}$.

Since in our approximation we have that $\mathscr{G}_\varepsilon(X, Y) = \alpha f_\varepsilon e^{\alpha(g_\varepsilon *(X+Y))-\frac{\alpha^2}{2}c_\varepsilon}$, then the sum $\varphi^\varepsilon_t = X^\varepsilon_t + Y^\varepsilon_t$ solves the SPDE

$$(\partial_t - \Delta + m^2)\varphi^\varepsilon_t = -\alpha f_\varepsilon e^{\alpha g_\varepsilon * \varphi^\varepsilon_t - \frac{\alpha^2}{2}c_\varepsilon} + \xi_t.$$

Therefore, if $F \in \mathrm{Cyl}^b_{B_X \times B_Y}$ is of the form $F = \bar{F} \circ P^{X+Y}$ for some $\bar{F} \in \mathrm{Cyl}^b_E$, we have

$$G^\lambda_\varepsilon(X_0, Y_0) = \bar{G}^\lambda_\varepsilon(X_0 + Y_0) = \mathbb{E}_\xi\left[\int_0^{+\infty} e^{-\lambda t}\bar{F}(\varphi^\varepsilon_t)\,\mathrm{d}t \,\bigg|\, \varphi_0 = X_0 + Y_0\right].$$

**Proof of point ii.** In the following, for $j = 1, 2$, we denote by $(X_t^{\varepsilon,j}, Y_t^{\varepsilon,j})$ the solution to the system of equations (3.3)–(3.4) with initial conditions $(X_0, Y_0) \sim \mu_j$. Sometimes we also write $(X_0^j, Y_0^j)$ to indicate that $(X_0, Y_0) \sim \mu_j$, for $j = 1, 2$. By definition of $G_\varepsilon^\lambda$ (i.e. by equation (3.2)), we have

$$\lambda \int G_\varepsilon^\lambda(X_0, Y_0)\, \mu_1(dX_0, dY_0) - \lambda \int G_\varepsilon^\lambda(X_0, Y_0)\, \mu_2(dX_0, dY_0)$$
$$= \lambda \iint_0^{+\infty} e^{-\lambda t} \mathbb{E}_\xi[F(X_t^{\varepsilon,1}, Y_t^{\varepsilon,1})]\, dt\, d\mu_1 - \lambda \iint_0^{+\infty} e^{-\lambda t} \mathbb{E}_\xi[F(X_t^{\varepsilon,2}, Y_t^{\varepsilon,2})]\, dt\, d\mu_2$$
$$= \lambda \int_0^{+\infty} e^{-\lambda t} \mathbb{E}[F(X_t^{\varepsilon,1}, Y_t^{\varepsilon,1}) - F(X_t^{\varepsilon,2}, Y_t^{\varepsilon,2})]\, dt,$$

Notice that, for any $\varsigma > 0$, there exist two compact sets $K_1, K_2$, such that $\mu_1(K_1), \mu_2(K_2) > 1 - \varsigma$. Therefore, the following inequality holds

$$\left| \lambda \int_0^{+\infty} e^{-\lambda t} \mathbb{E}[F(X_t^{\varepsilon,1}, Y_t^{\varepsilon,1}) - F(X_t^{\varepsilon,2}, Y_t^{\varepsilon,2})]\, dt \right|$$
$$\leq \left| \lambda \int_0^{+\infty} e^{-\lambda t} \mathbb{E}[|F(X_t^{\varepsilon,1}, Y_t^{\varepsilon,1}) - F(X_t^{\varepsilon,2}, Y_t^{\varepsilon,2})| \mathbb{I}_{(X_0^1, Y_0^1) \in K_1, (X_0^2, Y_0^2) \in K_2}]\, dt \right| + 2\varsigma \|F\|_{L^\infty} \lambda \int_0^{+\infty} e^{-\lambda \tau} d\tau.$$

Let us then focus on the case where $(X_0^1, Y_0^1) \in K_1$ and $(X_0^2, Y_0^2) \in K_2$. By Lagrange's theorem, we have

$$F(X_t^{\varepsilon,1}, Y_t^{\varepsilon,1}) - F(X_t^{\varepsilon,2}, Y_t^{\varepsilon,2}) = \int_0^1 \langle dF_\beta, \mathfrak{D}_\beta \rangle\, d\beta,$$

where, for $\beta \in (0, 1)$,

$$dF_\beta = dF(X_t^\varepsilon((1-\beta) X_0^1 + \beta X_0^2), Y_t^\varepsilon((1-\beta) Y_0^1 + \beta Y_0^2)),$$

with $dF(X, Y) = \sum_{j=1}^N \partial_j \tilde{F} \langle u_j, X \rangle + \sum_{j=1}^M \partial_j \tilde{F} \langle v_j, Y \rangle$, for $F \in \mathrm{Cyl}_{B_X \times B_Y}^b$ having the following form $F(X, Y) = \tilde{F}(\langle u_1, X \rangle, \ldots, \langle u_N, X \rangle, \langle v_1, Y \rangle, \ldots, \langle v_M, Y \rangle)$, with $\tilde{F} \in C_b^2(\mathbb{R}^{N+M})$ and $u_i \in \mathcal{S}(\mathbb{R}^2)$, $i = 1, \ldots, N$, and $v_i \in \mathcal{S}(\mathbb{R}^2)$, $i = 1, \ldots, M$, and, for $\beta \in (0, 1)$,

$$\mathfrak{D}_\beta = (\nabla_{X_0^1 - X_0^2} X_t^\varepsilon((1-\beta) X_0^1 + \beta X_0^2), 0, \nabla_{X_0^1 - X_0^2} Y_t^\varepsilon((1-\beta) Y_0^1 + \beta Y_0^2), \nabla_{Y_0^1 - Y_0^2} Y_t^\varepsilon((1-\beta) Y_0^1 + \beta Y_0^2)).$$

We want to get some bounds for $\left| \int_0^1 \langle dF_\beta, \mathfrak{D}_\beta \rangle d\beta \right|$.

For the term $|\langle dF_\beta, \nabla_{X_0^1 - X_0^2} X_t^\varepsilon((1-\beta) X_0^1 + \beta X_0^2) \rangle|$ we have

$$|\langle dF_\beta, \nabla_{X_0^1 - X_0^2} X_t^\varepsilon((1-\beta) X_0^1 + \beta X_0^2) \rangle| = |\langle dF_\beta, P_t((1-\beta) X_0^1 + \beta X_0^2) \rangle|$$
$$\leq e^{-m^2 t} \|dF_\beta\| \|(1-\beta) X_0^1 + \beta X_0^2\|_{C_t^{-\delta}}.$$

For the term $|\langle dF_\beta, \nabla_{Y_0^1 - Y_0^2} Y_t^\varepsilon((1-\beta) Y_0^1 + \beta Y_0^2) \rangle|$, considering equation (3.8), and multiplying it by $\tilde{g}_\varepsilon * \nabla_{Y_0} Y_t$ and exploiting the negativity of the non-linearity, we get the following a priori estimate for some constant $k \in (0, m^2)$

$$\partial_t \|\tilde{g}_\varepsilon * \nabla_{Y_0^1 - Y_0^2} Y_t^\varepsilon((1-\beta) Y_0^1 + \beta Y_0^2)\|_{L^2}^2 + k \|\tilde{g}_\varepsilon * \nabla_{Y_0^1 - Y_0^2} Y_t^\varepsilon((1-\beta) Y_0^1 + \beta Y_0^2)\|_{L^2}^2 \leq 0,$$

and therefore, by Gronwall lemma,

$$\|\tilde{g}_\varepsilon * \nabla_{Y_0^1 - Y_0^2} Y_t^\varepsilon((1-\beta) Y_0^1 + \beta Y_0^2)\|_{L^2}^2 \lesssim e^{-kt},$$

and the estimate is independent of $\varepsilon$.

Consider now the operator $f \mapsto \tilde{g}_\varepsilon * f$, and recall that in Fourier representation convolution corresponds to multiplication, i.e. the previous operator can be viewed as $\hat{f} \mapsto \hat{\tilde{g}}_\varepsilon \cdot \hat{f}$. When $\hat{f}$ has compact support, such an operator is invertible if $0 < \varepsilon \leq \varepsilon_0$ for some positive constant $\varepsilon_0$ depending on the size of the support of $\hat{f}$, with inverse given by

$$\hat{\mathcal{Q}}_\varepsilon : \hat{f} \mapsto \frac{1}{\hat{\tilde{g}}_\varepsilon} \cdot \hat{f},$$

which is a well-defined operation because $\hat{\tilde{g}}_\varepsilon(k) = \hat{\tilde{g}}_1(\varepsilon k)$, $\hat{\tilde{g}}_1(0) = 1$, and $\hat{\tilde{g}}_1$ is smooth since $\tilde{g}_1$ is a Schwartz function. Using the previous notation, we have

$$\begin{aligned}|\langle dF_\beta, \nabla_{Y_0^1 - Y_0^2} Y_t^\varepsilon((1-\beta)Y_0^1 + \beta Y_0^2)\rangle| &= |\langle dF_\beta \circ \mathcal{Q}_\varepsilon, \tilde{g}_\varepsilon * \nabla_{Y_0^1 - Y_0^2} Y_t^\varepsilon((1-\beta)Y_0^1 + \beta Y_0^2)\rangle| \\ &\lesssim e^{-kt} \|dF_\beta \circ \mathcal{Q}_\varepsilon\| \\ &\lesssim e^{-kt} \sup_{0 < \varepsilon \leq \varepsilon_0} \|dF_\beta \circ \mathcal{Q}_\varepsilon\|,\end{aligned}$$

where we used the fact that $F$ has compact support in Fourier variables and so the norm $\|dF_\beta \circ \mathcal{Q}_\varepsilon\|$ is bounded for any $0 < \varepsilon \leq \varepsilon_0$, and also that

$$\sup_{\substack{0 < \varepsilon < \varepsilon_0 \\ k \in \text{supp}(F)}} |\hat{\tilde{g}}_\varepsilon(k)|^{-1} < +\infty.$$

Now, focus on the term $|\langle dF_\beta, \nabla_{X_0^1 - X_0^2} Y_t^\varepsilon((1-\beta)Y_0^1 + \beta Y_0^2)\rangle|$. By equation (3.9), we have

$$\begin{aligned}&\partial_t \|\nabla_{X_0^1 - X_0^2} Y_t^\varepsilon((1-\beta)Y_0^1 + \beta Y_0^2)\|_{L^2}^2 + \|\nabla_{X_0^1 - X_0^2} Y_t^\varepsilon((1-\beta)Y_0^1 + \beta Y_0^2)\|_{H^1}^2 \\ &+ m^2 \|\nabla_{X_0^1 - X_0^2} Y_t^\varepsilon((1-\beta)Y_0^1 + \beta Y_0^2)\|_{L^2}^2 \\ &\lesssim \int_{\mathbb{R}^2} D_X \mathcal{G}_\varepsilon(X_t^\varepsilon((1-\beta)X_0^1 + \beta X_0^2), Y_t^\varepsilon((1-\beta)Y_0^1 + \beta Y_0^2)) \nabla_{X_0^1 - X_0^2} X_t^\varepsilon((1-\beta)X_0^1 + \beta X_0^2) \\ &\quad \times \nabla_{X_0^1 - X_0^2} Y_t^\varepsilon((1-\beta)Y_0^1 + \beta Y_0^2) \, dz,\end{aligned}$$

and by Hölder and Young inequalities, and re-absorbing the terms properly, we have that, for any $\varsigma > 0$, the following inequality holds

$$\begin{aligned}&\partial_t \|\nabla_{X_0^1 - X_0^2} Y_t^\varepsilon((1-\beta)Y_0^1 + \beta Y_0^2)\|_{L^2}^2 + k\|\nabla_{X_0^1 - X_0^2} Y_t^\varepsilon((1-\beta)Y_0^1 + \beta Y_0^2)\|_{L^2}^2 \\ &\lesssim \frac{1}{\varsigma} \|D_X \mathcal{G}_\varepsilon(X_t^\varepsilon((1-\beta)X_0^1 + \beta X_0^2), Y_t^\varepsilon((1-\beta)Y_0^1 + \beta Y_0^2))\|_{B_{p,p}^{-s}}^2 \|\nabla_{X_0^1 - X_0^2} X_t^\varepsilon((1-\beta)X_0^1 + \beta X_0^2)\|_{B_{q,q}^{s}}^2.\end{aligned}$$

Now, we can bound the norm $\|\nabla_{X_0^1 - X_0^2} X_t^\varepsilon((1-\beta)X_0^1 + \beta X_0^2)\|_{B_{q,q}^s}^2$ by a constant since $X_0^1$ and $X_0^2$ belong to compact sets. Moreover, by Lemma A.7, we have

$$\begin{aligned}&\partial_t \|\nabla_{X_0^1 - X_0^2} Y_t^\varepsilon((1-\beta)Y_0^1 + \beta Y_0^2)\|_{L^2}^2 + k\|\nabla_{X_0^1 - X_0^2} Y_t^\varepsilon((1-\beta)Y_0^1 + \beta Y_0^2)\|_{L^2}^2 \\ &\lesssim \frac{e^{-m^2 t}}{t^{\delta + s}} \|D_X \mathcal{G}_\varepsilon(X_t^\varepsilon((1-\beta)X_0^1 + \beta X_0^2), 0)\|_{B_{p,p}^{-s}}^2.\end{aligned}$$

Applying Gronwall lemma yields the following inequality

$$\|\nabla_{X_0^1 - X_0^2} Y_t^\varepsilon((1-\beta)Y_0^1 + \beta Y_0^2)\|_{L^2}^2 \lesssim e^{-kt} \int_0^t \frac{e^{(k - m^2)\tau}}{\tau^{\delta + s}} \|D_X \mathcal{G}_\varepsilon(X_\tau^\varepsilon((1-\beta)X_0^1 + \beta X_0^2), 0)\|_{B_{p,p}^{-s}}^2 d\tau.$$

Thus, we get, for some constant $C_\varsigma > 0$ depending on $\varsigma > 0$,

$$\mathbb{E}[\|\nabla_{X_0^1 - X_0^2} Y_t^\varepsilon((1-\beta)Y_0^1 + \beta Y_0^2)\|_{L^2} \mathbb{I}_{(X_0^1, Y_0^1) \in K_1, (X_0^2, Y_0^2) \in K_2}]$$
$$\leq \mathbb{E}[\|\nabla_{X_0^1 - X_0^2} Y_t^\varepsilon((1-\beta)Y_0^1 + \beta Y_0^2)\|_{L^2}^2 \mathbb{I}_{(X_0^1, Y_0^1) \in K_1, (X_0^2, Y_0^2) \in K_2}]^{1/2}$$
$$\lesssim C_\varsigma \left( e^{-kt} \int_0^t \frac{e^{(k-m^2)\tau}}{\tau^{\delta+s}} \mathbb{E}[\|D_X \mathcal{G}_\varepsilon(X_\tau^\varepsilon((1-\beta)X_0^1 + \beta X_0^2), 0)\|_{B_{p,p}^{-s}}^2] d\tau \right)^{1/2}$$
$$\lesssim C_\varsigma \left( e^{-kt} \int_0^t \frac{e^{(k-m^2)\tau}}{\tau^{\delta+s}} d\tau \right)^{1/2},$$

since the law of $X_\tau^\varepsilon((1-\beta)X_0^1 + \beta X_0^2)$ is that of a Gaussian free field independent of $\beta \in (0,1)$ and of $\varepsilon > 0$.

Putting everything together, we get

$$\mathbb{E}[(F(X_t^{\varepsilon,1}, Y_t^{\varepsilon,1}) - F(X_t^{\varepsilon,2}, Y_t^{\varepsilon,2})) \mathbb{I}_{(X_0^1, Y_0^1) \in K_1, (X_0^2, Y_0^2) \in K_2}]$$
$$\lesssim e^{-m^2 t} \|dF_\beta\| \mathbb{E}[\|(1-\beta)X_0^1 + \beta X_0^2\|_{C_\ell^{-\delta}}] + e^{-kt} \|dF_\beta\| + C_\varsigma \left( e^{-kt} \int_0^t \frac{e^{(k-m^2)\tau}}{\tau^{\delta+s}} d\tau \right)^{1/2}$$
$$\lesssim e^{-m^2 t} \|dF_\beta\| + e^{-kt} \|dF_\beta\| + C_\varsigma \left( e^{-kt} \int_0^t \frac{e^{(k-m^2)\tau}}{\tau^{\delta+s}} d\tau \right)^{1/2},$$

where we used that the law of $(1-\beta)X_0^1 + \beta X_0^2$ does not depend on $\beta \in (0,1)$.

Therefore, we get

$$\left| \lambda \int G_\varepsilon^\lambda(X_0, Y_0) \mu_1(dX_0, dY_0) - \lambda \int G_\varepsilon^\lambda(X_0, Y_0) \mu_2(dX_0, dY_0) \right|$$
$$\lesssim \left| \lambda \int_0^{+\infty} \left[ e^{-\lambda t} e^{-m^2 t} \|dF_\beta\| + e^{-kt} \sup_{0 < \varepsilon \leq \varepsilon_0} \|dF_\beta \circ \mathcal{Q}_\varepsilon\| (1 + C_\varsigma) \left( 1 + \int_0^t \frac{e^{(k-m^2)\tau}}{\tau^{\delta+s}} d\tau \right)^{1/2} \right] dt \right|$$
$$+ 2\varsigma \|F\|_{L^\infty} \lambda \int_0^{+\infty} e^{-\lambda \tau} d\tau.$$

For some constant $K > 0$, we have

$$\left| \lambda \int_0^{+\infty} \left[ e^{-\lambda t} e^{-m^2 t} \|dF_\beta\| + e^{-kt} \sup_{0 < \varepsilon \leq \varepsilon_0} \|dF_\beta \circ \mathcal{Q}_\varepsilon\| (1 + C_\varsigma) \left( 1 + \int_0^t \frac{e^{(k-m^2)\tau}}{\tau^{\delta+s}} d\tau \right)^{1/2} \right] dt \right|$$
$$\lesssim (1 + C_\varsigma) \lambda \int_0^{+\infty} e^{-\lambda t - Kt} dt$$
$$\lesssim (1 + C_\varsigma) \frac{\lambda}{\lambda + K},$$

and the last term converges to zero as $\lambda \to 0$. Thus, we have

$$\left| \lambda \int G_\varepsilon^\lambda(X_0, Y_0) \mu_1(dX_0, dY_0) - \lambda \int G_\varepsilon^\lambda(X_0, Y_0) \mu_2(dX_0, dY_0) \right| \lesssim \varsigma + (1 + C_\varsigma) \frac{\lambda}{\lambda + K}.$$

**Proof of point iii.** Let $\mathcal{K} \in B_{\exp}^{q, \ell/2}$ and note that $\nabla_{Y_0} Y_t^\varepsilon[\mathcal{K}]$ solves equation (3.8). We write

$$\nabla \mathcal{Y}_t^\varepsilon = \nabla_{Y_0} Y_t^\varepsilon[\mathcal{K}] - P_t \mathcal{K}. \qquad (3.15)$$

Recall that in the proof of point *i.* in Proposition 3.5 we saw that

$$\nabla_{Y_0} G_\varepsilon^\lambda(X_0, Y_0) = \mathbb{E}_\xi \left[ \int_0^{+\infty} e^{-\lambda t} dF(X_t, Y_t^\varepsilon) \nabla_{Y_0} Y_t^\varepsilon \, dt \right].$$

Now, recalling the representation (3.15) of $\nabla \mathscr{Y}_t^\varepsilon$, we have

$$\int \mathbb{E}_\xi \left[ \int_0^{+\infty} e^{-\lambda t} dF(X_t, Y_t^\varepsilon) \nabla_{Y_0} Y_t^\varepsilon [\mathscr{K}] \, dt \right] d\mu$$
$$\lesssim \|dF\| \int \mathbb{E}_\xi \left[ \int_0^{+\infty} e^{-\lambda t} \|\nabla_{Y_0} Y_t^\varepsilon [\mathscr{K}]\|_{B_{\exp}^{q,\ell/2}} \, dt \right] d\mu$$
$$\lesssim \int \mathbb{E}_\xi \left[ \int_0^{+\infty} e^{-\lambda t} \left( \|P_t \mathscr{K}\|_{B_{\exp}^{q,\ell/2}} + \|\nabla \mathscr{Y}_t^\varepsilon\|_{L_\ell^2(\mathbb{R}^2)} \right) dt \right] d\mu.$$

By Lemma 3.9 below, we have

$$\int \mathbb{E}_\xi \left[ \int_0^{+\infty} e^{-\lambda t} dF(X_t, Y_t^\varepsilon) \nabla_{Y_0} Y_t^\varepsilon [\mathscr{K}] \, dt \right] d\mu$$
$$\lesssim \int \mathbb{E}_\xi \left[ \int_0^{+\infty} e^{-\lambda t} \left( 1 + (1+t)^\sigma \|f_\varepsilon : e^{\alpha(g_\varepsilon * X_s)} : \|_{L_{\ell''}^r(\mathbb{R}_+, B_{\exp}^{r,\ell/2})} \right) \|\mathscr{K}\|_{B_{\exp}^{q,\ell/2}} \, dt \right] d\mu.$$

Hölder inequality with respect to all measures, together with the fact that $\lambda > 0$, yields

$$\int \mathbb{E}_\xi \left[ \int_0^{+\infty} e^{-\lambda t} dF(X_t, Y_t^\varepsilon) \nabla_{Y_0} Y_t^\varepsilon [\mathscr{K}] \, dt \right] d\mu$$
$$\lesssim_\lambda \left( \int \|\mathscr{K}\|_{B_{\exp}^{q,\ell/2}}^q \, d\mu \right)^{1/q} \left( \int \mathbb{E}_\xi \left[ \|f_\varepsilon : e^{\alpha(g_\varepsilon * X_t)} : \|_{L_{\ell''}^r(\mathbb{R}_+, B_{\exp}^{r,\ell/2})}^r \right] d\mu \right)^{1/r} \int_0^{+\infty} e^{-\lambda t}(1 + (1+t)^\sigma) \, dt.$$

Here, since $\lambda > 0$, the integral $\int_0^{+\infty} e^{-\lambda t}(1 + (1+t)^\sigma) \, dt$ is finite. Moreover, by Proposition B.3, we have that the term

$$\left( \int \mathbb{E}_\xi \left[ \|f_\varepsilon : e^{\alpha(g_\varepsilon * X_s)} : \|_{L_{\ell''}^r(\mathbb{R}_+, B_{\exp}^{r,\ell/2})}^r \right] d\mu \right)^{1/r}$$

is bounded with respect to $\varepsilon > 0$. This concludes the proof of Proposition 3.5. □

We close the section with an auxiliary lemma used in the previous proof.

**Lemma 3.9.** *Let $\alpha^2 < 4\pi\tilde{\gamma}_{\max}$. Then, for every $\varepsilon > 0$, $Y_0 \in B_Y^\leqslant$, and $X_0$ in a set of full measure with respect to the free field measure $\mu^{\mathrm{free}}$ with mass $m > 0$, we have that*

$$\|\nabla \mathscr{Y}_t^\varepsilon\|_{L_\ell^2(\mathbb{R}^2)} \lesssim (1+t)^\sigma \|f_\varepsilon : e^{\alpha(g_\varepsilon * X_s)} : \|_{L_{\ell''}^r(\mathbb{R}_+, B_{\exp}^{r,\ell/2})} \|\mathscr{K}\|_{B_{\exp}^{q,\ell/2}}.$$

**Proof.** Let $\gamma = \alpha^2/(4\pi) < \tilde{\gamma}_{\max}$. The proof is similar to the one of Theorem C.1. Indeed, by multiplying the equation for $\nabla \mathscr{Y}^\varepsilon$ by $g_\varepsilon * \nabla \mathscr{Y}^\varepsilon$, and integrating, we get

$$\rho_{\ell'}(t) \|\tilde{g}_\varepsilon * \nabla \mathscr{Y}_t^\varepsilon\|_{L_\ell^2}^2 + \int_0^t \rho_{\ell'}(s) \|\tilde{g}_\varepsilon * \nabla \mathscr{Y}_s^\varepsilon\|_{H_\ell^1}^2 \, ds$$
$$+ \int_0^t \rho_{2\ell}(z) \rho_{\ell'}(s) \alpha^2 f_\varepsilon : e^{\alpha(g_\varepsilon * X_s)(z)} : e^{\alpha(g_\varepsilon * Y_s)(z)} (g_\varepsilon * \nabla \mathscr{Y}_s^\varepsilon(z))^2 \, ds \, dz$$
$$\lesssim \int_0^t \rho_{\ell'}(s) (\rho_{2\ell}(z) f_\varepsilon : e^{\alpha(g_\varepsilon * X_s)(z)} : e^{\alpha(g_\varepsilon * Y_s)(z)} (g_\varepsilon * e^{-(-\Delta + m^2)s} \mathscr{K})) g_\varepsilon * \nabla \mathscr{Y}_s^\varepsilon(z) \, ds \, dz.$$

In the last line, we have

$$\int_0^t \rho_\ell(s)\rho_{2\ell}(z)f_\varepsilon \underbrace{:e^{\alpha(g_\varepsilon*X_s)}:}_{L^p(\mathbb{R},B^{r,\ell/2}_{\exp}(\mathbb{R}^2))} \underbrace{e^{\alpha(g_\varepsilon*Y_s)}}_{L^\infty(\mathbb{R}_+\times\mathbb{R}^2)} \underbrace{\left(g_\varepsilon * e^{-(-\Delta+m^2)s}\mathcal{K}\right)}_{L^\infty\left(t^{-\beta}\mathbb{R}_+, B^{-\gamma(q-1)-\delta+2\beta-2/q}_{q,q,\ell/2}(\mathbb{R}^2)\right)} \underbrace{g_\varepsilon * \nabla\mathcal{Y}^\varepsilon_s}_{\substack{L^\infty_{\ell/2}(\mathbb{R}_+,L^2_t(\mathbb{R}^2))\cap \\ \cap L^2_{\ell/2}(\mathbb{R}_+,H^1_t(\mathbb{R}^2))}} ds.$$

Now, whenever $\gamma < \tilde{\gamma}_{\max}$, there exist $p, q \geq 1$, $0 \leq \theta \leq 1$, $0 \leq \beta \leq 1$, $\delta > 0$, such that the following system has a solution

$$\begin{cases} -\gamma(q-1) + 2\beta - \frac{2}{q} > 0, \\ -\gamma(p-1) - \delta + \theta > 0, \\ \frac{1}{p} + \frac{1}{2} < 1, \\ \frac{1}{p} + \beta - \delta + \frac{\theta}{2} \leq 1. \end{cases}$$

The proof is complete. □

## 3.3 Proof of uniqueness of solutions to Problem B

In the case of Problem B a better result can be obtained. It is useful to introduce a new approximating operator $\mathscr{L}_{\varepsilon,\bar{\varepsilon}}$ for $\mathscr{L}$ of the form

$$\mathscr{L}_{\varepsilon,\bar{\varepsilon}}(\Phi) := \frac{1}{2}\mathrm{tr}(\nabla_X^2\Phi) - \langle(-\Delta+m^2)X, \nabla_X\Phi\rangle - \langle(-\Delta+m^2)Y + \mathscr{G}_{\varepsilon,\bar{\varepsilon}}(X,Y), \nabla_Y\Phi\rangle,$$

where

$$\mathscr{G}_{\varepsilon,\bar{\varepsilon}}(X,Y) := \alpha f_\varepsilon(:e^{\alpha(g_\varepsilon*X)}:e^{\alpha(g_\varepsilon*Y)}).$$

Here $f_\varepsilon$ and $g_\varepsilon$ are defined as in Section 2.2. Recall that by Remark 2.9, since $\lim_{\bar{\varepsilon},\varepsilon\to 0}\mathscr{L}_{\varepsilon,\bar{\varepsilon}} = \mathscr{L}$, Problem B with the operator $\mathscr{L}_{\varepsilon,\bar{\varepsilon}}$ is equivalent to the one with operator $\mathscr{L}_\varepsilon$ defined in equation (2.14). We can then consider the resolvent equation associated to $\mathscr{L}_{\varepsilon,\bar{\varepsilon}}$, namely

$$(\lambda - \mathscr{L}_{\varepsilon,\bar{\varepsilon}})G^\lambda_{\varepsilon,\bar{\varepsilon}} = F,$$

where $F \in \mathrm{Cyl}_{B_X \times B_Y}$ with compact support in Fourier variables. A solution to such equation is given by

$$G^\lambda_{\varepsilon,\bar{\varepsilon}}(X_0, Y_0) = \mathbb{E}_\xi\left[\int_0^{+\infty} e^{-\lambda t}F(X_t, Y^{\varepsilon,\bar{\varepsilon}}_t)dt\right], \tag{3.16}$$

where $X_t, Y^{\varepsilon,\bar{\varepsilon}}_t$ solves the following system of equations

$$(\partial_t - \Delta + m^2)X_t = \xi_t, \qquad X(0) = X_0, \tag{3.17}$$
$$(\partial_t - \Delta + m^2)Y^{\varepsilon,\bar{\varepsilon}}_t = -\mathscr{G}_{\varepsilon,\bar{\varepsilon}}(X_t, Y^{\varepsilon,\bar{\varepsilon}}_t), \qquad Y^{\varepsilon,\bar{\varepsilon}}(0) = Y_0. \tag{3.18}$$

It is easy to show that all the results in Proposition 3.8 hold also for equations (3.17)–(3.18), adapting the form of the equations for the derivatives. This implies that points *i.*, *ii.* stated in Proposition 3.5 hold true also for $G^\lambda_{\varepsilon,\bar{\varepsilon}}$. As far as point *iii* in Proposition 3.5 is concerned, we can get a slightly better result in the present scenario. Recall that $\gamma_{\max}$ is defined as in Remark 1.1.

Notice that the operator (3.16) and equations (3.17)–(3.18) can be defined also for the case $\bar{\varepsilon} = 0$. Moreover, (a suitable adaptation of) Proposition 3.8 holds also for the case $\bar{\varepsilon} = 0$, and point *iii.* of Proposition 3.5 holds also for $G^\lambda_{\varepsilon,0}$.

**Remark 3.10.** It is worth to note that the operator $\mathscr{L}_{\varepsilon,\bar{\varepsilon}}$ cannot be use directly to solve Problem A'', since the solution to the resolvent equation depends in general on $(X, Y)$ and not only on $X + Y$.

**Proposition 3.11.** *Let $F \in \mathrm{Cyl}_{B_X \times B_Y}$ and consider $G^\lambda_{\varepsilon,\bar{\varepsilon}}$ given by (3.16). If $\alpha^2 < 4\pi\gamma_{\max}$, then there exists $q > 1$ such that, for every $\mu \in \mathscr{M}_{B_X \times B_Y}$, we have*

$$\lim_{\bar{\varepsilon}\to 0}\int |\langle \nabla_{Y_0} G^\lambda_{\varepsilon,\bar{\varepsilon}}, \mathscr{G} - \mathscr{G}_{\varepsilon,\bar{\varepsilon}}\rangle|\,d\mu \lesssim_\lambda \left(\int \|\mathscr{G} - \mathscr{G}_{\varepsilon,0}\|^q_{B^{q,\ell/2}_{\exp}}\,d\mu\right)^{1/q},$$

*uniformly in $\varepsilon > 0$.*

In order to prove the previous result, we need some technical lemmas. First, we deal with the convergence as $\bar{\varepsilon} \to 0$. In particular, whenever an object has $\bar{\varepsilon}$ as one of its parameters (e.g. $\nabla_{Y_0} G^\lambda_{\varepsilon,\bar{\varepsilon}}$), the same notation with $\bar{\varepsilon} = 0$ indicates that it is the limiting object as $\bar{\varepsilon} \to 0$ (e.g. $\nabla_{Y_0} G^\lambda_{\varepsilon,0} = \lim_{\bar{\varepsilon}\to 0}\nabla_{Y_0} G^\lambda_{\varepsilon,\bar{\varepsilon}}$), whenever it exists.

**Lemma 3.12.** *For every $\varepsilon > 0$, $Y_0 \in B^\leqslant_Y$, and $X_0$ in a set of full measure with respect to the free field measure $\mu^{\mathrm{free}}$ with mass $m > 0$, we have that $\nabla_{Y_0} Y^{\varepsilon,\bar{\varepsilon}}_t(\mathscr{G} - \mathscr{G}_{\varepsilon,\bar{\varepsilon}})$ converges to $\nabla_{Y_0} Y^{\varepsilon,0}_t(\mathscr{G} - \mathscr{G}_{\varepsilon,0})$ in $B_{\exp}$, as $\bar{\varepsilon} \to 0$.*

**Proof.** Recall that $\nabla_{Y_0} Y^{\varepsilon,\bar{\varepsilon}}_t(\mathscr{G} - \mathscr{G}_{\varepsilon,\bar{\varepsilon}})$ solves equation (3.8). This means that

$$\nabla \mathscr{Y}^{\varepsilon,\bar{\varepsilon}}_t = \nabla_{Y_0} Y^{\varepsilon,\bar{\varepsilon}}_t(\mathscr{G} - \mathscr{G}_{\varepsilon,\bar{\varepsilon}}) - P_t(\mathscr{G} - \mathscr{G}_{\varepsilon,\bar{\varepsilon}})$$

solves the equation

$$\begin{aligned}(\partial_t - \Delta + m^2)\nabla\mathscr{Y}^{\varepsilon,\bar{\varepsilon}}_t &= -\mathrm{D}_Y\mathscr{G}_{\varepsilon,\bar{\varepsilon}}(X_t, Y^{\varepsilon,\bar{\varepsilon}}_t)[\nabla\mathscr{Y}^{\varepsilon,\bar{\varepsilon}}_t] - \mathrm{D}_Y\mathscr{G}_{\varepsilon,\bar{\varepsilon}}(X_t, Y^{\varepsilon,\bar{\varepsilon}}_t)[P_t(\mathscr{G} - \mathscr{G}_{\varepsilon,\bar{\varepsilon}})], \\ \nabla\mathscr{Y}^{\varepsilon,\bar{\varepsilon}}_0 &= 0.\end{aligned} \quad (3.19)$$

By Corollary C.3, we get the estimate

$$\|\nabla\mathscr{Y}^{\varepsilon,\bar{\varepsilon}}_t\| \lesssim \rho_{\ell''}(t) P_2(\|f_\varepsilon : e^{\alpha(g_\varepsilon * X_t)}:\|_{L^\infty}, \|\mathscr{G} - \mathscr{G}_{\varepsilon,\bar{\varepsilon}}\|_{B_{\exp}}).$$

Such a bound is uniform with respect to $\bar{\varepsilon} > 0$ and therefore we get a converging subsequence whose limit $\nabla\mathscr{Y}^{\varepsilon,0}_t$ solves equation (3.19) with $\bar{\varepsilon} = 0$. Since the term $P_t(\mathscr{G} - \mathscr{G}_{\varepsilon,\bar{\varepsilon}})$ converges pointwise to $P_t(\mathscr{G} - \mathscr{G}_{\varepsilon,0})$, the result is proved. □

Let us prove that it makes sense to consider some parameters satisfying certain conditions that will be useful in upcoming results.

**Lemma 3.13.** *Suppose that $\gamma = \alpha^2/(4\pi) < \gamma_{\max}$. Then there exist $q > 1$ such that $\gamma q < 2$, $\kappa > 0$, and $r$ and $\delta$ as in Definition 2.6, such that the following inequalities are satisfied*

$$-\gamma(r-1) - \gamma(q-1) - 2\delta + \kappa > 0, \qquad \frac{1}{q} + \frac{1}{r} < 1, \qquad \frac{\kappa}{2}q < 1. \tag{3.20}$$

**Proof.** In order to prove the result, it is enough to show that there exists some solution with the previous properties to the system of equations

$$-\gamma(r-1) - \gamma\frac{1}{r-1} + \frac{2(r-1)}{r} = 0, \qquad \frac{1}{q} + \frac{1}{r} = 1, \qquad \frac{\kappa}{2}q = 1.$$

From such relations we get the equality

$$\gamma = \frac{2(r-1)^2}{r((r-1)^2+1)}. \tag{3.21}$$

We get the maximum value $\gamma_{\max}$ of $\gamma$ for some $r = \bar{r} \approx 2.52$. Since, with this choices of parameters, we have $q \approx 1.21$, the result is proved. □

In the proof of Lemma 3.12, we introduced the object

$$\nabla\mathcal{Y}_t^{\varepsilon,\bar{\varepsilon}} = \nabla_{Y_0} Y_t^{\varepsilon,\bar{\varepsilon}}(\mathcal{G} - \mathcal{G}_{\varepsilon,\bar{\varepsilon}}) - P_t(\mathcal{G} - \mathcal{G}_{\varepsilon,\bar{\varepsilon}}), \tag{3.22}$$

satisfying equation (3.19) and admitting a limit $\nabla\mathcal{Y}_t^{\varepsilon,0}$ as $\bar{\varepsilon} \to 0$.

**Lemma 3.14.** *Suppose that $\alpha^2 < 4\pi\gamma_{\max}$, consider the parameters $q$ and $r$ defined as in Lemma 3.13, and $\ell, \ell'''$ such that $\ell q/2 > 2$, $\ell r/2 > 2$ and $\ell''' r > 1$. Then, for every $\varepsilon > 0$, $Y_0 \in B_Y^{\leqslant}$, and $X_0$ in a set of full measure with respect to the free field measure $\mu^{\text{free}}$ with mass $m > 0$, we have that*

$$\|\nabla\mathcal{Y}_t^{\varepsilon,0}\|_{L_t^1(\mathbb{R}^2)} \lesssim (1+t)^\sigma \|f_\varepsilon : e^{\alpha(g_\varepsilon * X_s)}:\|_{L_{\ell'''}^r(\mathbb{R}_+, B_{\exp}^{r,\ell/2}(\mathbb{R}^2))} \|\mathcal{G} - \mathcal{G}_{\varepsilon,0}\|_{B_{\exp}^{q,\ell/2}(\mathbb{R}^2)}.$$

**Proof.** Let $\gamma = \alpha^2/(4\pi)$. With the same notation as in Lemma 3.12, we consider the limit $\nabla\mathcal{Y}_t^{\varepsilon,0}$ and show that it converges to zero. Indeed, it satisfies

$$(\partial_t - \Delta_z + m^2 + \alpha^2 f_\varepsilon : e^{\alpha(g_\varepsilon * X_t)(z)}: e^{\alpha Y_t(z)})\nabla\mathcal{Y}_t^{\varepsilon,0} = -\alpha^2 f_\varepsilon : e^{\alpha(g_\varepsilon * X_t)(z)}: e^{\alpha Y_t(z)} P_t(\mathcal{G} - \mathcal{G}_{\varepsilon,0}),$$

with $\nabla\mathcal{Y}_0^{\varepsilon,0} = 0$.

We now want to exploit a similar argument as the one used in Theorem C.2 to get some estimates concerning the solution to the previous equation. Let $w: \mathbb{R} \to \mathbb{R}$ be an increasing smooth function with bounded derivatives and such that $w(0) = 0$, $w(x) \to \pm 1$ as $x \to \pm\infty$, and define $W(x) = \int_0^x w(y)\,dy$. Then, multiplying the equation by $\rho_\ell(\theta\cdot)w(\tilde{\mathcal{K}}_t^{\varepsilon,0})$, where $\theta, \ell > 0$, and integrating, we have

$$\partial_t\|\rho_\ell(\theta\cdot)W(\nabla\mathcal{Y}_t^{\varepsilon,0})\|_{L^1} + m_\theta^2\|\rho_\ell(\theta\cdot)w(\nabla\mathcal{Y}_t^{\varepsilon,0})\tilde{\mathcal{K}}_t^{\varepsilon,0}\|_{L^1} + \|\rho_\ell(\theta\cdot)w'(\nabla\mathcal{Y}_t^{\varepsilon,0})(\nabla(\nabla\mathcal{Y}_t^{\varepsilon,0}))^2\|_{L^1}$$
$$+ \int \alpha^2\rho_\ell(\theta z) f_\varepsilon : e^{\alpha(g_\varepsilon * X_t)(z)}: e^{\alpha Y_t(z)} w(\nabla\mathcal{Y}_t^{\varepsilon,0})\nabla\mathcal{Y}_t^{\varepsilon,0}\,dz$$
$$\lesssim -\int \alpha^2\rho_{\ell/2}(\theta z) f_\varepsilon : e^{\alpha(g_\varepsilon * X_t)(z)}: e^{\alpha Y_t(z)} w(\nabla\mathcal{Y}_t^{\varepsilon,0})\rho_{\ell/2}(\theta z)(P_t(\mathcal{G} - \mathcal{G}_{\varepsilon,0}))\,dz, \tag{3.23}$$

where $0 < m_\theta \leqslant \sqrt{m^2 - \left|\frac{\nabla\rho_\ell(\theta z)}{\rho_\ell(\theta z)}\right|}$ for every $z \in \mathbb{R}^2$, which holds for $\theta > 0$ small enough.

If we choose the parameters as in Lemma 3.13, we have

$$\rho_{\ell/2}(\theta\cdot)f_\varepsilon\mathpunct{:}e^{\alpha(g_\varepsilon * X_t)}\mathpunct{:}e^{\alpha Y_t} \in L^r_{\ell'''}(\mathbb{R}_+, B^{r,0}_{\exp}),$$
$$w(\nabla \mathcal{Y}^{\varepsilon,0}_t) \in L^\infty(\mathbb{R}^2),$$
$$\rho_{\ell/2}(\theta\cdot)(P_t(\mathcal{G} - \mathcal{G}_{\varepsilon,0})) \in L^\infty(|\cdot|^{-\kappa/2}\mathbb{R}_+, B^{-\gamma(q-1)-\delta+\kappa}_{q,q}(\mathbb{R}^2)).$$

Therefore, by integrating with respect to time equation (3.23) and recalling that $\nabla \mathcal{Y}^{\varepsilon,0}_0 = 0$, we get

$$\|\rho_\ell(\theta\cdot)W(\nabla \mathcal{Y}^{\varepsilon,0}_t)\|_{L^1} \lesssim$$
$$\int_0^t \|\rho_{\ell/2}(\theta\cdot)f_\varepsilon\mathpunct{:}e^{\alpha(g_\varepsilon * X_s)}\mathpunct{:}e^{\alpha Y_s}w(\nabla \mathcal{Y}^{\varepsilon,0}_s)\|_{B^{r,0}_{\exp}}\|\rho_{\ell/2}(\theta\cdot)(e^{-(-\Delta+m^2)s}(\mathcal{G}-\mathcal{G}_{\varepsilon,0}))\|_{B^{-\gamma(q-1)-\delta+\kappa}_{q,q}}ds$$

where we used the fact that $\frac{1}{r} + \frac{1}{q} < 1$ and $-\gamma(r-1) - \gamma(q-1) - 2\delta + \kappa > 1$. By Proposition A.4, we have that

$$\|f_\varepsilon\mathpunct{:}e^{\alpha(g_\varepsilon * X_s)}\mathpunct{:}e^{\alpha Y_s}w(\nabla \mathcal{Y}^{\varepsilon,0}_s)\|_{B^{r,\ell/2}_{\exp}} \lesssim \|f_\varepsilon\mathpunct{:}e^{\alpha(g_\varepsilon * X_s)}\mathpunct{:}\|_{B^{r,\ell/2}_{\exp}}\|e^{\alpha Y_s}\|_{L^\infty}\|w(\nabla \mathcal{Y}^{\varepsilon,0}_s)\|_{L^\infty}.$$

Therefore, by Proposition A.6, we have

$$\|\rho_\ell(\theta\cdot)W(\nabla \mathcal{Y}^{\varepsilon,0}_t)\|_{L^1} \lesssim \int_0^t \|f_\varepsilon\mathpunct{:}e^{\alpha(g_\varepsilon * X_s)}\mathpunct{:}\|_{B^{r,\ell/2}_{\exp}} s^{\kappa/2}\|\mathcal{G} - \mathcal{G}_{\varepsilon,0}\|_{B^{q,\ell/2}_{\exp}}ds.$$

Multiply and divide then by $\rho_{\ell'''}(s)$ and apply Hölder inequality with respect to $s$ to get

$$\|\rho_\ell(\theta\cdot)W(\nabla \mathcal{Y}^{\varepsilon,0}_t)\|_{L^1}$$
$$\lesssim \|\rho_{\ell/2}(\theta\cdot)f_\varepsilon\mathpunct{:}e^{\alpha(g_\varepsilon * X_s)}\mathpunct{:}\|_{L^r_{\ell'''}(\mathbb{R}_+, B^{r,\ell/2}_{\exp})}\|\mathcal{G} - \mathcal{G}_{\varepsilon,0}\|_{B^{q,\ell/2}_{\exp}}\left(\int_0^t s^{q\kappa/2}\rho_{-\ell'''q}(s)ds\right)^{1/q}.$$

Since the expression $(\int_0^t s^{q\kappa/2}\rho_{-\ell'''q}(s)ds)^{1/q}$ is bounded for $t \to 0$ and grows polynomially in time as $t \to +\infty$, we deduce that the result holds. □

**Proof of Proposition 3.11.** As in the proof of point *i.* of Proposition 3.5 we can get

$$\nabla_{Y_0} G^\lambda_{\varepsilon,\bar{\varepsilon}}(X_0, Y_0) = \mathbb{E}_\xi\left[\int_0^{+\infty} e^{-\lambda t}dF(X_t, Y^{\varepsilon,\bar{\varepsilon}}_t)\nabla_{Y_0} Y^{\varepsilon,\bar{\varepsilon}}_t dt\right].$$

Therefore, by Lemma 3.12 and Lebesgue's dominated convergence theorem, we have

$$\lim_{\bar{\varepsilon}\to 0}\int \nabla_{Y_0} G^\lambda_{\varepsilon,\bar{\varepsilon}}(X_0, Y_0)(\mathcal{G} - \mathcal{G}_{\varepsilon,\bar{\varepsilon}})d\mu = \lim_{\bar{\varepsilon}\to 0}\int \mathbb{E}_\xi\left[\int_0^{+\infty} e^{-\lambda t}dF(X_t, Y^{\varepsilon,\bar{\varepsilon}}_t)\nabla_{Y_0} Y^{\varepsilon,\bar{\varepsilon}}_t(\mathcal{G}-\mathcal{G}_{\varepsilon,\bar{\varepsilon}})dt\right]d\mu$$
$$= \int \mathbb{E}_\xi\left[\int_0^{+\infty} e^{-\lambda t}dF(X_t, Y^{\varepsilon,0}_t)\nabla_{Y_0} Y^{\varepsilon,0}_t(\mathcal{G}-\mathcal{G}_{\varepsilon,0})dt\right]d\mu$$
$$= \int \nabla_{Y_0} G^\lambda_{\varepsilon,0}(X_0, Y_0)(\mathcal{G} - \mathcal{G}_{\varepsilon,0})d\mu.$$

Now, recalling the representation (3.22) of $\nabla \mathcal{Y}^{\varepsilon,0}_t$, we have

$$\int \mathbb{E}_\xi\left[\int_0^{+\infty} e^{-\lambda t}dF(X_t, Y^{\varepsilon,0}_t)\nabla_{Y_0} Y^{\varepsilon,0}_t(\mathcal{G}-\mathcal{G}_{\varepsilon,0})dt\right]d\mu$$
$$\lesssim \|dF\|_{L^\infty}\int \mathbb{E}_\xi\left[\int_0^{+\infty} e^{-\lambda t}\|\nabla_{Y_0} Y^{\varepsilon,0}_t(\mathcal{G}-\mathcal{G}_{\varepsilon,0})\|_{B^{q,\ell/2}_{\exp}}dt\right]d\mu$$
$$\lesssim \int \mathbb{E}_\xi\left[\int_0^{+\infty} e^{-\lambda t}\left(\|P_t(\mathcal{G}-\mathcal{G}_{\varepsilon,0})\|_{B^{q,\ell/2}_{\exp}} + \|\nabla \mathcal{Y}^{\varepsilon,0}_t\|_{L^1_t(\mathbb{R}^2)}\right)dt\right]d\mu.$$

At this point, the proof follows the same steps as the one of point *iv.* in Proposition 3.5. □

# 4 Existence of solutions via Lyapunov functions

In this section, we prove existence of solutions to Problem A" and B. Since existence for the Problem B implies existence for Problem A" (see point *iii.* in Theorem 2.5), we only prove the first. We actually prove a stronger statement than the one given in Problem B, namely that equation (2.5) holds for every $\Phi \in \mathscr{F}$ and not only for $\Phi \in \mathrm{Cyl}_{B_X \times B_Y}$. We exploit a strategy based on Lyapunov functions.

Let us introduce a finite-dimensional approximation to the operator $\mathscr{L}$. Let $M, N \in \mathbb{N}$, $\varepsilon > 0$, and consider $\mathbb{T}_M^2$ the two-dimensional torus of size $M$ which we identify hereafter with the subset $(-M\pi, M\pi]^2 \subset \mathbb{R}^2$. Consider the *Fejér operator* $Q_{N,M} : \mathscr{S}'(\mathbb{T}_M^2) \to C^\infty(\mathbb{T}_M^2)$ defined as, for $F \in \mathscr{S}'(\mathbb{T}_M^2)$,

$$Q_{N,M}(F) = \mathrm{Fej}_{N,M} * F, \qquad \mathrm{Fej}_{N,M}(x) = \sum_{|j_1|, |j_2| \leq N-1} \left(1 - \frac{|j_1|}{N}\right)\left(1 - \frac{|j_2|}{N}\right) e^{ijx/M}, \qquad x \in \mathbb{T}_M^2.$$

Let us stress that the operator $Q_{N,M}$ is both a positive operator, i.e. $\langle Q_{N,M}(F), F \rangle \geq 0$, and a positive preserving operator, namely, if $F$ is a positive distribution, then $Q_{N,M}(F)$ is a positive function. The latter property following from the positivity of the kernel $\mathrm{Fej}_{N,M}$.

The new approximating operator $\mathscr{L}_{N,M,\varepsilon}$ is then given by, for $\Phi \in \mathscr{F}$, by

$$\mathscr{L}_{N,M,\varepsilon}(\Phi) = \frac{1}{2} \mathrm{tr}_{L^2(\mathbb{T}_M^2)}(P_{\mathbb{T}_M^2}(\mathrm{Per}_{\mathbb{T}_M^2} \nabla_X^2 \Phi) P_{\mathbb{T}_M^2}) - \langle (-\Delta + m^2)X, \nabla_X \Phi \rangle$$
$$- \langle (-\Delta + m^2)Y + \mathscr{G}_{N,M,\varepsilon}(X, Y), \nabla_Y \Phi \rangle,$$

where $\mathscr{G}_{N,M,\varepsilon}$ is the following approximation for the non-linearity in equation (3.4),

$$\mathscr{G}_{N,M,\varepsilon}(X, Y) = \alpha Q_{N,M}(g_\varepsilon * (:e^{\alpha Q_{N,M}(g_\varepsilon * X)}: e^{\alpha Q_{N,M}(g_\varepsilon * Y)})). \qquad (4.1)$$

Here, $g_\varepsilon$ is defined as in Section 2.2, $P_{\mathbb{T}_M^2}$ is the natural projection of $L^2(\mathbb{R}^2)$ on the space $L^2(\mathbb{T}_M^2)$, and $\mathrm{Per}_{\mathbb{T}_M^2}$ is the periodicization on the torus $\mathbb{T}_M^2$. Moreover, the term $:\exp(\alpha Q_{N,M}(g_\varepsilon * X)):$ is defined as follows

$$:e^{\alpha Q_{N,M}(g_\varepsilon * X)}: = e^{\alpha Q_{N,M}(g_\varepsilon * X)} - \frac{\alpha^2}{2} c_{N,M,\varepsilon}, \qquad (4.2)$$

where $c_{N,M,\varepsilon} = \int_{\mathbb{T}_M^2} Q_{N,M}(g_\varepsilon)(z)(-\Delta + m^2)^{-1}(Q_{N,M}g_\varepsilon)(z)\, dz$, where the inverse $(-\Delta + m^2)^{-1}$ is taken with periodic boundary conditions on $\mathbb{T}_M^2$.

The system of equations for the flow is then given by

$$(\partial_t - \Delta + m^2)X_t^M = \xi_t^M, \qquad (4.3)$$
$$(\partial_t - \Delta + m^2)Y_t^{N,M,\varepsilon} = -\mathscr{G}_{N,M,\varepsilon}(X_t^M, Y_t^{N,M,\varepsilon}), \qquad (4.4)$$

where $Y^{N,M,\varepsilon}$ is negative, periodic, and it belongs to a subspace of $\mathfrak{H}_{N,M} := \mathrm{span}\{e^{in \cdot/M}; |n| \leq N\}$. As usual, we drop the dependence on the parameters $N, M, \varepsilon$ when no ambiguity occurs.

## 4.1 Lyapunov functions

We introduce some Lyapunov functions $V_1$, $V_2$, and $V_3$ that will be crucial in the proof of the estimates for $X$ and $Y$. In particular, we take them such that $V_1$ and $V_2$ depend on both $X$ and $Y$, while $V_3$ depends on $X$ only.

**Lemma 4.1.** *Consider two Banach spaces $B_1$ and $B_2$. Let $\tilde{\mathscr{L}}$ be an operator taking values in some space of functions $\mathscr{F}$ and consider a measure $\tilde{\mu}$ such that $\int \tilde{\mathscr{L}}\Phi \, d\tilde{\mu} = 0$ for every $\Phi \in \mathscr{F}$. Suppose that there exist some functions $V_1, V_2 : B_1 \times B_2 \to \mathbb{R}$, and $V_3 : B_1 \to \mathbb{R}$, such that we have $V_1, V_2, V_3 \in \mathscr{F}$ and*

  i. *$V_2$ and $V_3$ are positive,*

  ii. *The inequality*

$$\tilde{\mathscr{L}} V_1(X, Y) \leq -V_2(X, Y) + V_3(X)$$

  *holds true.*

*Then, we have*

$$\int V_2(X, Y) \tilde{\mu}(dX, dY) \leq \int V_3(X) \tilde{\mu}(dX, dY). \tag{4.5}$$

**Proof.** The statement follows from the fact that $\int \tilde{\mathscr{L}} V_1(X, Y) \tilde{\mu}(dX, dY) = 0$. □

We now work on the estimate for $X$. We will have to choose the aforementioned Lyapunov functions, we start with $V_1 : B_X \to \mathbb{R}$, that we take to be, for $s' > 0$,

$$V_1(X) = \|X\|_{B_{p,p,\ell}^{-s'}}^p.$$

Note that, since $X$ corresponds to the free field, the previous norm is finite. Notice also that the previous function is not cylindric. Moreover, if we consider the representation of Besov spaces given in Proposition A.9, then we have, for $-ps'/2 + 1 < 0$,

$$\|X\|_{B_{p,p,\ell}^{-s'}}^p \simeq \|\varphi_0(D)X\|_{L_\ell^p}^p + \int_0^{+\infty} \int_{\mathbb{R}^2} \frac{|\rho_\ell^k|^p}{t^{-ps'/2+1}} |P_t X(z)|^p \, dz \, dt =: V_{1,1}(X) + V_{1,2}(X),$$

where $k$ is to be determined, similarly as in the proof of Proposition 2.10.

**Proposition 4.2.** *Let $s' > 0$, $p > 1$ and $\ell > 0$ such that $p\ell > 2$ and $-ps'/2 + 1 < 0$. Consider $X$ to be the solution to equation (4.3) on the torus $\mathbb{T}_M^2$. For some $\delta \in (0, 1)$, we have the inequality*

$$\mathscr{L}_{N,M,\varepsilon} V_1(X) = \mathscr{L}_{N,M,\varepsilon}(\|X\|_{B_{p,p,\ell}^{-s'}}^p) \lesssim -(1-\delta)\|X\|_{B_{p,p,\ell}^{-s'+2/p}}^p + \frac{1}{\delta}C. \tag{4.6}$$

**Proof.** We want to evaluate $\mathscr{L}_{N,M,\varepsilon} V_1(X)$. The gradient of $V_1$ is given by

$$\nabla V_1(X)(h) = \nabla \|\varphi_0(D)X\|_{L_\ell^p}^p(h) + \int_0^{+\infty} \int_{\mathbb{R}^2} \frac{|\rho_\ell^k|^p}{t^{-ps'/2+1}} \nabla(|P_t X(z)|^p)(h) \, dz \, dt$$

$$= p\int_{\mathbb{R}^2} |\rho_\ell^k|^p |\varphi_0(D)X(z)|^{p-1} \text{sign}(\varphi_0(D)X) \varphi_0(D) h \, dz$$

$$+ p\int_0^{+\infty} \int_{\mathbb{R}^2} \frac{|\rho_\ell^k|^p}{t^{-ps'/2+1}} |P_t X(z)|^{p-1} \text{sign}(P_t X(z)) P_t h \, dz \, dt,$$

while its second derivative is

$$\begin{aligned}
&\nabla^2 V_1(X)(h,h') \\
&= p\int_{\mathbb{R}^2} |\rho_\ell^k|^p \nabla(|\varphi_0(D)X(z)|^{p-1}\operatorname{sign}(\varphi_0(D)X)\,\varphi_0(D)h)(h')\,dz \\
&\quad + p\int_0^{+\infty}\int_{\mathbb{R}^2} \frac{|\rho_\ell^k|^p}{t^{-ps'/2+1}} \nabla(|P_t X(z)|^{p-1}\operatorname{sign}(P_t X(z))\,P_t h)(h')\,dz\,dt \\
&= p(p-1)\int_{\mathbb{R}^2} |\rho_\ell^k|^p |\varphi_0(D)X(z)|^{p-2}(\varphi_0(D)h)(\varphi_0(D)h')\,dz \\
&\quad + p(p-1)\int_0^{+\infty}\int_{\mathbb{R}^2} \frac{|\rho_\ell^k|^p}{t^{-ps'/2+1}} |P_t X(z)|^{p-2}(P_t h)(P_t h')\,dz\,dt.
\end{aligned}$$

The first term to deal with is $\langle(-\Delta+m^2)X,\nabla_X V_1(X)\rangle$. We have

$$\begin{aligned}
-\langle(-\Delta+m^2)X,\nabla_X V_{1,1}(X)\rangle &= -p\int_{\mathbb{R}^2}|\rho_\ell^k|^p |\varphi_0(D)X(z)|^{p-1}\operatorname{sign}(\varphi_0(D)X)(\varphi_0(D)(-\Delta+m^2)X)\,dz \\
&= -p\int_{\mathbb{R}^2} \nabla(|\rho_\ell^k|^p |\varphi_0(D)X(z)|^{p-1}\operatorname{sign}(\varphi_0(D)X))\nabla(\varphi_0(D)X)\,dz \\
&\quad - pm^2\int_{\mathbb{R}^2}|\rho_\ell^k|^p |\varphi_0(D)X(z)|^p\,dz \\
&= -p(p-1)\int_{\mathbb{R}^2}|\rho_\ell^k|^p |\varphi_0(D)X(z)|^{p-2}|\nabla(\varphi_0(D)X)|^2\,dz \\
&\quad - p^2\int_{\mathbb{R}^2}|\rho_\ell^k|^{p-1}|\varphi_0(D)X(z)|^{p-1}(\nabla\rho_\ell^k)\operatorname{sign}(\varphi_0(D)X)\nabla(\varphi_0(D)X)\,dz \\
&\quad - pm^2\int_{\mathbb{R}^2}|\rho_\ell^k|^p |\varphi_0(D)X(z)|^p\,dz.
\end{aligned}$$

Rearranging the second term on the right-hand side and exploiting Young's inequality, we have, for all $\sigma>0$,

$$\begin{aligned}
-\langle(-\Delta+m^2)X,\nabla_X V_{1,1}(X)\rangle &\lesssim -p(p-1)\int_{\mathbb{R}^2}|\rho_\ell^k|^p |\varphi_0(D)X(z)|^{p-2}|\nabla(\varphi_0(D)X)|^2\,dz \\
&\quad + p^2\int_{\mathbb{R}^2}|\rho_\ell^k|^p |\varphi_0(D)X(z)|^{p/2}\frac{|\nabla\rho_\ell^k|}{\rho_\ell^k}|\varphi_0(D)X|^{(p-2)/2}|\nabla(\varphi_0(D)X)|\,dz \\
&\quad - pm^2\int_{\mathbb{R}^2}|\rho_\ell^k|^p |\varphi_0(D)X(z)|^p\,dz \\
&\lesssim -p(p-1)\int_{\mathbb{R}^2}|\rho_\ell^k|^p |\varphi_0(D)X(z)|^{p-2}|\nabla(\varphi_0(D)X)|^2\,dz \\
&\quad + \frac{p^2}{4\sigma}\int_{\mathbb{R}^2}|\rho_\ell^k|^p |\varphi_0(D)X(z)|^p\frac{|\nabla\rho_\ell^k|^2}{\rho_{2\ell}^k}\,dz \\
&\quad + p^2\sigma\int_{\mathbb{R}^2}|\rho_\ell^k|^p |\varphi_0(D)X|^{p-2}|\nabla(\varphi_0(D)X)|^2\,dz \\
&\quad - pm^2\int_{\mathbb{R}^2}|\rho_\ell^k|^p |\varphi_0(D)X(z)|^p\,dz,
\end{aligned}$$

from which, choosing $\sigma>0$ small enough and $k>0$ large enough, and appropriately rearranging the terms, we get

$$-\langle(-\Delta+m^2)X,\nabla_X V_{1,1}(X)\rangle \lesssim -\|\varphi_0(D)X\|_{L_\ell^p}^p.$$

We also have

$$-\langle(-\Delta + m^2)X, \nabla_X V_{1,2}(X)\rangle$$
$$= -p\int_0^{+\infty}\int_{\mathbb{R}^2} \frac{|\rho_\ell|^p}{t^{-ps'/2+1}} |P_t X(z)|^{p-1}\mathrm{sign}(P_t X(z))(P_t(-\Delta+m^2)X(z))\,\mathrm{d}z\,\mathrm{d}t$$
$$= \int_0^{+\infty}\int_{\mathbb{R}^2} \frac{|\rho_\ell|^p}{t^{-ps'/2+1}} \partial_t |P_t X(z)|^p \,\mathrm{d}z\,\mathrm{d}t$$
$$= -\int_0^{+\infty}\int_{\mathbb{R}^2} \partial_t\left(\frac{|\rho_\ell|^p}{t^{-ps'/2+1}}\right) |P_t X(z)|^p \,\mathrm{d}z\,\mathrm{d}t$$
$$= -\int_0^{+\infty}\int_{\mathbb{R}^2}\left(\frac{ps'}{2}-1\right)\frac{|\rho_\ell|^p}{t^{-ps'/2+2}} |P_t X(z)|^p \,\mathrm{d}z\,\mathrm{d}t$$
$$= -\int_0^{+\infty}\int_{\mathbb{R}^2}\left(\frac{ps'}{2}-1\right)\frac{|\rho_\ell|^p}{t^{p(-s'+2/p)/2+1}} |P_t X(z)|^p \,\mathrm{d}z\,\mathrm{d}t$$
$$\simeq -(\|X\|^p_{B^{-s'+2/p}_{p,p,\ell}} - \|\varphi_0(D)X\|^p_{L^p_\ell}).$$

The second term to deal with is the trace of the second derivative, i.e. $\mathrm{tr}_{L^2(\mathbb{T}_M^2)}(\nabla^2 V_1(X))$. By the results presented in Section XI.V of [RS79], it suffices to consider the second derivative with $h' = h$, i.e. $\nabla^2 V_1(X)(h, h)$, and integrate with respect to $h$. Then, exploiting Young inequality for products and using the fact that $p\ell > 2$, we have for any $\sigma > 0$

$$\mathrm{tr}(\nabla^2 V_{1,1}(X)) = p(p-1)\int_{\mathbb{T}_M^2}\int_{\mathbb{R}^2} |\rho_\ell^k|^p |\varphi_0(D)X(z)|^{p-2} |\hat{\varphi}_0(z-y)|^2 \,\mathrm{d}z\,\mathrm{d}y$$
$$= p(p-1)\int_{\mathbb{R}^2} |\rho_\ell^k|^p |\varphi_0(D)X(z)|^{p-2} \int_{\mathbb{T}_M^2} |\hat{\varphi}_0(z-y)|^2 \,\mathrm{d}y\,\mathrm{d}z$$
$$\lesssim p(p-1)\int_{\mathbb{R}^2} |\rho_\ell^k|^p |\varphi_0(D)X(z)|^{p-2} \,\mathrm{d}z$$
$$\lesssim \frac{1}{\sigma}p(p-1)\int_{\mathbb{R}^2} |\rho_\ell^k|^p \,\mathrm{d}z + \sigma p(p-1)\int_{\mathbb{R}^2} |\rho_\ell^k|^p |\varphi_0(D)X(z)|^p \,\mathrm{d}z$$
$$\lesssim \frac{1}{\sigma}C + \sigma \|\varphi_0(D)X\|^p_{L^p_\ell}.$$

We get, if $\mathscr{P}_t^M$ is the heat kernel on the torus $\mathbb{T}_M^2$, that, for any $\sigma > 0$,

$$\mathrm{tr}(\nabla^2 V_{1,2}(X)) = p(p-1)\int_{\mathbb{T}_M^2}\int_0^{+\infty}\int_{\mathbb{R}^2} \frac{|\rho_\ell|^p}{t^{-ps'/2+1}} |P_t X(z)|^{p-2} |\mathscr{P}_t^M(z-y)|^2 \,\mathrm{d}z\,\mathrm{d}t\,\mathrm{d}y$$
$$= p(p-1)\int_0^{+\infty}\int_{\mathbb{R}^2} \frac{|\rho_\ell|^p}{t^{-ps'/2+1}} |P_t X(z)|^{p-2}\left(\int_{\mathbb{T}_M^2} |\mathscr{P}_t^M(z-y)|^2 \,\mathrm{d}y\right)\,\mathrm{d}z\,\mathrm{d}t$$
$$= p(p-1)\int_0^{+\infty}\int_{\mathbb{R}^2} \frac{|\rho_\ell|^p}{t^{-ps'/2+1}} |P_t X(z)|^{p-2}\frac{Ce^{-m^2 t}}{t}\,\mathrm{d}z\,\mathrm{d}t$$
$$= p(p-1)\int_0^{+\infty}\int_{\mathbb{R}^2} \frac{Ce^{-m^2 t}|\rho_\ell|^p}{t^{-ps'/2+2}} |P_t X(z)|^{p-2}\,\mathrm{d}z\,\mathrm{d}t$$
$$\lesssim \frac{1}{\sigma}C'\int_0^{+\infty}\int_{\mathbb{R}^2} \frac{|\rho_\ell|^p e^{-q(p)m^2 t}}{t^{-ps'/2+2}}\,\mathrm{d}z\,\mathrm{d}t + \sigma\int_0^{+\infty}\int_{\mathbb{R}^2} \frac{|\rho_\ell|^p}{t^{-ps'/2+2}} |P_t X(z)|^p \,\mathrm{d}z\,\mathrm{d}t$$
$$\lesssim \frac{1}{\sigma}C + \sigma\left(\|X\|^p_{B^{-s'+2/p}_{p,p,\ell}} - \|\varphi_0(D)X\|^p_{L^p_\ell}\right),$$

where we used Young inequality as above.

Choosing $\sigma > 0$ small enough, we get

$$\mathscr{L}V_1(X) \lesssim -(1-\sigma)\|X\|^p_{B^{-s'+2/p}_{p,p,\ell}} + \frac{1}{\sigma}C,$$

which gives the result. $\square$

**Remark 4.3.** The result in Lemma 4.1 with $B_2 = B_X$, no dependence on $B_1$ whatsoever, and with Lyapunov functions chosen as

$$\begin{aligned}
V_1(X) &= \|X\|^p_{B^{-s'}_{p,p,\ell}}, \\
V_2(X) &= (1-\sigma)\|X\|^p_{B^{-s'+2/p}_{p,p,\ell}}, \\
V_3 &\equiv \frac{1}{\sigma}C,
\end{aligned}$$

for some constant $C > 0$, gives the estimate

$$(1-\sigma)\int \|X\|^p_{B^{-s'+2/p}_{p,p,\ell}} \mu_{N,M,\varepsilon}(\mathrm{d}X,\mathrm{d}Y) \lesssim \frac{1}{\sigma}C.$$

Notice further that the bound (4.6) can be chosen to be uniform with respect to the size $M$ of the torus $\mathbb{T}^2_M$ since the integrals $\int_{\mathbb{T}^2_M} |\hat{\varphi}_0(z-y)|^2 \, \mathrm{d}y$ and $\int_{\mathbb{T}^2_M} \|\mathscr{P}^M_t(z-y)|(z-y)|^2 \, \mathrm{d}y$ are uniformly bounded in $M$. It is worth to note also that by the Besov embedding (Proposition A.2) the space $B^{-s'+2/p}_{p,p,\ell}(\mathbb{R}^2)$ is embedded in $B_X = C^{-\delta}_\ell(\mathbb{R}^2)$.

Let us choose the Lyapunov functions in order to get an estimate for $Y$. We may need $X \mapsto V_1(X, \cdot)$ to be $C^1$ and $Y \mapsto V_1(\cdot, Y)$ to be $C^2$.

Let $s > 0$. Referring to the representation of Besov spaces given in Proposition A.9 and Remark A.10, take $(k-s/2)p > 1$, and consider

$$V_1(Y) = \|Y\|^p_{B^s_{p,p,\ell}} \simeq \|Y\|^p_{L^p_\ell} + \int_0^{+\infty}\!\!\int_{\mathbb{R}^2} \frac{t^{kp}|\rho_\ell|^p}{t^{ps/2+1}} |\partial_{t^k} P_t Y(z)|^p \, \mathrm{d}z\,\mathrm{d}t =: V_{1,1}(Y) + V_{1,2}(Y). \quad (4.7)$$

**Theorem 4.4.** *Suppose that $Y$ is a solution to equation (4.4). Let $\gamma = \alpha^2/(4\pi) < 2$, and $1 < p, q < +\infty$. Take $s$ and $r$ such that*

$$0 < s < \gamma + 2 - \sqrt{8\gamma}, \quad \frac{2+\gamma-s-\sqrt{(s-\gamma-2)^2-8\gamma}}{2\gamma} < r < \frac{2+\gamma-s+\sqrt{(s-\gamma-2)^2-8\gamma}}{2\gamma}, \quad \gamma r < 2. \quad (4.8)$$

*If $: e^{\alpha X}: \in B^{r,\ell}_{\exp} = B^{-\gamma(r-1)-\delta}_{r,r,\ell}(\mathbb{R}^2)$, for $\ell$ large enough, then we have, for any $\sigma > 0$,*

$$\mathscr{L}_{N,M,\varepsilon}V_1(Y) = \mathscr{L}_{N,M,\varepsilon}(\|Y\|^p_{B^s_{p,p,\ell}}) \lesssim -(1-\sigma)\|Y\|^p_{B^{s+2/p}_{p,p,\ell}} + \frac{1}{\sigma}\|\mathscr{G}_{N,M,\varepsilon}(X,Y)\|^{(pr-r+1)/(pr^2)}_{B^{r,\ell}_{\exp}}. \quad (4.9)$$

**Proof.** In the following proof, we neglect the term $V_{1,1}(Y)$ appearing in (4.7), which can be dealt with in a similar way as in the proof of Proposition 4.2. The gradient of $V_1$ is given by

$$\nabla V_{1,2}(Y)(h) = \int_0^{+\infty} \int_{\mathbb{R}^2} \frac{t^{pk} |\rho_\ell|^p}{t^{ps/2+1}} \nabla(|\partial_{t^k} P_t Y(z)|^p)(h) \, dz \, dt$$

$$= \int_0^{+\infty} \int_{\mathbb{R}^2} \frac{t^{pk} |\rho_\ell|^p}{t^{ps/2+1}} p |\partial_{t^k} P_t Y(z)|^{p-1} \partial_{t^k} P_t h \, dz \, dt.$$

Therefore, taking $h = -(-\Delta + m^2)Y - \mathcal{G}_{N,M,\varepsilon}(X, Y)$, we have

$$-\langle (-\Delta + m^2)Y + \mathcal{G}_{N,M,\varepsilon}(X, Y), \nabla V_{1,2}(Y)\rangle$$

$$= \int_0^{+\infty} \int_{\mathbb{R}^2} \frac{t^{pk} |\rho_\ell|^p}{t^{ps/2+1}} p |\partial_{t^k} P_t Y(z)|^{p-1} \partial_{t^k} P_t(-(-\Delta + m^2)Y - \mathcal{G}_{N,M,\varepsilon}(X,Y)) \, dz \, dt$$

$$= \int_0^{+\infty} \int_{\mathbb{R}^2} \frac{t^{pk} |\rho_\ell|^p}{t^{ps/2+1}} \partial_t(|\partial_{t^k} P_t Y|^p) \, dz \, dt$$

$$+ \int_0^{+\infty} \int_{\mathbb{R}^2} \frac{t^{pk} |\rho_\ell|^p}{t^{ps/2+1}} p |\partial_{t^k} P_t Y(z)|^{p-1} \partial_{t^k} P_t(-\mathcal{G}_{N,M,\varepsilon}(X,Y)) \, dz \, dt$$

$$=: I_1 + I_2.$$

Let us focus on $I_1$. Integrating by parts, we get

$$I_1 = -\int_0^{+\infty} \int_{\mathbb{R}^2} \partial_t \left( \frac{t^{pk} |\rho_\ell|^p}{t^{ps/2+1}} \right) |\partial_{t^k} P_t Y|^p \, dz \, dt$$

$$= -\int_0^{+\infty} \int_{\mathbb{R}^2} \left( pk - \frac{ps}{2} - 1 \right) \frac{t^{pk} |\rho_\ell|^p}{t^{ps/2+2}} |\partial_{t^k} P_t Y|^p \, dz \, dt$$

$$= -\int_0^{+\infty} \int_{\mathbb{R}^2} \left( pk - \frac{ps}{2} - 1 \right) \frac{t^{pk} |\rho_\ell|^p}{t^{p(s+2/p)/2+1}} |\partial_{t^k} P_t Y|^p \, dz \, dt$$

$$\simeq -(\|Y\|_{B^{s+2/p}_{p,p,\ell}}^p - \|Y\|_{L^p_\ell}^p).$$

Consider now $I_2$. Exploiting Hölder inequality with $q_1, p_1$ and introducing $s' > 0$ yields

$$I_2 = \int_0^{+\infty} \int_{\mathbb{R}^2} \frac{t^{pk} |\rho_\ell|^p}{t^{ps/2+1}} p |\partial_{t^k} P_t Y(z)|^{p-1} \partial_{t^k} P_t(-\mathcal{G}_{N,M,\varepsilon}(X,Y)) \, dz \, dt$$

$$= \int_0^{+\infty} \frac{1}{t} \int_{\mathbb{R}^2} \left[ \frac{t^{(p-1)k} |\rho_\ell|^{p-1}}{t^{ps/2+s'/2}} p |\partial_{t^k} P_t Y(z)|^{p-1} \right]$$

$$\times [t^{k+s'/2} |\rho_\ell| \partial_{t^k} P_t(-\mathcal{G}_{N,M,\varepsilon}(X,Y))] \, dz \, dt$$

$$\leq \left( \int_0^{+\infty} \int_{\mathbb{R}^2} \left| \frac{t^{(p-1)k} |\rho_\ell|^{p-1}}{t^{ps/2+s'/2}} p |\partial_{t^k} P_t Y(z)|^{p-1} \right|^{p_1} dz \frac{dt}{t} \right)^{1/p_1}$$

$$\times \left( \int_0^{+\infty} \int_{\mathbb{R}^2} |t^{s'/2}(t^k |\rho_\ell|) \partial_{t^k} P_t(-\mathcal{G}_{N,M,\varepsilon}(X,Y))|^{q_1} dz \frac{dt}{t} \right)^{1/q_1}.$$

If we apply Young inequality, we then have, for any $\sigma > 0$,

$$I_2 \leq \sigma \left( \int_0^{+\infty} \int_{\mathbb{R}^2} \left| \frac{t^{(p-1)k} |\rho_\ell|^{p-1}}{t^{ps/2+s'/2}} p |\partial_{t^k} P_t Y(z)|^{p-1} \right|^{p_1} dz \frac{dt}{t} \right)^{p_2/p_1}$$

$$+ \frac{1}{\sigma} \left( \int_0^{+\infty} \int_{\mathbb{R}^2} |t^{s'/2}(t^k |\rho_\ell|) \partial_{t^k} P_t(-\mathcal{G}_{N,M,\varepsilon}(X,Y))|^{q_1} dz \frac{dt}{t} \right)^{q''/q_1}$$

$$= \sigma \left( \int_0^{+\infty} \int_{\mathbb{R}^2} \left( \frac{t^{(p-1)k} |\rho_\ell|^{p-1}}{t^{ps/2+s'/2}} \right)^{p_1} p^{p_1} |\partial_{t^k} P_t Y(z)|^{p_1 p - p_1} dz \frac{dt}{t} \right)^{p_2/p_1}$$

$$+\frac{1}{\sigma}\left(\int_0^{+\infty}\int_{\mathbb{R}^2}|t^{s'/2}(t^k|\rho_\ell|)\partial_{t^k}P_t(-\mathscr{G}_{N,M,\varepsilon}(X,Y))|^{q_1}\mathrm{d}z\frac{\mathrm{d}t}{t}\right)^{q''/q_1}$$

$$= \sigma\left(\int_0^{+\infty}\int_{\mathbb{R}^2}\frac{t^{p_1(p-1)k}}{t^{psp_1/2+s'p_1/2}}\left(\frac{t^k}{t^{(ps+s')/(2(p-1))}}\right)^{p_1(p-1)}p^{p_1}\right.$$

$$\left.\times|\rho_\ell|^{p_1(p-1)}|\partial_{t^k}P_tY(z)|^{p_1(p-1)}\mathrm{d}z\frac{\mathrm{d}t}{t}\right)^{p_2/p_1}$$

$$+\frac{1}{\sigma}\left(\int_0^{+\infty}\int_{\mathbb{R}^2}|t^{k+s'/2}|\rho_\ell|\partial_{t^k}P_t(-\mathscr{G}_{N,M,\varepsilon}(X,Y))|^{q_1}\mathrm{d}z\frac{\mathrm{d}t}{t}\right)^{q_2/q_1}$$

$$= \sigma\|Y\|_{B^{(ps+s')/(p-1)}_{p_1(p-1),p_1(p-1),\ell}}^{p_2/p_1}+\frac{1}{\sigma}\left(\int_0^{+\infty}\int_{\mathbb{R}^2}|t^{s'/2}(t^k|\rho_\ell|)\partial_{t^k}P_t(-\mathscr{G}_{N,M,\varepsilon}(X,Y))|^{q_1}\mathrm{d}z\frac{\mathrm{d}t}{t}\right)^{q_2/q_1}$$

$$=: I_{2,1}+I_{2,2}.$$

We want to apply Besov embedding (see Proposition A.2) in order to reabsorb the term $I_{2,1}$ in $I_1$: considering the parameters of the involved norms, we get the condition

$$s > \frac{ps+s'}{p-1} - \frac{2}{p_1(p-1)}, \qquad (4.10)$$

so that the reabsorbing procedure follows by a suitable choice of $p_2$.

On the other hand, we have by Proposition A.4

$$I_{2,2} \lesssim \frac{1}{\sigma}\|\mathscr{G}_{N,M,\varepsilon}(X,Y)\|_{B^{-s'}_{q_1,q_1,\ell}}^{q_2/q_1}.$$

All in all, we have the inequality

$$\mathscr{L}_{N,M,\varepsilon}V_1(Y) = \mathscr{L}_{N,M,\varepsilon}\|Y\|_{B^s_{p,p,\ell}}^p \lesssim -(1-\sigma)\|Y\|_{B^{s+2/p}_{p,p,\ell}}^p+\frac{1}{\sigma}\|\mathscr{G}_{N,M,\varepsilon}(X,Y)\|_{B^{-s'}_{q_1,q_1,\ell}}^{q_2/q_1},$$

provided condition (4.10) holds true and $p_2$ is chosen in an appropriate way.

Now, in order to obtain (4.9), we consider $q_1 = r$ and $s' > \gamma(r-1)$, and have to check the condition

$$\gamma r^2 + r(s-\gamma-2) + 2 < 0,$$

which is satisfied with the choice of parameters given in (4.8). □

**Remark 4.5.** If $Y \leq 0$, we have the bound

$$\|\mathscr{G}_{N,M,\varepsilon}(X,Y)\|_{B^{r,\ell}_{\exp}}^{(pr-r+1)/(pr^2)} \leq \|\mathscr{G}_{N,M,\varepsilon}(X,0)\|_{B^{r,\ell}_{\exp}}^{(pr-r+1)/(pr^2)},$$

on the right-hand side of inequality (4.9). In that case, choosing the Lyapunov functions as follows

$$V_1(Y) = \|Y\|_{B^s_{p,p,\ell}}^p,$$
$$V_2(Y) = C_\sigma\|Y\|_{B^{s+2/p}_{p,p,\ell}}^p,$$
$$V_3(X) = \frac{1}{\sigma}\|\mathscr{G}_{N,M,\varepsilon}(X,0)\|_{B^{r,\ell}_{\exp}}^{(pr-r+1)/(pr^2)},$$

we get that the conditions in Lemma 4.1 are satisfied thanks to Theorem 4.4. Therefore, we have

$$\int \|Y\|^p_{B^{s+2/p}_{p,p,\ell}} \mu_{N,M,\varepsilon}(\mathrm{d}X, \mathrm{d}Y) \leqslant \frac{1}{\sigma} \int \|\mathscr{G}_{N,M,\varepsilon}(X,0)\|^{(pr-r+1)/(pr^2)}_{B^{r,\ell}_{\exp}} \mu_{N,M,\varepsilon}(\mathrm{d}X, \mathrm{d}Y).$$

## 4.2 Measure of the approximating problem

Equations (4.3) and (4.4) induce a Markov process on $B_X \times B_Y^{\leqslant}$. This is due to the fact that, if we start from a negative initial condition for $Y$, then $Y$ remains negative throughout all the evolution.

**Proposition 4.6.** *The operator $(\mathscr{L}_{N,M,\varepsilon}, \mathscr{F})$ is the restriction of the generator of the Markov process associated to equations (4.3) and (4.4) to the space of functions $\mathscr{F}$.*

**Proof.** The proof is based on the fact that we can apply Itô formula to the functions in $\mathscr{F}$. Exploiting Itô formula, the proof follows the argument of the proof of point *ii.* in Proposition 3.5. □

By linearity, the invariant measure for equation (4.3) is given by the free field measure with mass $m$ on the torus $\mathbb{T}^2_M$, that we can identify as usual with the periodic free field of mass $m$ on the whole space $\mathbb{R}^2$.

We are interested in studying (infinitesimally) invariant measures for the aforementioned pair of equations, i.e. equations (4.3) and (4.4). Indeed, by Proposition 4.6, an invariant measure for those two equations is a solution to a FPK equation associated with $\mathscr{L}_{N,M,\varepsilon}$ in the sense of Problem B. Let us start off by giving an argument for the existence of an infinitesimally invariant measure (for its definition and the relation with invariant measures see Chapter 5 in [BKRS15]) by means of the Lyapunov functions introduced in Section 4.1.

**Proposition 4.7.** *Equations (4.3) and (4.4) admit an infinitesimally invariant measure.*

**Proof.** The proof follows the same argument as Lemma 3.3 in [ADRU22]. We start equations (4.3) and (4.4) from deterministic initial conditions $(x_0, y_0) \in B_X \times B_Y^{\leqslant}$, and let $\mu_{(x_0, y_0), t}$ denote the measure of the solution of the equations starting at $(x_0, y_0)$. Let

$$\tilde{\mu}_T = T^{-1} \int_0^T \mu_{(x_0, y_0), t} \, \mathrm{d}t.$$

If such a measure is tight with respect to $T$, and $\tilde{\mu}$ is its weak limit as $T \to +\infty$, then $\tilde{\mu}$ is an infinitesimal invariant measure of the equation. We consider the following Lyapunov functions, consisting of the sum of the one considered in Section 4.1 where we dealt with $X$ and $Y$ separately (see Remarks 4.3 and 4.5), namely

$$V_1(X, Y) = \|X\|^p_{B^{-s}_{p,p,\ell}} + \|Y\|^p_{B^s_{p,p,\ell}}, \tag{4.11}$$

$$V_2(X, Y) = (1-\sigma) \|X\|^p_{B^{-s+2/p}_{p,p,\ell}} + C_\sigma \|Y\|^p_{B^{s+2/p}_{p,p,\ell}}, \tag{4.12}$$

$$V_3(X) = \frac{1}{\sigma} \Big( C + \|\mathscr{G}_{N,M,\varepsilon}(X,0)\|^{(pr-r+1)/(pr^2)}_{B^{r,\ell}_{\exp}} \Big). \tag{4.13}$$

Then, by Proposition 4.2 and Theorem 4.4, we have

$$\begin{aligned}\mathbb{E}[V_1(X_t, Y_t) - V_1(x_0, y_0)] &= \mathbb{E}\int_0^t \mathscr{L}_{N,M,\varepsilon} V_1(X_\tau, Y_\tau)\,\mathrm{d}\tau \\ &\leq \int_0^t \mathbb{E}[-V_2(X_\tau, Y_\tau) + V_3(X_\tau)]\,\mathrm{d}\tau.\end{aligned}$$

Recall that $X_\tau = e^{-(-\Delta + m^2)\tau}x_0 + \int_0^\tau e^{-(-\Delta + m^2)(\tau - s)}\xi_s\,\mathrm{d}s$, then $\mathbb{E}[\|\mathscr{G}_{N,M,\varepsilon}(X_\tau, 0)\|_{B_{\exp}^{r,\ell}}^{(pr-r+1)/(pr^2)}]$ converges exponentially to a constant, and hence

$$K = \sup_{t \in [0,+\infty)} \frac{1}{t}\int_0^t \mathbb{E}\left[\|\mathscr{G}_{N,M,\varepsilon}(X_\tau, 0)\|_{B_{\exp}^{r,\ell}}^{(pr-r+1)/(pr^2)}\right]\mathrm{d}\tau < +\infty.$$

Therefore, we have, for any $t > 0$,

$$\frac{1}{t}\mathbb{E}[V_1(X_t, Y_t) - V_1(x_0, y_0)] + \frac{1}{t}\int_0^t \mathbb{E}[V_2(X_\tau, Y_\tau)]\,\mathrm{d}\tau \lesssim K,$$

which yields

$$\frac{1}{t}\int_0^t \mathbb{E}[V_2(X_\tau, Y_\tau)]\,\mathrm{d}\tau \lesssim K + \frac{1}{t}\mathbb{E}[V_1(x_0, y_0)].$$

Taking the supremum over $t \in [0, +\infty)$, we have

$$\sup_{t \in [0,+\infty)} \frac{1}{t}\int_0^t \mathbb{E}[V_2(X_\tau, Y_\tau)]\,\mathrm{d}\tau < +\infty.$$

On the other hand, for any $t > 0$,

$$\int V_2(X, Y)\tilde{\mu}_t(\mathrm{d}X, \mathrm{d}Y) = \frac{1}{t}\int_0^t \mathbb{E}[V_2(X_\tau, Y_\tau)]\,\mathrm{d}s.$$

Since $V_2(X, Y)$ has compact sub-levels on $B_X \times B_Y^{\leqslant}$, tightness of $\tilde{\mu}_t$ is implied. □

We can actually prove a stronger result on the form of one invariant measure for equations (4.3) and (4.4). Indeed, if we start from $x_0 = 0$ and $y_0 = 0$ in the proof of Proposition 4.7, then we have that the invariant measure that has been built there must be the measure of the process solving the equation

$$X_t = \int_{-\infty}^t P_{t-s}\xi_s\,\mathrm{d}s, \qquad Y_t = \int_{-\infty}^t P_{t-s}\mathscr{G}_{N,M,\varepsilon}(X_s, Y_s)\,\mathrm{d}s. \tag{4.14}$$

It is possible, in fact, to prove that the previous integral equations admit a unique solution.

**Proposition 4.8.** *Equation (4.14) admits a unique solution in $B_X \times \mathfrak{H}_{N,M} \subset B_X \times B_Y$. Moreover, for every $t \in \mathbb{R}$, the law of $(X_t, Y_t)$ is a solution to the FPK equation for $\mathscr{L}_{N,M,\varepsilon}$ in the sense of Problem B.*

**Proof.** Existence and uniqueness for equation (4.14) can be proved in the same way as in the proof of Theorem E.1.

In order to prove the second part of the result, we exploit the fact that

$$\text{law}(X_t, Y_t) = \lim_{t \to +\infty} \mu_{((0,0),t)} = \lim_{t \to +\infty} \frac{1}{t}\int_0^t \mu_{((0,0),s)} \, ds.$$

Indeed, if $X_t^{0,-T}$ and $Y_t^{0,-T}$ are the processes solution to equations (4.3) and (4.4) starting at time $-T$ with initial condition zero, then we have that

$$\text{law}(X_t^{0,-T}, Y_t^{0,-T}) = \mu_{((0,0),t)},$$

and that $(X_t^{0,-T}, Y_t^{0,-T}) \to (X_t, Y_t)$ as $T \to +\infty$, almost surely. Therefore, by the proof of Proposition 4.7, the law of $(X_t, Y_t)$ is an invariant measure for equations (4.3) and (4.4), and thus a solution to the FPK equation for $\mathscr{L}_{N,M,\varepsilon}$. □

Hereafter, we denote by $\mu_{N,M,\varepsilon}$ the law at fixed time $t$ of the process $(X_t, Y_t)$ solution to equation (4.14), which is a solution to the FPK equation for $\mathscr{L}_{N,M,\varepsilon}$. Thanks to the representation given by Proposition 4.8, we are able to establish some more precise properties of the measure $\mu_{N,M,\varepsilon}$.

**Proposition 4.9.** *For every $M, N \in \mathbb{N}$, and every $\varepsilon > 0$, the measure $\mu_{N,M,\varepsilon}$ satisfies the following properties*

  i.  $\text{supp}(\mu_{N,M,\varepsilon}) \subset B_{2,2,\ell}^2(\mathbb{T}_M^2) \times \mathfrak{H}_{N,M}$,

  ii. *For every $F \in \mathscr{F}$, we have*

$$\int \mathscr{L}_{N,M,\varepsilon} F \, d\mu_{N,M,\varepsilon} = 0,$$

  iii. *For every $F, G \in \mathfrak{F}$, we have*

$$\int \mathscr{L}_{N,M,\varepsilon} F(X+Y) G(X+Y) \mu_{N,M,\varepsilon}(dX, dY) = \int F(X+Y) \mathscr{L}_{N,M,\varepsilon} G(X+Y) \mu_{N,M,\varepsilon}(dX, dY).$$

**Proof.** Property *i.* follows from the fact that the solution $Y_t$ to the second equation in (4.14) is supported on the image of the projection of the operator $Q_{N,M}$, which is exactly $\mathfrak{H}_{N,M}$.

Property *ii.* is due to the fact that $\mu_{N,M,\varepsilon}$ solves the FPK equation for $\mathscr{L}_{N,M,\varepsilon}$. Indeed, the system (4.14) can be split in an infinite-dimensional linear equation and an independent finite-dimensional non-linear equation. The statement then follows from the results in Section 2.3 in [Da 04] for the infinite-dimensional part and Theorem 5.2.2 in [BKRS15] for the non-linear finite-dimensional part.

As far as point *iii.* is concerned, letting $Z_t = X_t + Y_t$ and calling $Z_t^1$ the projection of $Z_t$ on $\mathfrak{H}_{N,M}$ and $Z_t^2$ the projection on $\mathfrak{H}_{N,M}^\perp$, we have that $Z_t^1$ and $Z_t^2$ solve two independent equations: the equation for $Z_t^1$ is a non-linear finite-dimensional equation, since $\mathfrak{H}_{N,M}$ is a finite-dimensional Hilbert space a linear equation, while the equation for $Z_t^2$ is a linear equation in infinite dimensions. More precisely, we have

$$\begin{aligned}\partial_t Z_t^1 &= (\Delta - m^2) Z_t^1 + \frac{\delta V_{N,M,\varepsilon}}{\delta Z}(Z_t^1) + \xi_t^1, \\ \partial_t Z_t^2 &= (\Delta - m^2) Z_t^2 + \xi_t^2,\end{aligned}$$

where $\delta$ is the functional derivative of

$$V_{N,M,\varepsilon}(Z) = \int_{\mathbb{T}_M^2} \exp\left(\alpha Q_{N,M} Z(z) - \frac{\alpha^2}{2} c_{N,M,\varepsilon}\right) dz,$$

with $c_{N,M,\varepsilon}$ is the constant appearing in the definition (4.2) of the Wick product appearing in the functional $\mathcal{G}_{N,M,\varepsilon}$, and $\xi_t^1$ and $\xi_t^2$ are the projection of $\xi$ on the spaces $\mathfrak{H}_{N,M}$ and $\mathfrak{H}_{N,M}^\perp$, respectively.

Since the equation for $Z_t^1$ is a finite-dimensional equation with drift given by the gradient of a function and the equation for $Z_t^2$ is linear with self-adjoint drift, then the unique invariant measure of the process $(Z^1, Z^2)$ satisfies property *iii.*, and hence also their sum $Z^1 + Z^2 = X + Y$ does. For further details on the relation between symmetric processes and the integration by parts formula see e.g. [AR91, Bog10, Föl88]. □

**Remark 4.10.** From the proof of the previous result it is evident that the system consists of two independent equations: an infinite-dimensional linear one and a finite-dimensional non-linear one. This means that $\mu_{M,N,\varepsilon}$ is the unique invariant measure of the system.

## 4.3 Tightness of the measure

We prove tightness of the measures $\mu_{N,M,\varepsilon}$.

**Theorem 4.11.** *The family of measures $(\mu_{N,M,\varepsilon})_{N,M,\varepsilon}$ is tight in $B_X \times B_Y^\leqslant$.*

**Proof.** Let $V_1$, $V_2$, and $V_3$ be the Lyapunov functions defined in equations (4.11), (4.12), and (4.13), respectively. Since the measures $\mu_{N,M,\varepsilon}$ are solutions to the FPK equation associated with $\mathscr{L}_{N,M,\varepsilon}$ in the sense of Problem B and $V_1 \in \mathscr{F}$, then

$$\int \mathscr{L}_{N,M,\varepsilon} V_1(X,Y) \mu_{N,M,\varepsilon}(dX, dY) = 0.$$

Moreover, we have by Theorem 4.2 together with Remark 4.3 and by Theorem 4.4 together with Remark 4.5, that $V_1$ is a Lyapunov function for $(X, Y)$, where $X$ and $Y$ solves equation (4.3) and (4.4), respectively. Thus, by Lemma 4.1, we have

$$\int V_2(X,Y) \mu_{N,M,\varepsilon}(dX, dY) \lesssim \int V_3(X) \mu_{N,M,\varepsilon}(dX, dY) = \int V_3(X) \mu_M^{\text{free}}(dX).$$

Since $\sup_{M \in \mathbb{N}} \int V_3(X) \mu_M^{\text{free}}(dX) < +\infty$, and since $V_2$ has compact sub-levels, the thesis follows. □

We get the following consequence.

**Corollary 4.12.** *There exist three sequences $(N_k)_{k \in \mathbb{N}} \subset \mathbb{N}$, $(M_l)_{l \in \mathbb{N}} \subset \mathbb{N}$, and $(\varepsilon_j)_j \subset \mathbb{R}_+$ such that $N_k \to +\infty$, $M_l \to +\infty$, and $\varepsilon_j \to 0$ as $k, l, j \to +\infty$, respectively, and some probability measures $\mu_{M_l, \varepsilon_j}$, $\mu_{\varepsilon_j}$, and $\mu$ such that, for any $M_l$ and any $\varepsilon_j$, we have*

$$\lim_{k \to +\infty} \mu_{N_k, M_l, \varepsilon_j} = \mu_{M_l, \varepsilon_j},$$

*for any $\varepsilon_j$ we have*

$$\lim_{l \to +\infty} \mu_{M_l, \varepsilon_j} = \mu_{\varepsilon_j},$$

*and finally*

$$\lim_{j \to +\infty} \mu_{\varepsilon_j} = \mu.$$

**Proof.** The result follows from Theorem 4.11 and a diagonalization argument. □

## 4.4 Limit of the operator

We want now to prove that any limit measure $\mu$ built in Corollary 4.12 solves the FPK equation in the sense of Problem B. We first prove that any measure appearing in Corollary 4.12 solves the FPK equation associated with the corresponding operator. For simplicity, let us drop the dependence on the parameters $k$, $l$, and $j$ introduced in Corollary 4.12.

The aim of the section is to prove the following result.

**Theorem 4.13.** *Let $\Phi \in \mathcal{F}$. We have the limit as $N, M \to +\infty$ and $\varepsilon \to 0$,*

$$\int \mathscr{L}_{N,M,\varepsilon} \Phi \, d\mu_{N,M,\varepsilon} \to \int \mathscr{L} \Phi \, d\mu. \tag{4.15}$$

In order to prove the limit (4.15), we proceed by first taking $N \to +\infty$, then $M \to +\infty$, and eventually $\varepsilon \to 0$, showing the convergence step by step.

### 4.4.1 Limit as $N \to +\infty$

As mentioned above, we start off by taking $N \to +\infty$, namely we want to show that

$$\int \mathscr{L}_{N,M,\varepsilon}(\Phi) \, d\mu_{N,M,\varepsilon} \to \int \mathscr{L}_{M,\varepsilon}(\Phi) \, d\mu_{M,\varepsilon}, \qquad \text{as } N \to +\infty.$$

The term in the operator $\mathscr{L}_{M,N,\varepsilon}$ involving $X$ only is independent of $N$, in particular we do not need to consider the trace. Instead, we need to work on the terms involving only the derivatives with respect to $Y$. Since, by Corollary 4.12, $\mu_{N,M,\varepsilon}$ is tight and converges weakly up to subsequences to some limit $\mu_{M,\varepsilon}$, we have to show that

$$\int [\mathscr{L}_{N,M,\varepsilon}(\Phi(X,Y)) - \mathscr{L}_{M,\varepsilon}(\Phi(X,Y))] \mu_{N,M,\varepsilon}(dX, dY) \to 0, \qquad \text{as } N \to +\infty.$$

Let us compute the integrand and rearrange its terms. We have

$$\begin{aligned}
&\mathscr{L}_{N,M,\varepsilon}(\Phi(X,Y)) - \mathscr{L}_{M,\varepsilon}(\Phi(X,Y)) \\
&= \alpha \langle Q_{N,M}(g_{\varepsilon} * (:e^{\alpha Q_{N,M}(g_{\varepsilon} * X)} : e^{\alpha Q_{N,M}(g_{\varepsilon} * Y)})), \nabla_Y \Phi \rangle - \alpha \langle g_{\varepsilon} * (:e^{\alpha(g_{\varepsilon} * X)} : e^{\alpha(g_{\varepsilon} * Y)}), \nabla_Y \Phi \rangle \\
&= \mathrm{I} + \mathrm{II} + \mathrm{III}
\end{aligned}$$

with

$$\begin{aligned}\mathbb{I} &:= \alpha\langle Q_{N,M}(g_\varepsilon*(:e^{\alpha Q_{N,M}(g_\varepsilon*X)}:e^{\alpha Q_{N,M}(g_\varepsilon*Y)})), \nabla_Y\Phi\rangle - \alpha\langle g_\varepsilon*(:e^{\alpha Q_{N,M}(g_\varepsilon*X)}:e^{\alpha Q_{N,M}(g_\varepsilon*Y)}), \nabla_Y\Phi\rangle \\ &= \langle g_\varepsilon*(:e^{\alpha Q_{N,M}(g_\varepsilon*X)}:e^{\alpha Q_{N,M}(g_\varepsilon*Y)}), (Q_{N,M}-I)\operatorname{Per}_{\mathbb{T}_M^2}\nabla_Y\Phi\rangle\end{aligned} \quad (4.16)$$

$$\mathbb{II} := \alpha\langle g_\varepsilon*(:e^{\alpha Q_{N,M}(g_\varepsilon*X)}:e^{\alpha Q_{N,M}(g_\varepsilon*Y)}), \nabla_Y\Phi\rangle - \alpha\langle g_\varepsilon*(:e^{\alpha(g_\varepsilon*X)}:e^{\alpha Q_{N,M}(g_\varepsilon*Y)}), \nabla_Y\Phi\rangle \quad (4.17)$$

and

$$\mathbb{III} := \alpha\langle g_\varepsilon*(:e^{\alpha(g_\varepsilon*X)}:e^{\alpha Q_{N,M}(g_\varepsilon*Y)}), \nabla_Y\Phi\rangle - \alpha\langle g_\varepsilon*(:e^{\alpha(g_\varepsilon*X)}:e^{\alpha(g_\varepsilon*Y)}), \nabla_Y\Phi\rangle \quad (4.18)$$

Let us deal with the term $\mathbb{I}$ (4.16). We need, for $p>1$, the bound

$$\int \|:e^{\alpha Q_{N,M}(g_\varepsilon*X)}:e^{\alpha Q_{N,M}(g_\varepsilon*Y)}\|_{L^\infty}^p \mathrm{d}\mu_{N,M,\varepsilon} < C_{M,\varepsilon}. \quad (4.19)$$

Indeed, if estimate (4.19) holds, then, by the regularization property of $g_\varepsilon$ – namely the continuity for any $s>0$ of the operator $g_\varepsilon*(\cdot)\colon L^\infty(\mathbb{T}_M^2) \to B_{\infty,\infty}^s(\mathbb{T}_M^2)$ – and exploiting the fact that the norm $\|Q_{N,M}-I\|_{L(B_{\infty,\infty}^s, B_{\infty,\infty}^{s-\delta})}$ converges to zero as $N\to+\infty$ for any $s>0$ and $\delta>0$, we get the convergence of the term $\mathbb{I}$ to 0 as $N\to+\infty$. Let us prove the bound (4.19). The exponential involving $Y$ disappears, since $Y\leq 0$. Moreover, since $\varepsilon>0$, we have that the Wick exponential can be written as the exponential divided by some positive constant $C_{N,M,\varepsilon}$, i.e.

$$:e^{\alpha Q_{N,M}(g_\varepsilon*X)}:= C_{N,M,\varepsilon}^{-1} e^{\alpha Q_{N,M}(g_\varepsilon*X)},$$

which converges to some finite number as $N\to+\infty$. Then, we have the inequality

$$\|:e^{\alpha Q_{N,M}(g_\varepsilon*X)}:\|_{L^\infty} \leq C_{N,M,\varepsilon}^{-1} e^{|\alpha Q_{N,M}(g_\varepsilon*X)|_{L^\infty}},$$

and hence, by positivity of $Q_{N,M}$, we obtain

$$\|:e^{\alpha Q_{N,M}(g_\varepsilon*X)}:\|_{L^\infty} \leq C_{N,M,\varepsilon}^{-1} e^{|\alpha(g_\varepsilon*X)|_{L^\infty}},$$

which is in $L^p(\mu_M^{\text{free}})$ by Fernique's theorem (see Theorem 2.8.5 in [Bog98]).

We deal with the term $\mathbb{II}$ given by (4.17). We have

$$\mathbb{II} \leq \alpha\|g_\varepsilon\|_{L^1}\|:e^{\alpha Q_{N,M}(g_\varepsilon*X)}:-:e^{\alpha(g_\varepsilon*X)}:\|_{L^p}\|e^{\alpha Q_{N,M}(g_\varepsilon*Y)}\|_{L^\infty}\|\operatorname{Per}_{\mathbb{T}_M^2}\nabla_Y\Phi\|_{L^q},$$

and, by stochastic estimates (see Proposition B.2), we have convergence of the mean of the term $\mathbb{II}$ to zero.

For the term $\mathbb{III}$ given by (4.18) we exploit the non-positivity of $Y$ and replace the exponential by a bounded smooth function with bounded derivatives. Then we have that

$$e^{\alpha Q_{N,M}(g_\varepsilon*Y)} \to e^{\alpha(g_\varepsilon*Y)}, \qquad \text{in } L^\infty.$$

provided

$$Q_{N,M}(g_\varepsilon * Y) \to g_\varepsilon * Y, \qquad \text{in } L^\infty.$$

On the other hand, we have the bound

$$\|Q_{N,M}(g_\varepsilon * Y) - g_\varepsilon * Y\|_{L^\infty} \lesssim \|Q_{N,M} - I\|_{L(C^\delta, L^\infty)} \|g_\varepsilon * Y\|_{C^\delta}.$$

But $\|Q_{N,M} - I\|_{L(C^\delta, L^\infty)} \to 0$ depending only on $M$, while $\|g_\varepsilon * Y\|_{C^\delta} \in L^p$ uniformly since by the proof of Theorem 4.11 we have that the integral $\int V_2(X,Y) \mu_{N,M,\varepsilon}(dX, dY)$ is bounded uniformly with respect to $M$.

### 4.4.2 Limit as $M \to +\infty$

We now have to take $M \to +\infty$, namely we want to show that

$$\int \mathscr{L}_{M,\varepsilon}(\Phi) \, d\mu_{M,\varepsilon} \to \int \mathscr{L}_\varepsilon(\Phi) \, d\mu_\varepsilon, \qquad \text{as } M \to +\infty. \tag{4.20}$$

Notice that the only $M$-dependent term is now the trace-term $\operatorname{tr}_{L^2(\mathbb{T}_M^2)}(P_{\mathbb{T}_M^2} \nabla_X^2 \Phi P_{\mathbb{T}_M^2})$. By the same considerations as above, we have to show

$$\int [\mathscr{L}_{M,\varepsilon}(\Phi(X,Y)) - \mathscr{L}_\varepsilon(\Phi(X,Y))] \mu_{M,\varepsilon}(dX, dY) \to 0, \qquad \text{as } M \to +\infty.$$

Since $\Phi \in \mathscr{F}$ (see Definition 3.1 and Remark 3.2), we have

$$\operatorname{tr}_{L^2}(|\nabla_X^2 \Phi|) < f_\Phi(X), \qquad \text{and} \qquad \operatorname{tr}_{L^2}(|\rho_{-\ell} \nabla_X^2 \Phi|) < f_\Phi(X), \quad \ell > 1,$$

where $\int f_\Phi \, d\mu^{\text{free}} < +\infty$. First, note that

$$\begin{aligned}
\operatorname{tr}_{L^2}(P_M \nabla_X^2 \Phi P_M) - \operatorname{tr}_{L^2}(\nabla_X^2 \Phi) &= \operatorname{tr}_{L^2}((P_M - I)\nabla_X^2 \Phi P_M) + \operatorname{tr}_{L^2}(\nabla_X^2 \Phi (P_M - I)) \\
&= \operatorname{tr}_{L^2}((P_M - I)\nabla_X^2 \Phi P_M) + \operatorname{tr}_{L^2}(\nabla_X^2 \Phi \rho_\ell \rho_{-\ell}(P_M - I)) \\
&= \operatorname{tr}_{L^2}(\nabla_X^2 \Phi \rho_\ell \rho_{-\ell}(P_M - I)).
\end{aligned}$$

Therefore, we have the bound

$$|\operatorname{tr}_{L^2}(P_M \nabla_X^2 \Phi P_M) - \operatorname{tr}_{L^2}(\nabla_X^2 \Phi)| \lesssim \operatorname{tr}_{L^2}(|\nabla_X^2 \Phi \rho_{-\ell}|) \|\rho_\ell (P_M - I)\|_{L(L^2(\mathbb{R}^2), L^2(\mathbb{R}^2))}.$$

Now, let $h \in L^2(\mathbb{R}^2)$, then we have

$$(\rho_\ell(P_M - I))(h) = \rho_\ell \mathbb{I}_{\mathbb{R}^2 \setminus \mathbb{T}_M^2} (\operatorname{Per}_{\mathbb{T}_M^2} h - h).$$

Take $\ell > \ell' > 1$, then, considering the operator $L^2$ norm, we have

$$\|\rho_\ell(P_M - I)(h)\|_{L^2}^2 \leq 2 \sum_{y \in 2\pi \mathbb{Z}^2 \setminus \{0\}} \int_{\mathbb{T}_M^2} \rho_\ell(z - y)^2 h(z)^2 \, dz + 2 \int_{\mathbb{R}^2 \setminus \mathbb{T}_M^2} \rho_\ell(z)^2 h(z)^2 \, dz.$$

Since $\rho_\ell(z-y) = \rho_{\ell'}(z-y)\rho_{\ell-\ell'}(z-y) \leq \rho_{\ell'}(z-y)(1+M)^{-(\ell-\ell')}$, we get

$$\|\rho_\ell(P_M - I)(h)\|_{L^2}^2$$
$$\lesssim 2(1+M)^{-2(\ell-\ell')} \sum_{y \in 2\pi\mathbb{Z}^2 \setminus \{0\}} \int_{\mathbb{T}_M^2} \rho_{\ell'}(z-y)^2 h(z)^2 \,\mathrm{d}z + 2(1+M)^{-2\ell} \int_{\mathbb{R}^2 \setminus \mathbb{T}_M^2} h(z)^2 \,\mathrm{d}z.$$

Now note

$$\sum_{y \in 2\pi\mathbb{Z}^2 \setminus \{0\}} \int_{\mathbb{T}_M^2} \rho_{\ell'}(z-y)^2 h(z)^2 \,\mathrm{d}z \lesssim \int_{\mathbb{T}_M^2} (\mathrm{Per}_{\mathbb{T}_M^2}(\rho_{\ell'})(z))^2 h(z)^2 \,\mathrm{d}z,$$

and $\mathrm{Per}_{\mathbb{T}_M^2}(\rho_{\ell'})(z) \leq C_{\ell'} := \sum_{y \in \mathbb{Z}^2} \rho_{\ell'}(y)$, to get

$$\|\rho_\ell(P_M - I)(h)\|_{L^2}^2 \lesssim 2(1+M)^{-2(\ell-\ell')} C_{\ell'} \int_{\mathbb{T}_M^2} h(z)^2 \,\mathrm{d}z + 2(1+M)^{-2\ell} \int_{\mathbb{R}^2 \setminus \mathbb{T}_M^2} h(z)^2 \,\mathrm{d}z$$
$$\lesssim 2(1+M)^{-2(\ell-\ell')}(C_{\ell'}+1) \int_{\mathbb{T}_M^2} h(z)^2 \,\mathrm{d}z.$$

By taking the supremum over all $h$ with $\|h\|_{L^2(\mathbb{R}^2)} \leq 1$, we then obtain

$$\|\rho_\ell(P_M - I)\|_{L(L^2(\mathbb{R}^2), L^2(\mathbb{R}^2))} \lesssim (1+M)^{-2(\ell-\ell')}.$$

Letting $M \to +\infty$, we get the limit (4.20).

### 4.4.3 Limit as $\varepsilon \to 0$

Now, we are left to show the convergence

$$\int \mathcal{L}_\varepsilon(\Phi) \,\mathrm{d}\mu_\varepsilon \to \int \mathcal{L}(\Phi) \,\mathrm{d}\mu, \qquad \text{as } \varepsilon \to 0.$$

As above, we have to show that

$$\int [\mathcal{L}_\varepsilon(\Phi(X,Y)) - \mathcal{L}(\Phi(X,Y))] \mu_\varepsilon(\mathrm{d}X, \mathrm{d}Y) \to 0, \qquad \text{as } \varepsilon \to 0.$$

Rewriting the integrand, we get

$$\mathcal{L}_\varepsilon(\Phi(X,Y)) - \mathcal{L}(\Phi(X,Y)) = \alpha \langle (g_\varepsilon * (:e^{\alpha(g_\varepsilon * X)}: e^{\alpha(g_\varepsilon * Y)})), \nabla_Y \Phi \rangle - \alpha \langle :e^{\alpha X}: e^{\alpha Y}, \nabla_Y \Phi \rangle$$
$$= \mathrm{I}' + \mathrm{II}' + \mathrm{III}' \tag{4.21}$$

with

$$\mathrm{I}' := \alpha \langle (g_\varepsilon * (:e^{\alpha(g_\varepsilon * X)}: e^{\alpha(g_\varepsilon * Y)})), \nabla_Y \Phi \rangle - \alpha \langle (:e^{\alpha(g_\varepsilon * X)}: e^{\alpha(g_\varepsilon * Y)}), \nabla_Y \Phi \rangle \tag{4.22}$$

$$\mathrm{II}' := \alpha \langle (:e^{\alpha(g_\varepsilon * X)}: e^{\alpha(g_\varepsilon * Y)}), \nabla_Y \Phi \rangle - \alpha \langle (:e^{\alpha X}: e^{\alpha(g_\varepsilon * Y)}), \nabla_Y \Phi \rangle \tag{4.23}$$

$$\mathrm{III}' := \alpha \langle :e^{\alpha X}: e^{\alpha(g_\varepsilon * Y)}), \nabla_Y \Phi \rangle - \alpha \langle (:e^{\alpha X}: e^{\alpha Y}), \nabla_Y \Phi \rangle \tag{4.24}$$

The stochastic estimates for $:\exp(\alpha(g_\varepsilon * X)):$ are done in Proposition B.1. We now deal with the term $\mathbb{I}'$ (4.22). We have the inequality

$$\begin{aligned} \mathbb{I}' &= \alpha\langle :e^{\alpha(g_\varepsilon * X)}: e^{\alpha(g_\varepsilon * Y)}, (g_\varepsilon - I) * \nabla_Y \Phi\rangle \\ &\leq \alpha \|:e^{\alpha(g_\varepsilon * X)}: e^{\alpha(g_\varepsilon * Y)}\|_{B^s_{p,p,\ell}} \|(g_\varepsilon - I) * \nabla_Y \Phi\|_{B^{-s}_{q,q,-\ell}} \\ &\lesssim \alpha \|:e^{\alpha(g_\varepsilon * X)}:\|_{B^s_{p,p,\ell}} \|e^{\alpha(g_\varepsilon * Y)}\|_{L^\infty} \|g_\varepsilon - I\|_{L(B^{-s+\delta}_{q,q,-\ell}, B^{-s}_{q,q,-\ell})} \|\nabla_Y \Phi\|_{B^{-s+\delta}_{q,q,-\ell}}. \end{aligned}$$

Taking the expectation and exploiting the negativity of $Y$, we have

$$\begin{aligned} &\mathbb{E}[\alpha\langle (g_\varepsilon * (:e^{\alpha(g_\varepsilon * X)}: e^{\alpha(g_\varepsilon * Y)})) - :e^{\alpha(g_\varepsilon * X)}: e^{\alpha(g_\varepsilon * Y)}, \nabla_Y \Phi\rangle] \\ &\lesssim \mathbb{E}[\|:e^{\alpha(g_\varepsilon * X)}:\|_{B^s_{p,p,\ell}}] \|g_\varepsilon - I\|_{L(B^{-s+\delta}_{q,q,-\ell}, B^{-s}_{q,q,-\ell})} \|\nabla_Y \Phi\|_{L^\infty(B_X \times B_Y, B^{-s+\delta}_{q,q,-\ell})} \\ &\lesssim \|(g_\varepsilon - 1)\|_{L(B^{-s+\delta}_{q,q,-\ell}, B^{-s}_{q,q,-\ell})}, \end{aligned}$$

and this last term converges to zero as $\varepsilon \to 0$. The convergence to zero of the term $\mathbb{II}'$ (4.23) follows from Proposition B.1. Finally for the term $\mathbb{III}'$ (4.24), we proceed as follows. We have the bound

$$\mathbb{III}' \lesssim \|:e^{\alpha X}:\|_{B^s_{p,p,\ell}} \left( \left( \|Y\|_{B^{-s+\delta}_{q,q,\ell}} \|g_\varepsilon - I\|_{L(B^{-s+\delta}_{q,q,-\ell}, B^{-s}_{q,q,-\ell})} \right) \wedge 2 \right) \|\nabla_Y \Phi\|_{B^{-s+\delta}_{q,q,-\ell}}. \tag{4.25}$$

Since $\mu_\varepsilon$ is tight and $\|:\exp(\alpha X):\|_{B^s_{p,p,\ell}}$ is uniformly integrable with respect to the measure $\mu_\varepsilon$, then, for any $\tau > 0$, there exists a Borel subset $\Omega_\tau \subset B_X \times B_Y$ such that

$$\int_{\Omega_\tau^c} \|:e^{\alpha X}:\|_{B^s_{p,p,\ell}} d\mu_\varepsilon < \tau,$$

and $R_{\Omega_\tau} = \sup_{\Omega_\tau} \|Y\|_{B^{-s+\delta}_{q,q,\ell}} < +\infty$. Therefore, we have

$$\begin{aligned} \int |\mathbb{III}'| d\mu_\varepsilon &\leq \int_{\Omega_\tau} |\alpha\langle :e^{\alpha X}: (e^{\alpha(g_\varepsilon * Y)} - e^{\alpha Y}), \nabla_Y \Phi\rangle| d\mu_\varepsilon + \int_{\Omega_\tau^c} |\alpha\langle :e^{\alpha X}: (e^{\alpha(g_\varepsilon * Y)} - e^{\alpha Y}), \nabla_Y \Phi\rangle| d\mu_\varepsilon \\ &\leq R_{\Omega_\tau} \|\nabla_Y \Phi\|_{B^{-s+\delta}_{q,q,-\ell}} \|g_\varepsilon - I\|_{L(B^{-s+\delta}_{q,q,-\ell}, B^{-s}_{q,q,-\ell})} \int_\Omega \|:e^{\alpha X}:\|_{B^s_{p,p,\ell}} d\mu + 2\tau \|\nabla_Y \Phi\|_{B^{-s+\delta}_{q,q,-\ell}}. \end{aligned}$$

Thus,

$$\lim_{\varepsilon \to 0} \int |\mathbb{III}'| d\mu_\varepsilon \leq 2\|\nabla_Y \Phi\|_{B^{-s+\delta}_{q,q,-\ell}} \tau,$$

which gives convergence to zero by arbitrary choice of $\delta$. In order to show inequality (4.25), note that

$$\begin{aligned} |\mathbb{III}'| &\leq \alpha \|:\exp(\alpha X): (\exp(\alpha(g_\varepsilon * Y)) - \exp(\alpha Y))\|_{B^{s-\delta}_{l,l,2\ell}} \|\nabla_Y \Phi\|_{B^{-s+\delta}_{l,l,-\ell}} \\ &\lesssim \alpha \|:\exp(\alpha X): (\exp(\alpha(g_\varepsilon * Y)) - \exp(\alpha Y))\|_{B^s_{1,1,2\ell}} \|\nabla_Y \Phi\|_{B^{-s+\delta}_{l,l,-\ell}} \\ &\lesssim \alpha [(\|:\exp(\alpha X):\|_{B^s_{p,p,\ell}} \|(g_\varepsilon * Y) - Y\|_{B^{-s}_{q,q,\ell}}) \wedge (2 \|:\exp(\alpha X):\|_{B^s_{1,1,\ell}})] \|\nabla_Y \Phi\|_{B^{-s+\delta}_{l,l,-\ell}} \\ &\lesssim (\|:\exp(\alpha X):\|_{B^s_{p,p,\ell}} + \|:\exp(\alpha X):\|_{B^s_{1,1,\ell}})(\|(g_\varepsilon * Y) - Y\|_{B^{-s}_{q,q,\ell}} \wedge 2) \|\nabla_Y \Phi\|_{B^{-s+\delta}_{l,l,-\ell}}, \end{aligned}$$

where $l \in (1, +\infty)$ and $1/l' + 1/l = 1$ such that $\delta > 2 - 2/l'$. Now, if $|s|$ is small enough such that $B_{q,q,\ell}^{-s+\delta} \subset B_Y$, then we have the

$$\left(\|(g_\varepsilon * Y) - Y\|_{B_{q,q,\ell}^{-s}} \wedge 2\right) \lesssim \left(\|Y\|_{B_{q,q,\ell}^{-s+\delta}} \|g_\varepsilon - I\|_{L(B_{q,q,-\ell}^{-s+\delta}, B_{q,q,-\ell}^{-s})}\right) \wedge 2,$$

which is converging to zero for every $Y \in B_{q,q,\ell}^{-s+\delta} \subset B_Y$. We conclude the argument by Lebesgue dominated convergence theorem and with the reasonings similar as above.

# Appendix A  Besov spaces and heat semigroup

In this section, we collect some results about weighted Besov spaces. While we only focus on spaces defined on the whole space $\mathbb{R}^n$, the results hold also for Besov spaces on the $n$-dimensional torus $\mathbb{T}^n$.

Let us start by introducing Littlewood-Paley blocks. Let $\chi$ and $\varphi$ be smooth non-negative functions from $\mathbb{R}^n$ into $\mathbb{R}$ satisfying the following properties:

- $\mathrm{supp}(\chi) \subset B_{4/3}(0)$ and $\mathrm{supp}(\varphi) \subset B_{8/3}(0) \setminus B_{3/4}(0)$,

- $\chi, \varphi \leq 1$ and $\chi(y) + \sum_{j \geq 0} \varphi(2^{-j}y) = 1$, for any $y \in \mathbb{R}^n$,

- $\mathrm{supp}(\chi) \cap \mathrm{supp}(\varphi(2^{-j}\cdot)) = \emptyset$, for $j \geq 1$,

- $\mathrm{supp}(\varphi(2^{-j}\cdot)) \cap \mathrm{supp}(\varphi(2^{-i}\cdot)) = \emptyset$, if $|i - j| > 1$,

where $B_r(x)$ denotes the ball centered at $x \in \mathbb{R}^n$ with radius $r > 0$.

We introduce the following notations: $\varphi_{-1} = \chi$, $\varphi_j(\cdot) = \varphi(2^{-j}\cdot)$, $D_j = \hat{\varphi}_j$, and for any $f \in \mathcal{S}'(\mathbb{R}^n)$ we put $\Delta_j(f) = D_j * f$. Moreover, we write, for any $\ell > 0$, $\rho_\ell(y) = (1 + |y|^2)^{-\ell/2}$, and let $L_\ell^p(\mathbb{R}^n)$ be the $L^p$-space with respect to the norm

$$\|f\|_{L_\ell^p} = \left(\int_{\mathbb{R}^n} (f(y)\rho_\ell(y))^p \, dy\right)^{1/p},$$

where $p \in [1, +\infty]$.

**Definition A.1. (Besov space $B_{p,q,\ell}^s$)** *Let $s \in \mathbb{R}$, $p, q \in [1, +\infty]$, and $\ell \in \mathbb{R}$. For $f \in \mathcal{S}'(\mathbb{R}^n)$, we define the norm*

$$\|f\|_{B_{p,q,\ell}^s} = \left(\sum_{j \geq -1} 2^{sqj} \|\Delta_j(f)\|_{L_\ell^p}^q\right)^{1/q}.$$

*The space $B_{p,q,\ell}^s(\mathbb{R}^n)$ is the subset of $\mathcal{S}'(\mathbb{R}^n)$ such that the norm $\|\cdot\|_{B_{p,q,\ell}^s}$ is finite.*

In the case where $p = q = +\infty$, the weighted Besov space $B_{\infty,\infty,\ell}^s(\mathbb{R}^n)$ is denoted by $C_\ell^s(\mathbb{R}^d)$ and it is called *weighted Besov-Hölder space* with regularity $s$. Moreover, if $s \in \mathbb{R}_+ \setminus \mathbb{Z}$, the space $C_\ell^s(\mathbb{R}^d)$ coincides with the Banach space of $s$-Hölder-continuous functions.

The relation between weighted Besov spaces is stated in the following result.

**Proposition A.2. (Besov embedding)** *Let $p_1, p_2, q_1, q_2 \in [1, +\infty]$, $\ell_1, \ell_2 \in \mathbb{R}$, and $s_1, s_2 \in \mathbb{R}$ be such that $s_1 - \frac{n}{p_1} > s_2 - \frac{n}{p_2}$ and $\ell_1 > \ell_2$. Then, we have the compact immersion*

$$B^{s_2}_{p_2, q_2, \ell_2} \subset B^{s_1}_{p_1, q_1, \ell_1}.$$

**Proof.** The proof can be found in Theorem 6.7 in [Tri06]. □

If $p = q = 2$, the Besov space $B^s_{2,2,\ell}$ coincides with the Sobolev space $H^s_\ell$, i.e. the space of measurable tempered distributions $f$ with bounded norm

$$\|f\|^2_{H^s_\ell} = \int_{\mathbb{R}^n} \rho^2_\ell(y) ((-\Delta + 1)^{s/2} f)^2(y) \, \mathrm{d}y.$$

The following theorem allows to extend products between a smooth function and a distribution to elements of Besov spaces.

**Theorem A.3. (Paraproduct)** *Let $p_1, p_2, p, q_1, q_2, q \in [1, +\infty]$, $\ell_1, \ell_2, \ell, s_1, s_2, s \in \mathbb{R}$, be such that*

$$\frac{1}{p} = \frac{1}{p_1} + \frac{1}{p_2}, \qquad \frac{1}{q} = \frac{1}{q_1} + \frac{1}{q_2}, \qquad \ell = \ell_1 + \ell_2, \qquad s_1 + s_2 > 0, \qquad s = s_1 \wedge s_2. \tag{A.1}$$

*Consider the bilinear map $\Pi : \mathcal{S}'(\mathbb{R}^n) \times \mathcal{S}(\mathbb{R}^n) \to \mathcal{S}'(\mathbb{R}^n)$, mapping $(f, g) \mapsto \Pi(f, g) = f \cdot g$. Then, there exists a unique continuous extension of $\Pi$ as the map*

$$\Pi : B^{s_1}_{p_1, q_1, \ell_1}(\mathbb{R}^n) \times B^{s_2}_{p_2, q_2, \ell_2}(\mathbb{R}^n) \to B^s_{p, q, \ell}(\mathbb{R}^n),$$

*and we have, for any $f \in B^{s_1}_{p_1, q_1, \ell_1}(\mathbb{R}^n)$, $g \in B^{s_2}_{p_2, q_2, \ell_2}(\mathbb{R}^n)$,*

$$\|\Pi(f, g)\|_{B^s_{p, q, \ell}} \lesssim \|f\|_{B^{s_1}_{p_1, q_1, \ell_1}} \|g\|_{B^{s_2}_{p_2, q_2, \ell_2}}.$$

**Proof.** See Section 3.3 in [MW17] for Besov spaces with exponential weights. The proof for polynomial weights follows in a similar way. □

In the case of products between a positive measure and an element of a Besov space, the previous result can be improved as follows.

**Proposition A.4.** *Consider the same parameters as in Theorem A.3 satisfying (A.1) and $s_1 > 0$, $s_2 \leq 0$. Suppose that $f \in B^{s_1}_{p_1, q_1, \ell_1}(\mathbb{R}^2) \cap L^\infty(\mathbb{R}^2)$ and that $\mu \in B^{s_2}_{p_2, q_2, \ell_2}(\mathbb{R}^2)$ is a positive measure, then we have*

$$\|f \cdot \mu\|_{B^{s_2}_{p_2, q_2, \ell_2}} \lesssim \|f\|_{L^\infty} \cdot \|\mu\|_{B^{s_2}_{p_2, q_2, \ell_2}}.$$

**Proof.** See, e.g., Lemma 28 in [ADG21]. □

The next result is an interpolation estimate for Besov spaces.

**Proposition A.5.** *Consider $p_1, p_2, q_1, q_2 \in [1, +\infty]$, $\ell_1, \ell_2 \in \mathbb{R}$ and $s_1, s_2 \in \mathbb{R}$, and write, for any $\theta \in [0,1]$,*

$$\frac{1}{p_\theta} = \frac{\theta}{p_1} + \frac{1-\theta}{p_2}, \qquad \frac{1}{q_\theta} = \frac{\theta}{q_1} + \frac{1-\theta}{q_2}, \qquad \ell_\theta = \theta \ell_1 + (1-\theta)\ell_2, \qquad s_\theta = \theta s_1 + (1-\theta)s_2.$$

*If $f \in B^{s_1}_{p_1,q_1,\ell_1}(\mathbb{R}^n) \cap B^{s_2}_{p_2,q_2,\ell_2}(\mathbb{R}^n)$, then $f \in B^{s_\theta}_{p_\theta,q_\theta,\ell_\theta}(\mathbb{R}^n)$, and furthermore*

$$\|f\|_{B^{s_\theta}_{p_\theta,q_\theta,\ell_\theta}} \leq \|f\|^\theta_{B^{s_1}_{p_1,q_1,\ell_1}} \|f\|^{1-\theta}_{B^{s_2}_{p_2,q_2,\ell_2}}.$$

**Proof.** The proof is based on the fact that the complex interpolation of the two spaces $B^{s_1}_{p_1,q_1,\ell_1}(\mathbb{R}^n)$ and $B^{s_2}_{p_2,q_2,\ell_2}(\mathbb{R}^n)$ is given by $B^{s_\theta}_{p_\theta,q_\theta,\ell_\theta}(\mathbb{R}^n)$. Such an interpolation is shown in Theorem 6.4.5 in [BL76] for unweighted Besov spaces. The proof for weighted spaces follows from the fact that $f \in B^s_{p,q,\ell}$ if and only if $f \cdot \rho_\ell \in B^s_{p,q}$ (see Theorem 6.5 in [Tri06]). □

We introduce now the heat kernel and present some of its properties. Let $P_t = e^{-(-\Delta+m^2)t}$, we consider $f \in L^r_{\ell_1}(\mathbb{R}, B^s_{p,q,\ell_2}(\mathbb{R}^2))$ and define the *heat kernel* on $f$ as

$$e^{-(-\Delta+m^2)}f(t) = \int_{-\infty}^{t} P_{t-\tau}f(\tau)\,\mathrm{d}\tau.$$

Notice that $e^{-(-\Delta+m^2)}f(t)$ is a distribution.

More precisely, if $s > 0$ and if $f(t, x)$, $(t, x) \in \mathbb{R} \times \mathbb{R}^2$, is a measurable function then

$$e^{-(-\Delta+m^2)}f(t,x) = \int_{-\infty}^{t}\int_{\mathbb{R}^2} \frac{1}{4\pi(t-\tau)} e^{-\frac{(x-y)^2}{4(t-\tau)} + m^2(t-\tau)} f(\tau,y)\,\mathrm{d}y\mathrm{d}\tau.$$

We have the following regularization property for $e^{-(-\Delta+m^2)}$. Let us remark that all the following results hold also in the case where $f \in L^r_{\ell_1}([t_1, t_2], B^s_{p,q,\ell_2}(\mathbb{R}^2))$, where $-\infty \leq t_1 < t_2 \leq +\infty$, and the operator $e^{-(-\Delta+m^2)}f(t)$ is defined as the integral from $t_1$ to $t \in [t_1, t_2]$.

**Theorem A.6.** *Consider $r \in [1, +\infty]$, $p, q \in [1, +\infty]$, $s \in \mathbb{R}$, and let $f \in L^r_{\ell_1}(\mathbb{R}, B^s_{p,q,\ell_2}(\mathbb{R}^2))$. Then, for any $\beta_1, \beta_2 > 0$ such that $\beta_1 + \beta_2 < 1$, we have*

$$e^{-(-\Delta+m^2)}f \in B^{\beta_2}_{r,r,\ell_1}(\mathbb{R}, B^{s+2\beta_1}_{p,q,\ell_2}(\mathbb{R}^2)). \tag{A.2}$$

Notice that (A.2) states that we are gaining regularity $\beta_2$ in time and $2\beta_1$ in space.

In order to prove Theorem A.6, we need the following result saying that when we apply the heat kernel at time $t$ we gain $2\beta_1$ in space-regularity, but we have to pay with a multiplicative factor of $t^{-\beta_1}$.

**Lemma A.7.** *Let $m > 0$ and consider $g \in B_{p,q,\ell}^s(\mathbb{R}^2)$. We have, for every $t > 0$,*

$$\|P_t g\|_{B_{p,q,\ell}^{s+2\beta_1}(\mathbb{R}^2)} \lesssim t^{-\beta_1} e^{-m^2 t} \|g\|_{B_{p,q,\ell}^s}.$$

**Proof.** See Proposition 5 in [MW17]. □

We will also need the following lemma saying that giving up some space-regularity we can gain a factor $t^\beta$ on the right-hand side.

**Lemma A.8.** *Consider $0 < \beta < 1$ and $g \in B_{p,q,\ell}^s(\mathbb{R}^2)$. Then, for any $t > 0$,*

$$\|(1 - P_t)g\|_{B_{p,q,\ell}^{s-2\beta}} \lesssim t^\beta \|g\|_{B_{p,q,\ell}^s}.$$

**Proof.** See Proposition 6 in [MW17]. □

**Proof of Theorem A.6.** We give the proof for unweighted Besov spaces, the general case follows the same lines. We use the difference characterization of space-time Besov spaces (see e.g. Theorem 2.36 in [BCD11] or Chapter 2.6.1 in [Tri92]), which yields

$$\begin{aligned}\|e^{-(-\Delta-m^2)}f\|^r_{B_{r,r}^{\beta_2}(\mathbb{R}, B_{p,q}^{s+2\beta_1}(\mathbb{R}^2))} &\sim \|e^{-(-\Delta-m^2)}f\|^r_{L^r(\mathbb{R}, B_{p,q}^{s+2\beta_1}(\mathbb{R}^2))} \\ &\quad + \int_\mathbb{R} \int_{|\Delta t| \leq 1} \frac{\|e^{-(-\Delta-m^2)}f(t+\Delta t) - e^{-(-\Delta-m^2)}f(t)\|^r_{B_{p,q}^s}}{|\Delta t|^{1+r\beta_2}} \, d(\Delta t) \, dt.\end{aligned}$$

First, we prove that the first term on the right-hand side is finite. Write $\tilde{f} = e^{-(-\Delta+m^2)}f$, we have

$$\begin{aligned}\|\tilde{f}\|^r_{L^r(\mathbb{R}, B_{p,q}^{s+2\beta_1}(\mathbb{R}^2))} &= \int_\mathbb{R} \|\tilde{f}(t)\|^r_{B_{p,q}^{s+2\beta_1}} \, dt \\ &= \int_\mathbb{R} \left\|\int_{-\infty}^t P_{t-k} f(k) \, dk\right\|^r_{B_{p,q}^{s+2\beta_1}} \, dt \\ &\lesssim \int_\mathbb{R} \left(\int_{-\infty}^t e^{-m^2(t-k)} \|e^{\Delta(t-k)} f(k)\|_{B_{p,q}^{s+2\beta_1}} \, dk\right)^r \, dt\end{aligned}$$

Lemma A.7 and Young inequality yield

$$\begin{aligned}\|\tilde{f}\|^r_{L^r(\mathbb{R}, B_{p,q}^{s+2\beta_1}(\mathbb{R}^2))} &\lesssim \int_\mathbb{R} \left(\int_\mathbb{R} \frac{\mathbb{I}_{[0,+\infty]}(t-k)}{(t-k)^{\beta_1}} e^{-m^2(t-k)} \|f(k)\|_{B_{p,q}^s} \, dk\right)^r \, dt \\ &\lesssim \left(\int_\mathbb{R} \frac{\mathbb{I}_{[0,+\infty]}(t-k)}{(t-k)^{\beta_1}} e^{-m^2(t-k)} \, dk\right)^r + \|f\|^r_{L^r(\mathbb{R}, B_{p,q}^s(\mathbb{R}^n))},\end{aligned}$$

where the first integral on the last step is finite if and only if $\beta_1 < 1$.

Consider now the difference term. We have

$$\begin{aligned}\tilde{f}(t + \Delta t) - \tilde{f}(t) &= \int_t^{t+\Delta t} P_{t+\Delta t - k} f(k) \, dk + (1 - e^{-(-\Delta+m^2)\Delta t}) \int_{-\infty}^t P_{t-k} f(k) \, dk \\ &=: I_1 + I_2.\end{aligned}$$

Now, by Lemma A.7 and Young inequality,

$$\left\|\left(\|I_1\|_{B_{p,q}^{s+2\beta_1}}\right)\right\|_{L^r}^r \lesssim \int_{\mathbb{R}}\left(\int_{\mathbb{R}} \frac{\mathbb{I}_{[-\Delta t,0]}(t-k)}{(t+\Delta t-k)^{\beta_1}}\|f(k)\|_{B_{p,q}^s}\right)^r$$

$$\lesssim (\Delta t)^{1-\beta_1}\|f\|_{L^r(\mathbb{R},B_{p,q}^s(\mathbb{R}^n))}$$

Consider $\beta_2 < \tilde{\beta} < 1 - \beta_1$, then, by Lemma A.8,

$$\left\|\left(\|I_2\|_{B_{p,q}^{s+2\beta_1}}\right)\right\|_{L^r}^r \lesssim (\Delta t)^{\tilde{\beta}}\|e^{-(-\Delta+m^2)}f\|_{L^r(\mathbb{R},B^{s+2\beta_1+2\tilde{\beta}})}$$

$$\lesssim (\Delta t)^{\tilde{\beta}}\|f\|_{L^r(\mathbb{R},B_{p,q}^s)},$$

where we used the first part of the proof in the last step and the fact that $\beta_1 + \tilde{\beta} < 1$. Putting everything together, we get

$$\int_{\mathbb{R}}\int_{|\Delta t|\leq 1} \frac{\|\tilde{f}(t+\Delta t)-\tilde{f}(t)\|_{B_{p,q}^s}^r}{|\Delta t|^{1+r\beta_2}} d(\Delta t) dt \lesssim \|f\|_{L^1(\mathbb{R},B_{p,q}^s)}\int_{|\Delta t|\leq 1} \frac{(\Delta t)^{(1-\beta_1)r}+(\Delta t)^{\tilde{\beta}r}}{|\Delta t|^{1+r\beta_2}} d(\Delta t)$$

$$\lesssim \|f\|_{L^1(\mathbb{R},B_{p,q}^s)}\int_{|\Delta t|\leq 1} \frac{1}{|\Delta t|^{1-(\tilde{\beta}-\beta_2)r}} d(\Delta t)$$

$$\lesssim \|f\|_{L^r(\mathbb{R},B_{p,q}^s)},$$

which gives the result. □

Thanks to the heat kernel, we have another representation for weighted Besov spaces.

**Proposition A.9.** *Let $s \in \mathbb{R}$, $p, q \in (0, +\infty]$, and $k \in \mathbb{N}_0$ be such that*

$$k > \frac{s}{2}.$$

*Consider a smooth and compactly supported function $\varphi_0$. Then, for any $m > 0$, we have the following equivalence between norms*

$$\|f\|_{B_{p,q,\ell}^s} \simeq \|\mathscr{F}^{-1}(\varphi_0\mathscr{F}(f))\|_{L_\ell^p} + \left(\int_0^{+\infty} t^{(k-\frac{s}{2})q}\|\partial_{t^k}(P_t f)\|_{L_\ell^p}^q \frac{dt}{t}\right)^{1/q}. \tag{A.3}$$

**Proof.** See Section 2.6.4 in [Tri92] for a version of the theorem with a mass-less heat kernel. In particular, such a result differs from the one presented above since the integral with respect to $t$ in equation (A.3) goes from 0 to 1 instead of going from zero to $+\infty$. This extension is possible thanks to the regularization properties of the heat kernel and the exponential decay thereof, due to the presence of the positive mass $m > 0$. □

**Remark A.10.** Under the same hypotheses as in Proposition A.9, if we assume further that $s>0$, then the norm $\|\mathcal{F}^{-1}(\varphi_0\mathcal{F}(f))\|_{L^p_\ell}$ appearing in equation (A.3) can be substituted by the weighted $L^p$-norm of $f$, since we have

$$\|\mathcal{F}^{-1}(\varphi_0\mathcal{F}(f))\|_{L^p_\ell} \lesssim \|f\|_{L^p_\ell} \lesssim \|f\|_{B^s_{p,q,\ell}}.$$

Namely, for $s>0$, equivalence (A.3) becomes

$$\|f\|_{B^s_{p,q,\ell}} \simeq \|f\|_{L^p_\ell} + \left(\int_0^1 t^{(k-\frac{s}{2})q} \|\partial_{t^k}(P_t f)\|^q_{L^p_\ell} \frac{\mathrm{d}t}{t}\right)^{1/q}.$$

# Appendix B  Stochastic estimates for the Wick exponential

We prove here some stochastic estimates for the term $:e^{\alpha g_\varepsilon * X}:$. We need two different kind of estimates. One has to be at the initial time with respect to the Gaussian free field, while the second one needs to estimate the term in some $L^p$-space with respect to the variable $t$. Given $g_\varepsilon$ as in Section 2.2 and $m>0$, and $\mu^{\text{free}}$ being the Gaussian free field with mass $m$, we define the Wick exponential $:\exp(\alpha g_\varepsilon * X):$ as in equation (2.13).

The previous expression coincides with the standard Wick exponential in the case where we equip the space $B_X$ with the free field measure $\mu^{\text{free}}$.

**Proposition B.1.** *Let $\alpha^2 < 8\pi$, $\varepsilon > 0$. Then, for every $r>1$ such that $\alpha^2 r/(4\pi) < 2$, and for every $\ell>0$ such that $r\ell > 2$, and for every $\delta > 0$, we have that the sequence*

$$:\exp(\alpha g_\varepsilon * X): \rightarrow :\exp(\alpha X):, \qquad \text{in } L^r((B_X, \mu^{\text{free}}), B^{-\gamma(r-1)-\delta}_{r,r,\ell}(\mathbb{R}^2)),$$

*where $:\exp(\alpha X):$ is the unique (positive) limit random distribution and $\gamma = \alpha^2/(4\pi)$.*

**Proof.** We report here the proof for the case $\alpha^2 < 4\pi$, the general case can be obtained mixing the method presented here and the techniques by [HKK22]. We consider first the case $r=2$. We have, for $K_j = \mathcal{F}^{-1}(\varphi_j)$,

$$\mathbb{E}\left[\int_{\mathbb{R}^2} |[(:e^{\alpha(g_\varepsilon * X)}:-:e^{\alpha X}:) * K_j](z)|^2 (\rho_\ell(z))^2 \mathrm{d}z\right] = \int_{\mathbb{R}^2} \mathbb{E}[|[(:e^{\alpha(g_\varepsilon * X)}:-:e^{\alpha X}:) * K_j](z)|^2](\rho_\ell(z))^2 \mathrm{d}z.$$

By translation invariance and orthogonality of Wick polynomials, we can consider

$$\begin{aligned}\mathbb{E}[|[(:e^{\alpha(g_\varepsilon * X)}:-:e^{\alpha X}:) * K_j](0)|^2] &= \mathbb{E}[|\langle(:e^{\alpha(g_\varepsilon * X)}:-:e^{\alpha X}:), K_j\rangle|^2] \\ &= \sum_{n=0}^{+\infty} \frac{\alpha^{2n}}{(n!)^2} \mathbb{E}[|\langle :(g_\varepsilon * X)^n: - :X^n:, K_j\rangle|^2].\end{aligned}$$

It then suffices to show that, for each $n \in \mathbb{N}$,

$$\mathbb{E}[|\langle :(g_\varepsilon * X)^n: - :X^n:, K_j\rangle|^2] \to 0. \tag{B.1}$$

Indeed, for $|\alpha| < \sqrt{4\pi}$, we have

$$\sum_{n=0}^{+\infty} \frac{\alpha^{2n}}{(n!)^2} \mathbb{E}[|\langle :(g_\varepsilon * X)^n: - :X^n:, K_j\rangle|^2] \leq 2\sum_{n=0}^{+\infty} \frac{\alpha^{2n}}{(n!)^2} \mathbb{E}[|\langle :(g_\varepsilon * X)^n:, K_j\rangle|^2 + |\langle :X^n:, K_j\rangle|^2].$$

We have then to show
$$\left| \frac{\alpha^{2n}}{(n!)^2} \mathbb{E}[|\langle :(g_\varepsilon * X)^n:, K_j \rangle|^2 + |\langle :(X)^n:, K_j \rangle|^2] \right| \le c_n,$$
for some $\{c_n\} \in \ell^1(\mathbb{N})$.

We want an uniform bound on $\|\Delta_j :(g_\varepsilon * X)^n: \|_{L^2}$, which would then give, as $\varepsilon \to 0$, the convergence $\Delta_j :(g_\varepsilon * X)^n: \to \Delta_j :X^n:$ almost surely.

Notice that
$$\|\Delta_j :(g_\varepsilon * X)^n: \|_{L^2}^2 = n! \int_{\mathbb{T}_M^2} K_j(x) K_j(x') (\mathcal{X}_\varepsilon(x-x'))^n \, dx \, dx',$$
where
$$\mathcal{X}_\varepsilon = g_\varepsilon^{*2} * \mathcal{X} = \mathcal{F}^{-1}\left( \frac{|\hat{g}_\varepsilon|^2}{|\cdot|^2 + m^2} \right).$$

By Lemma 2.2 in [ORW21], we have, for some constant $C > 0$,
$$\mathcal{X}_\varepsilon(x) \lesssim -\frac{1}{2\pi} \log(|x| \wedge \varepsilon) + C,$$
$$\mathcal{X}(x) \lesssim -\frac{1}{2\pi} \log |x| + C.$$

Moreover, if $x \ne 0$, we have the point-wise limits
$$\mathcal{X}_\varepsilon(x) \xrightarrow{\varepsilon \to 0} \mathcal{X}(x).$$

Now,
$$\|\Delta_j :(g_\varepsilon * X)^n: \|_{L^2}^2 = n! \int_{\mathbb{R}^2} K_j(x) K_j(x') (\mathcal{X}_\varepsilon(x-x'))^n \, dx \, dx'$$
$$\lesssim n! \left( \frac{1}{2\pi} \log |2^j| \right)^n,$$
$$\|\Delta_j :X^n: \|_{L^2}^2 = n! \int_{\mathbb{R}^2} K_j(x) K_j(x') (\mathcal{X}(x-x'))^n \, dx \, dx'$$
$$\lesssim n! \left( \frac{1}{2\pi} \log |2^j| \right)^n,$$
where the multiplicative constant absorbed in the symbol $\lesssim$ does not depend on $j$, $\varepsilon$, and $n$.

Summing up, we have
$$\left| \frac{\alpha^{2n}}{(n!)^2} \mathbb{E}[|\langle :(X^{N,M})^n:, K_j \rangle|^2 + |\langle :(X^M)^n:, K_j \rangle|^2] \right| \lesssim \left| \frac{\alpha^{2n}}{n!} \left( \frac{1}{2\pi} \log |2^j| \right)^n \right| =: c_n^j.$$

Then,
$$c^j = \sum_n c_n^j \lesssim \sum_n \left| \frac{\alpha^{2n}}{n!} \left( \frac{1}{2\pi} \log |2^j| \right)^n \right| \lesssim \exp\left( \frac{\alpha^2}{2\pi} \log |2^j| \right) \lesssim 2^{j\alpha^2/(2\pi)},$$
which yields
$$\sum_{j \ge -1} 2^{js} c^j \lesssim \sum_{j \ge -1} 2^{j(\alpha^2/(2\pi) + s)}.$$

Therefore, we need
$$s < -\frac{\alpha^2}{4\pi}.$$

Now, we have

$$\begin{aligned}
\mathbb{E}[\|{:}\exp(\alpha(g_\varepsilon * X)){:} - {:}\exp(\alpha X){:}\|^2_{B^s_{2,2,\ell}}] &= \sum_{j \geq -1} 2^{2js} \mathbb{E}\|\Delta_j {:}\exp(\alpha(g_\varepsilon * X)){:} - \Delta_j {:}\exp(\alpha X){:}\|^2_{L^2} \\
&\lesssim \sum_{j \geq -1} 2^{2js} \mathbb{E}[\langle K_j, {:}\exp(\alpha(g_\varepsilon * X)){:} - {:}\exp(\alpha X){:}\rangle^2] \\
&\lesssim \sum_{j \geq -1} \sum_n 2^{2js} \frac{\alpha^{2n}}{(n!)^2} \mathbb{E}[\langle K_j, {:}(g_\varepsilon * X)^n{:} - {:}X^n{:}\rangle^2].
\end{aligned}$$

Notice that each term in the sum converges to 0 in $n$ and $j$. Moreover, we have an uniform bound in $n$ and $j$, since

$$2^{2js} \frac{\alpha^{2n}}{(n!)^2} \mathbb{E}[\langle K_j, {:}(g_\varepsilon * X)^n{:} - {:}X^n{:}\rangle^2] \lesssim 2^{2js} c_n^j,$$

where the term on the right-hand side is summable, and therefore it is in $\ell^1(\mathbb{N}^2)$. By Lebesgue dominated convergence theorem, everything converges to zero.

We follow the proof of Theorem 3.8 in [BD21] to address the case $r > 2$. Let $r > 2$, take $\gamma > 1$, and recall hypercontractivity of Gaussian Wick monomials (see e.g. Chapter III in [Üst95]):

$$\mathbb{E}[({:}X{:})^r] \leq (r-1)^{n/2} (\mathbb{E}[({:}X{:})^2])^{r/2}.$$

We have

$$\begin{aligned}
&\mathbb{E}[|[({:}e^{\alpha(g_\varepsilon X)}{:} - {:}e^{\alpha X}{:}) * K_j](0)|^r]^{1/r} \\
&= \mathbb{E}[|\langle ({:}e^{\alpha(g_\varepsilon X)}{:} - {:}e^{\alpha X}{:}), K_j\rangle|^r]^{1/r} \\
&= \sum_{n=0}^{+\infty} \frac{\alpha^n}{n!} \mathbb{E}[|\langle {:}(g_\varepsilon * X)^n{:} - {:}X^n{:}, K_j\rangle|^r]^{1/r} \\
&\leq \sum_{n=0}^{+\infty} \frac{(\alpha\sqrt{r-1})^n}{n!} \mathbb{E}[|\langle {:}(g_\varepsilon * X)^n{:} - {:}X^n{:}, K_j\rangle|^2]^{1/2} \\
&\lesssim \left(\frac{\kappa^2}{\kappa^2 - 1}\right)^{1/2} \left(\sum_{n=0}^{+\infty} \frac{\kappa^2 \alpha^{2n}(r-1)^n}{(n!)^2} \mathbb{E}[|\langle {:}(g_\varepsilon * X)^n{:} - {:}X^n{:}, K_j\rangle|^2]\right)^{1/2} \\
&\lesssim \left(\frac{\kappa^2}{\kappa^2 - 1}\right)^{1/2} \mathbb{E}[|\langle {:}e^{\alpha\kappa\sqrt{r-1}(g_\varepsilon X)}{:} - {:}e^{\alpha\kappa\sqrt{r-1}X}{:}, K_j\rangle|^2]^{1/2}.
\end{aligned}$$

Taking $\kappa$ such that

$$s < -\frac{\alpha^2 \kappa^2 (r-1)}{4\pi},$$

we conclude the proof in the same way as in the case $r = 2$. □

We need also a result for the periodic setting. In particular, we prove the convergence of the Wick exponential ${:}e^{\alpha Q_{N,M}(g_\varepsilon * X)}{:}$ introduced in equation (4.2), as $N \to +\infty$..

**Proposition B.2.** *Recall that $\gamma = \alpha^2/(4\pi)$. If $2 \leq p < 2/\gamma$, $\delta > 0$, $\ell > 0$, and $\ell' > \ell'_0(p)$, for a positive constant depending on $p$, then we have the convergence, as $N \to +\infty$,*

$$({:}e^{\alpha Q_{N,M}(g_\varepsilon * X)}{:} - {:}e^{\alpha(g_\varepsilon * X)}{:}) \to 0, \qquad \text{in } L^p((C^{-\delta}(\mathbb{T}^2_M), \mu^M_F), L^p(\mathbb{T}^2_M)).$$

**Proof.** Notice that, if $Z$ is a Gaussian random variable, then

$$:e^{\beta Z}: = e^{\beta Z - \frac{\beta^2}{2}\mathbb{E}[Z^2]}.$$

Therefore, applying the previous reasoning for $Z = Q_{N,M}(g_\varepsilon * X)$ and first $\beta = \alpha$, and then $\beta = p\alpha$, then we have

$$(: e^{\alpha Q_{N,M}(g_\varepsilon * X)}:)^p = e^{\alpha p Q_{N,M}(g_\varepsilon * X) - \frac{1}{2}p\alpha^2 \mathbb{E}[(Q_{N,M}(g_\varepsilon * X))^2]} =: e^{\alpha p Q_{N,M}(g_\varepsilon * X)}: e^{\frac{\alpha^2}{2}(p^2 - p)\mathbb{E}[(Q_{N,M}(g_\varepsilon * X))^2]}.$$

We observe that, if $X \in C^{-\delta}(\mathbb{T}_M^2)$, then, for every $x \in \mathbb{T}_M^2$,

$$Q_{N,M}(g_\varepsilon * X(x)) \to g_\varepsilon * X(x), \qquad \text{as } N \to +\infty.$$

Moreover, we have that

$$\mathbb{E}[(Q_{N,M}(g_\varepsilon * X))^2] \to \mathbb{E}[(g_\varepsilon * X)^2],$$

and also

$$: e^{\alpha Q_{N,M} g_\varepsilon * X}: \in L^p_{\mu_M^{\text{free}}}, \qquad \text{uniformly}.$$

These last properties imply that

$$: e^{\alpha Q_{N,M}(g_\varepsilon * X)}: \longrightarrow : e^{\alpha(g_\varepsilon * X)}:, \qquad \text{in } L^p_{\mu_M^{\text{free}}}. \tag{B.2}$$

Taking the norms, we have by translation invariance

$$\mathbb{E}\left[\int_{\mathbb{T}_M^2} |(: e^{\alpha Q_{N,M}(g_\varepsilon * X)}: -: e^{\alpha(g_\varepsilon * X)}:)(z)|^p dz\right] = (2\pi M)^2 \mathbb{E}[|(: e^{\alpha Q_{N,M}(g_\varepsilon * X)}: -: e^{\alpha(g_\varepsilon * X)}:)(0)|^p].$$

and the last term converges to zero by (B.2). □

We now prove convergence in space-time of the Wick exponential.

**Proposition B.3.** *Consider the same parameters $\alpha, p, \ell, \gamma$ as in Proposition B.1. Let $(X_t)_{t \in \mathbb{R}_+}$ be the solution of equation (3.3) with $X_0$ distributed as a Gaussian free field of mass $m$. Then we have*

$$:\exp(\alpha g_\varepsilon * X_t): \to :\exp(\alpha X_t):, \qquad \text{in } L^p(B_X \times \Omega, L^p_{\ell'}(\mathbb{R}_+, B_{p,p,\ell}^{-\gamma(p-1)-\delta}(\mathbb{R}^2))),$$

*where $\ell' > 0$ is such that $\ell' p > 1$.*

**Proof.** The proof follows closely the one of Proposition B.1. See also Theorem 3.2 in [HKK22] and Lemma 2.5 in [HKK21]. □

# Appendix C  Estimates on linearized PDEs

Consider the partial differential equation

$$(\partial_t - \Delta + m^2)\psi(t, z) = -B(t, z, \psi(t, z))(A\psi(t, z) + C(t, z)), \qquad (t, z) \in \mathbb{R}_+ \times \mathbb{R}^2, \tag{C.1}$$

where $B \colon \mathbb{R}_+ \times \mathbb{R}^2 \times H^1(\mathbb{R}^2) \to \mathbb{R}$ is a positive function with compact support with respect to $z \in \mathbb{R}^2$ independently of $t \in \mathbb{R}_+$, $A \colon B^s_{2,2,\ell}(\mathbb{R}^2) \to C^\infty(\mathbb{R}^2)$ is a linear and bounded operator which is self-adjoint with respect to the $L^2(\mathbb{R}^2)$ Hilbert space structure and which commutes with the Laplacian $-\Delta$, while $C \colon \mathbb{R}_+ \times \mathbb{R}^2 \to \mathbb{R}$ is a measurable function such that, for any $t \geq 0$, we have

$$\int_0^t \|B(s, \cdot, \psi(s, \cdot)) C(s, \cdot)\|_{L^2}^2 \, ds < +\infty,$$

where $\psi \in L^2(\mathbb{R}_+, H^1(\mathbb{R}^2)) \cap L^\infty(\mathbb{R}_+, L^2(\mathbb{R}^2))$.

In this section, we prove some a priori estimates for equation (C.1).

**Theorem C.1.** *Consider equation (C.1) with $\psi(0, z) \equiv 0$. We have the bound, for some constant $K > 0$ and for every $0 < \sigma < 1$,*

$$\|A^{1/2} \psi(t, \cdot)\|_{L^2}^2 + \int_0^t (\|A^{1/2} \psi(s, \cdot)\|_{H^1}^2 + (m^2 - K\sigma) \|A^{1/2} \psi(s, \cdot)\|_{L^2}^2) \, ds \lesssim \int_0^t \|B(s, \cdot, \psi(s, \cdot)) C(s, \cdot)\|_{L^2}^2 \, ds,$$

*where the constant implied in the symbol $\lesssim$ does not depend neither on $B$ nor on $C$.*

**Proof.** Multiplying equation (C.1) by $A\psi(t, z)$ and integrating we have,

$$\|A^{1/2} \psi(t, \cdot)\|_{L^2}^2 + \int_0^t (\|A^{1/2} \psi(s, \cdot)\|_{H^1}^2 + m^2 \|A^{1/2} \psi(s, \cdot)\|_{L^2}^2) \, ds$$

$$= -\int_0^t \int B(s, z, \psi(s, z))(A\psi(s, z) + C(s, z)) A\psi(s, z) \, dz \, ds$$

$$= -\int_0^t \int B(s, z, \psi(s, z))(A\psi(s, z))^2 \, dz \, ds$$

$$\quad - \int_0^t \int B(s, z, \psi(s, z)) C(s, z) A\psi(s, z) \, dz \, ds$$

$$\leq -\int_0^t \int B(s, z, \psi(s, z)) C(s, z) A\psi(s, z) \, dz \, ds.$$

Then, exploiting Young's inequality,

$$\|A^{1/2} \psi(t, \cdot)\|_{L^2}^2 + \int_0^t (\|A^{1/2} \psi(s, \cdot)\|_{H^1}^2 + m^2 \|A^{1/2} \psi(s, \cdot)\|_{L^2}^2) \, ds$$

$$\leq -\int_0^t \int B(s, z, \psi(s, z)) C(s, z) A\psi(s, z) \, dz \, ds$$

$$\lesssim C_\sigma \int_0^t \int (B(s, z, \psi(s, z)) C(s, z))^2 \, dz \, ds + C\sigma \int_0^t \|A \psi(s, \cdot)\|_{L^2}^2 \, ds.$$

Since $\|A\psi(s, \cdot)\|_{L^2} = \|A^{1/2} A^{1/2} \psi(s, \cdot)\|_{L^2} \lesssim \|A^{1/2} \psi(s, \cdot)\|_{L^2}$, we can reabsorb the last term on the right-hand side and apply Young's convolution inequality to the remaining term to get

$$\|A^{1/2} \psi(t, \cdot)\|_{L^2}^2 + \int_0^t (\|A^{1/2} \psi(s, \cdot)\|_{H^1}^2 + (m^2 - C\sigma) \|A^{1/2} \psi(s, \cdot)\|_{L^2}^2) \, ds$$

$$\lesssim C_\sigma \int_0^t \int (B(s, z, \psi(s, z)) C(s, z))^2 \, dz \, ds$$

$$\lesssim \int_0^t \|B(s, \cdot, \psi(s, \cdot)) C(s, \cdot)\|_{L^2}^2 \, ds.$$

This concludes the proof. □

Let us modify the previous result in order to deal with weighted norms. In particular, we consider the case where $A\psi = g_\varepsilon * \psi$, where $g_\varepsilon$ is defined as in Section 3, with $A^{1/2} \psi = \tilde{g}_\varepsilon * \psi$. Let us also recall the definition $\rho_\ell^k(z) = (1 + k|z|^2)^{-\ell/2}$, for $z \in \mathbb{R}^2$ and $k > 0$.

**Theorem C.2.** *Consider equation ([C.1](#)) with $\psi(0, z) \equiv 0$. Assume further that*

$$\|B(\cdot, \cdot, \psi)\|_{L^\infty([0,t], L^\infty_{-\ell}(\mathbb{R}^2))} < +\infty.$$

*We have the bound, for some constant $C_1, C_2 > 0$ and for every $0 < \sigma, \tilde{\sigma} < 1$,*

$$\|A^{1/2}\psi(t, \cdot)\|^2_{L^2_{-\ell/2}} + \int_0^t \left((1 - \sigma C_1)\|A^{1/2}\psi(s, \cdot)\|^2_{H^1_{-\ell/2}} + (m^2 - \tilde{\sigma} C_2)\|A^{1/2}\psi(s, \cdot)\|^2_{L^2_{-\ell/2}}\right) ds$$

$$\lesssim \int_0^t \|B(s, \cdot, \psi(s, \cdot)) C(s, \cdot)\|^2_{L^2_{-\ell/2}} ds.$$

**Proof.** Instead of multiplying equation ([C.1](#)) by $A\psi(t, z)$, we multiply it by $A^{1/2}(\rho^k_{-\ell} A^{1/2}\psi(t, z))$. Proceeding as in the previous proof and noticing that

$$|A^{1/2}(\rho^k_{-\ell} A^{1/2}\psi(t, z))| \lesssim \rho^k_{-\ell}(z)|A^{1/2}\psi(t, z)|,$$

we have, integrating by parts,

$$\|\rho^k_{-\ell/2} A^{1/2}\psi(t, \cdot)\|^2_{L^2} + \int_0^t (\|\rho^k_{-\ell/2} A^{1/2}\psi(s, \cdot)\|^2_{H^1} + m^2\|\rho^k_{-\ell/2} A^{1/2}\psi(s, \cdot)\|^2_{L^2}) ds$$

$$+ \int_0^t \int (\nabla A^{1/2}\psi(s, z))(A^{1/2}\psi(s, z))\nabla \rho^k_{-\ell}(z) dz ds$$

$$= -\int_0^t \int B(s, z, \psi(s, z))(A\psi(s, z) + C(s, z))A^{1/2}(\rho^k_{-\ell} A^{1/2}\psi(s, z)) dz ds$$

$$= -\int_0^t \int B(s, z, \psi(s, z))(A\psi(s, z))(A^{1/2}(\rho^k_{-\ell}(z) A^{1/2}\psi(s, z))) dz ds$$

$$- \int_0^t \int B(s, z, \psi(s, z))C(s, z)(A^{1/2}(\rho^k_{-\ell}(z) A^{1/2}\psi(s, z))) dz ds.$$

Let us focus on the term

$$\int_0^t \int (\nabla A^{1/2}\psi(s, z))(A^{1/2}\psi(s, z))\nabla \rho^k_{-\ell}(z) dz ds.$$

Multiplying and dividing by $\rho^k_{-\ell}(z)$ inside the integrals gives

$$\int_0^t \int \rho^k_{-\ell}(z) (\nabla A^{1/2}\psi(s, z)) (A^{1/2}\psi(s, z)) \frac{\nabla \rho^k_{-\ell}(z)}{\rho^k_{-\ell}(z)} dz ds.$$

By Young's inequality we have

$$\left|\int_0^t \int \rho^k_{-\ell}(z) (\nabla A^{1/2}\psi(s, z)) (A^{1/2}\psi(s, z)) \frac{\nabla \rho^k_{-\ell}(z)}{\rho^k_{-\ell}(z)} dz ds\right|$$

$$\lesssim C_\sigma \int_0^t \int \rho^k_{-\ell}(z) (A^{1/2}\psi(s, z))^2 \left(\frac{\nabla \rho^k_{-\ell}(z)}{\rho^k_{-\ell}(z)}\right)^2 dz ds$$

$$+ \sigma C_1 \int_0^t \int \rho^k_{-\ell}(z)(\nabla A^{1/2}\psi(s, z))^2 dz ds,$$

and now we have

$$\frac{\nabla \rho^k_{-\ell}(z)}{\rho^k_{-\ell}(z)} \simeq \frac{\ell k |z_1|}{1 + k|z|^2} \leq \ell \sqrt{k} \sup_{z \in \mathbb{R}^2} \frac{|z_1|}{1 + |z|^2} \leq \sqrt{k}\, C_2.$$

This yields

$$\|\rho^k_{-\ell/2}A^{1/2}\psi(t,\cdot)\|^2_{L^2} + \int_0^t \left((1-\sigma C_1)\|\rho^k_{-\ell/2}A^{1/2}\psi(s,\cdot)\|^2_{H^1} + (m^2 - \sqrt{k}\,C_\sigma)\|\rho^k_{-\ell/2}A^{1/2}\psi(s,\cdot)\|^2_{L^2}\right) ds$$

$$\simeq -\int_0^t\!\!\int B(s,z,\psi(s,z))(A\psi(s,z))(A^{1/2}(\rho^k_{-\ell}(z)A^{1/2}\psi(s,z)))\,dz\,ds$$
$$- \int_0^t\!\!\int B(s,z,\psi(s,z))C(s,z)(A^{1/2}(\rho^k_{-\ell}(z)A^{1/2}\psi(s,z)))\,dz\,ds.$$

Applying Young's inequality we get

$$\|\rho^k_{-\ell/2}A^{1/2}\psi(t,\cdot)\|^2_{L^2} + \int_0^t \left((1-\sigma C_1)\|\rho^k_{-\ell/2}A^{1/2}\psi(s,\cdot)\|^2_{H^1} + (m^2 - \sqrt{k}\,C_\sigma)\|\rho^k_{-\ell/2}A^{1/2}\psi(s,\cdot)\|^2_{L^2}\right) ds$$

$$\lesssim C_{\sigma_2}\int_0^t\!\!\int \rho^k_{-\ell}(z)|B(s,z,\psi(s,z))|^2 |A\psi(s,z)|^2\,dz\,ds + \sigma_2\int_0^t\!\!\int \rho^k_{-\ell}(z)|A^{1/2}\psi(s,z)|^2\,dz\,ds$$
$$+ C_{\sigma_3}\int_0^t\!\!\int (\rho^k_{-\ell/2}(z)B(s,z,\psi(s,z))C(s,z))^2\,dz\,ds + \sigma_3\int_0^t\!\!\int (\rho^k_{-\ell/2}(z)A^{1/2}\psi(s,z))^2\,dz\,ds$$

$$\lesssim C_{\sigma_2}\|B(\cdot,\cdot,\psi)\|_{L^\infty([0,t],L^\infty_{-\ell}(\mathbb{R}^2))}\int_0^t\!\!\int |B(s,z,\psi(s,z))|\,|A\psi(s,z)|^2\,dz\,ds$$
$$+ \sigma_2\int_0^t \|\rho^k_{-\ell/2}A^{1/2}\psi(s,\cdot)\|^2_{L^2}\,ds$$
$$+ C_{\sigma_3}\int_0^t \|\rho^k_{-\ell/2}B(s,\cdot,\psi(s,\cdot))\,C(s,\cdot)\|^2_{L^2}\,ds + \sigma_3\int_0^t \|\rho^k_{-\ell/2}A^{1/2}\psi(s,\cdot)\|^2_{L^2}\,ds,$$

reabsorbing the terms multiplied by $\sigma_2$ and $\sigma_3$ respectively, to the left-hand side and noticing that from the proof of Theorem C.1 we also get

$$\int_0^t\!\!\int B(s,z,\psi(s,z))(A\psi(s,z))^2\,dz\,ds \lesssim \int_0^t \|B(s,\cdot,\psi(s,\cdot))C(s,\cdot)\|^2_{L^2}\,ds,$$

we have, renaming the constants and introducing $\tilde{\sigma}$,

$$\|A^{1/2}\psi(t,\cdot)\|^2_{L^2_{-\ell/2}} + \int_0^t \left((1-\sigma C_1)\|A^{1/2}\psi(s,\cdot)\|^2_{H^1_{-\ell/2}} + (m^2 - \tilde{\sigma}C_2)\|A^{1/2}\psi(s,\cdot)\|^2_{L^2_{-\ell/2}}\right) ds$$
$$\lesssim \int_0^t \|B(s,\cdot,\psi(s,\cdot))\,C(s,\cdot)\|^2_{L^2_{-\ell/2}}\,ds.$$

This concludes the proof. □

We now apply a bootstrap argument to the previous result to get the following statement.

**Corollary C.3.** *Under the same hypotheses of Theorem C.2, we have*

$$\|\psi\|_{B^\beta_{q,q,\ell}(\mathbb{R}_+,B^{2-2\beta-\delta}_{p,p,-\ell}(\mathbb{R}^2))} \lesssim \mathfrak{P}_2(\|B\|_{L^\infty_\ell(\mathbb{R}_+,L^\infty_{-\ell}(\mathbb{R}^2))}, \|C\|_{L^\infty_{\ell'}(\mathbb{R}_+,L^\infty_{-\ell}(\mathbb{R}^2))}),$$

*where $\mathfrak{P}_2$ is a second degree polynomial.*

**Proof.** From Theorem C.2, we know that $A^{1/2}\psi$ lives in $L^2_\ell(\mathbb{R}_+, H^1_{-\ell}(\mathbb{R}^2)) \cap L^\infty_{\ell'}(\mathbb{R}_+, L^2_{-\ell}(\mathbb{R}^2))$. By interpolation and using the fact that $H^1_{-\ell}(\mathbb{R}^2) \subset B^0_{p,p,-\ell}(\mathbb{R}^2)$, for every $p < +\infty$, we have

$$\|A^{1/2}\psi\|^2_{L^q_{\ell/2}(\mathbb{R}_+,B^0_{p,p,-\ell/2}(\mathbb{R}^2))} \lesssim \int_{\mathbb{R}_+} \rho_{\ell'}(s)\|B(s,\cdot,\psi(s,\cdot))\,C(s,\cdot)\|^2_{L^2_{-\ell/2}}\,ds.$$

Applying the heat kernel to the equation (C.1), we have

$$\psi(t,z) = \int_0^t P_{t-s}(B(s,z,\psi(s,z))(A\psi(s,z) + C(s,z)))\,ds,$$

and therefore, if $1/q + 1/p = 1$ and $\beta > 0$, we have by Theorem A.6,

$$\|\psi\|_{B^{\beta}_{q,q,\ell}(\mathbb{R}_+, B^{2-2\beta-\delta}_{p,p,-\ell}(\mathbb{R}^2))} \lesssim \|B(\cdot, \cdot, \psi(\cdot, \cdot))(A\psi(\cdot, \cdot) + C(\cdot, \cdot))\|_{L^q_\ell(\mathbb{R}_+, B^0_{p,p,-\ell}(\mathbb{R}^2))}.$$

On the other hand, for some second degree polynomial $\mathfrak{P}_2$, we have

$$\|B(\cdot, \cdot, \psi(\cdot, \cdot))(A\psi(\cdot, \cdot) + C(\cdot, \cdot))\|_{L^q_\ell(\mathbb{R}_+, B^0_{p,p,-\ell}(\mathbb{R}^2))}$$
$$\leq \mathfrak{P}_2(\|B\|_{L^\infty_\ell(\mathbb{R}_+, L^\infty_{-\ell}(\mathbb{R}^2))}, \|C\|_{L^\infty_\ell(\mathbb{R}_+, L^\infty_{-\ell}(\mathbb{R}^2))}).$$

This concludes the proof. □

# Appendix D  Proof of Lemma 3.4

We give here the proof of Lemma 3.4. We only show the first part of the result, the second one following in a straightforward way. It is sufficient to prove the following property: for any $G \in \mathcal{F}$ there exists a sequence of cylinder functions $(G_k)_{k \in \mathbb{N}}$ (here we consider a sequence $(G_{N,M})_{N,M \in \mathbb{N}}$ depending on two parameters) such that $\mathscr{L} G_k \to \mathscr{L} G$ point-wise as $k \to +\infty$, and we have the uniform bound

$$|\mathscr{L} G_k| \leq F_G(X), \tag{D.1}$$

for some measurable $F_G \in L^1(\mu^{\text{free}})$, e.g. some polynomial of $X$. And secondly we have to show that $\mathscr{L}_\varepsilon G \to \mathscr{L} G$ point-wise as $\varepsilon \to 0$ with a bound $|\mathscr{L}_\varepsilon G| \leq F_G(X)$, for some $F_G \in L^1(\mu^{\text{free}})$. Since we take $\mu \in \mathcal{M}_{B_Y}$, the result follows from Lebesgue's dominated convergence theorem.

Take $G \in \mathcal{F}$ and $N, M, \in \mathbb{N}$, and define

$$G_{N,M}(X, Y) = G(f_M Q_{N,M}(f_M X), f_M Q_{N,M}(f_M Y)),$$

where $f_M = f(\cdot/M)$, $f: \mathbb{R}^2 \to [0, 1]$ being a compactly supported smooth function in $\mathbb{T}^2_1$ which is identically 1 in a neighborhood of the origin, and $Q_{N,M}$ is the operator defined in Section 4. We need to show that

$$\nabla_X G_{N,M}(X, Y) \to \nabla_X G(X, Y), \quad \text{point-wise in } B^{2-\delta}_{1,1,-\ell}(\mathbb{R}^2), \tag{D.2}$$
$$\nabla_Y G_{N,M}(X, Y) \to \nabla_Y G(X, Y), \quad \text{point-wise in } B^{(2-s)\wedge(\gamma(r-1))+\delta}_{l,l,-\ell}(\mathbb{R}^2) \text{ for } l \in (1, \infty), \tag{D.3}$$
$$\text{tr}(\nabla^2_X G_{N,M}(X, Y)) \to \text{tr}(\nabla^2_X G(X, Y)), \quad \text{point-wise.} \tag{D.4}$$

We focus only on the proof of (D.2) and (D.4) since the limit (D.3) follows with similar arguments. We have,

$$\nabla_X G_{N,M}(X, Y)[h]$$
$$= \nabla_X G(f_M Q_{N,M}(f_M X), f_M Q_{N,M}(f_M Y))[f_M Q_{N,M}(f_M h)],$$

On the other hand, since the integrability parameters in the Besov spaces $\tilde{B}_X$, $B_Y$, $B^{2-\delta}_{1,1,-\ell}(\mathbb{R}^2)$, and $B^{(2-s)\wedge(\gamma(r-1))+\delta}_{l,l,-\ell}(\mathbb{R}^2)$ are finite, the linear operator $Z \mapsto f_M Q_{N,M}(f_M Z)$, where $Z \in \mathfrak{B}$, where $\mathfrak{B}$ is anyone of the Besov spaces listed before, strongly converges to the identity $\text{id}_\mathfrak{B}$ on $\mathfrak{B}$. Therefore, by continuity of $\nabla_X G$, we get that the following convergence holds as $N, M \to +\infty$:

$$\nabla_X G(f_M Q_{N,M}(f_M X), f_M Q_{N,M}(f_M Y)) \to \nabla_X G(X, Y), \quad \text{in } B^{2-\delta}_{1,1,-\ell}(\mathbb{R}^2).$$

By strong convergence of the operator $Z \mapsto f_M Q_{N,M}(f_M Z)$ in $B^{2-\delta}_{1,1,-\ell}(\mathbb{R}^2)$ and the fact that the composition of strongly convergent operators is strongly convergent, we get the limit in (D.2).

Let us now prove (D.4). We have, with the notation $p_{N,M} = (f_M Q_{N,M}(f_M X), f_M Q_{N,M}(f_M Y))$, and $p = (X, Y)$,

$$\begin{aligned}
& \operatorname{tr}(f_M Q_{N,M} f_M \nabla^2 G(p_{N,M}) f_M Q_{N,M} f_M) - \operatorname{tr}(\nabla^2 G(p)) \\
= \; & \operatorname{tr}(f_M Q_{N,M} f_M \nabla^2 G(p_{N,M}) f_M Q_{N,M} f_M) - \operatorname{tr}(f_M Q_{N,M} f_M \nabla^2 G(p) f_M Q_{N,M} f_M) \\
& + \operatorname{tr}(f_M Q_{N,M} f_M \nabla^2 G(p) f_M Q_{N,M} f_M) - \operatorname{tr}(\nabla^2 G(p)).
\end{aligned} \tag{D.5}$$

Let us deal first with the second term of equation (D.5), that is

$$\operatorname{tr}(f_M Q_{N,M} f_M \nabla^2 G(p) f_M Q_{N,M} f_M) - \operatorname{tr}(\nabla^2 G(p)). \tag{D.6}$$

We have

$$\begin{aligned}
\operatorname{tr}(f_M Q_{N,M} f_M \nabla^2 G(p) f_M Q_{N,M} f_M) & = \operatorname{tr}(\nabla^2 G(p)(f_M Q_{N,M} f_M)^2) \\
& = \sum_{n \in \mathbb{N}} (\nabla^2 G(p)(f_M Q_{N,M} f_M)^2 h_n, h_n)_{L^2} \\
& = \sum_{n \in \mathbb{N}} ((f_M Q_{N,M} f_M)^2 h_n, \nabla^2 G(p) h_n)_{L^2} \\
& = \sum_{n \in \mathbb{N}} ((f_M Q_{N,M} f_M)^2 h_n^G, \nabla^2 G(p) h_n^G)_{L^2} \\
& = \sum_{n \in \mathbb{N}} \lambda_n ((f_M Q_{N,M} f_M)^2 h_n^G, h_n^G)_{L^2},
\end{aligned}$$

where we use that $f_M$, $Q_{N,M}$, and $\nabla^2 G(p)$ are self-adjoint operators in $L^2(\mathbb{R}^2)$, and $(h_n^G)_{n \in \mathbb{N}}$ is an orthonormal basis of $L^2(\mathbb{R}^2)$ of eigenvectors related to the eigenvalues $(\lambda_n)_{n \in \mathbb{N}}$ of the operator $\nabla^2 G(p)$ (which exists being $\nabla^2 G(p)$ a compact operator). Since $\nabla^2 G(p)$ is also a trace-class operator, we have that $(|\lambda_n|)_{n \in \mathbb{N}} \in \ell^1(\mathbb{N})$. Furthermore, we get

$$|\lambda_n ((f_M Q_{N,M} f_M)^2 h_n^Q, h_n^Q)| \leq |\lambda_n| \|(f_M Q_{N,M} f_M)^2\|_{L(L^2, L^2)} \leq |\lambda_n| \left( \sup_{N,M \in \mathbb{N}} \|(f_M Q_{N,M} f_M)^2\|_{L(L^2, L^2)} \right) \lesssim |\lambda_n|.$$

Since we have the strong convergences $f_M \to \operatorname{id}_{L^2(\mathbb{R}^2)}$ as $M \to +\infty$, and $Q_{N,M} \to \operatorname{id}_{L^2(\mathbb{T}^2_M)}$ as $N \to +\infty$, and so $\sup_M \|f_M\|_{L(L^2, L^2)}, \sup_{N,M} \|Q_{N,M}\|_{L(L^2, L^2)} < +\infty$, then we have

$$(f_M Q_{N,M} f_M)^2 \to \operatorname{id}_{L^2(\mathbb{R}^2)} \text{ strongly.}$$

Thus, we have $(f_M Q_{N,M} f_M)^2 h_n^G \to h_n^G$ in $L^2$ as $N, M \to +\infty$, and therefore

$$\{\lambda_n ((f_M Q_{N,M} f_M)^2 h_n^G, h_n^G)\}_{n \in \mathbb{N}} \to \{\lambda_n\}_{n \in \mathbb{N}} \text{ point-wise.}$$

Then the term (D.6) converges to zero by Lebesgue's dominated convergence theorem.

We now show that the first term in equation (D.5), i.e.

$$\operatorname{tr}(f_M Q_{N,M} f_M \nabla^2 G(p_{N,M}) f_M Q_{N,M} f_M) - \operatorname{tr}(f_M Q_{N,M} f_M \nabla^2 G(p) f_M Q_{N,M} f_M). \tag{D.7}$$

converges to zero. Let $\mathfrak{M}: H_\ell^\kappa(\mathbb{R}^2) \to L^2(\mathbb{R}^2)$ be the natural isomorphism between the two spaces, where $\kappa$ and $\ell$ are the same parameter as in point *iii.* of Definition 3.1. The natural identification of $L^2(\mathbb{R}^2)$ with its dual, allows us to identify the dual map $\mathfrak{M}^*$ of $\mathfrak{M}$ with the natural isomorphism between $L^2(\mathbb{R}^2)$ and $H_\ell^{-\kappa}(\mathbb{R}^2)$. We can then write

$$\nabla^2 G(p) = \mathfrak{M}^{-1} \mathfrak{M} \nabla^2 G(p) \mathfrak{M}^* (\mathfrak{M}^*)^{-1},$$

and therefore, for two points $p, p_{N,M} \in \tilde{B}_X \times B_Y$,

$$\begin{aligned}
&\operatorname{tr}(|\nabla^2 G(p) - \nabla^2 G(p_{N,M})|) \\
&\leq \operatorname{tr}(\mathfrak{M}^{-1}|\mathfrak{M}\nabla^2 G(p)\mathfrak{M}^* - \mathfrak{M}\nabla^2 G(p_{N,M})\mathfrak{M}^*|(\mathfrak{M}^*)^{-1}) \\
&= \sum_{n\in\mathbb{N}} (\mathfrak{M}^{-1}|\mathfrak{M}\nabla^2 G(p)\mathfrak{M}^* - \mathfrak{M}\nabla^2 G(p_{N,M})\mathfrak{M}^*|(\mathfrak{M}^*)^{-1} h_n, h_n)_{L^2} \\
&= \sum_{n\in\mathbb{N}} (|\mathfrak{M}\nabla^2 G(p)\mathfrak{M}^* - \mathfrak{M}\nabla^2 G(p_{N,M})\mathfrak{M}^*|(\mathfrak{M}^*)^{-1} h_n, (\mathfrak{M}^*)^{-1} h_n)_{L^2} \\
&\leq \|\mathfrak{M}\nabla^2 G(p)\mathfrak{M}^* - \mathfrak{M}\nabla^2 G(p_{N,M})\mathfrak{M}^*\|_{L(L^2,L^2)} \sum_{n\in\mathbb{N}} \|(\mathfrak{M}^*)^{-1} h_n\|_{L^2}^2 \\
&\leq \|\mathfrak{M}\nabla^2 G(p)\mathfrak{M}^* - \mathfrak{M}\nabla^2 G(p_{N,M})\mathfrak{M}^*\|_{L(L^2,L^2)} \operatorname{tr}(\iota_{H_\ell^\kappa \hookrightarrow L^2} \mathfrak{M}^{-1} (\mathfrak{M}^*)^{-1} \iota_{L^2 \hookrightarrow H_\ell^{-\kappa}}),
\end{aligned}$$

which is finite since point *iii.* in Definition 3.1 holds and the operator $\iota_{H_\ell^\kappa \hookrightarrow L^2} \mathfrak{M}^{-1} (\mathfrak{M}^*)^{-1} \iota_{L^2 \hookrightarrow H_\ell^{-\kappa}}$ is trace class because $\kappa > 1$ and $\ell > 1$ (see Remark 3.2). By continuity of the map $\nabla^2 G: \tilde{B}_X \times B_Y \to L(H_\ell^{-\kappa}(\mathbb{R}^2), H_\ell^\kappa(\mathbb{R}^2))$ and the fact that, by similar arguments as the one exploited to show (D.2), $p_{N,M} \to p$ as $N, M \to +\infty$, we have that $\operatorname{tr}(|\nabla^2 G(p) - \nabla^2 G(p_{N,M})|) \to 0$ as $N, M \to +\infty$.

The argument is then concluded by the following observation: if we take the absolute value of the difference in (D.7), then

$$\begin{aligned}
&|\operatorname{tr}(f_M Q_{N,M} f_M (\nabla^2 G(p_{N,M}) - \nabla^2 G(p)) f_M Q_{N,M} f_M)| \\
&\leq \operatorname{tr}(f_M Q_{N,M} f_M |\nabla^2 G(p_{N,M}) - \nabla^2 G(p)| f_M Q_{N,M} f_M) \\
&\leq \sup_{N,M\in\mathbb{N}} \|f_M Q_{N,M} f_M\|^2 \operatorname{tr}(|\nabla^2 G(p_{N,M}) - \nabla^2 G(p)|),
\end{aligned}$$

which converges to zero since $\sup_{N,M\in\mathbb{N}} \|f_M Q_{N,M} f_M\|^2$ is finite.

In order to get inequality (D.1), we notice that

$$\|\nabla_X G_{N,M}(X,Y)\| \lesssim \|\nabla_X G(X,Y)\| \leq F_G(X), \tag{D.8}$$

where $F_G$ plays the role of $f_\Phi$ in Definition 3.1, as well as similar inequalities for $\nabla_Y G$ and $\operatorname{tr}(\nabla^2 G)$, which is due to the fact that the norm of operator $f_M Q_{N,M} f_M$ is uniformly bounded in $N, M \in \mathbb{N}$.

The convergence $\mathscr{L}_\varepsilon G \to \mathscr{L} G$ as $\varepsilon \to 0$ is proved in Section 4.4.3.

# Appendix E  Technical results for the resolvent equation

We consider here the system of equations (3.3)–(3.4), and give a proof of Proposition 3.8. For notation simplicity, we will write $Y$ instead of $Y^\varepsilon$ when no confusion occurs. Moreover, if not explicitly specified, all the appearing parameters are assumed to be taken as in Definition 2.6.

## E.1 Flow equations

We start with a result about existence and uniqueness of solutions to the system of equations (3.3)–(3.4), that is the first part of the statement in Proposition 3.8.

**Proposition E.1.** *For any $\varepsilon > 0$, if $(X_0, Y_0) \in \hat{B}_X \times \{B_Y \cup B_{\exp}^{r,\ell}\}$, then there exists a unique solution $(X, Y)$ to equations (3.3)–(3.4) such that*

$$(X_t, Y_t) \in \hat{B}_X \times \{B_Y \cup B_{\exp}^{r,\ell}\}, \qquad t \in \mathbb{R}_+.$$

**Proof.** The unique solution to equation (3.3) is given by the explicit formula

$$X_t = P_t X_0 + \int_0^t P_{t-s} \xi_s \, ds \tag{E.1}$$

(see, e.g., Theorem 5.4 in [DZ14]). Showing existence of a solution to equation (3.4) is equivalent to showing existence for

$$\tilde{Y}_t = Y_t - P_t Y_0.$$

We proceed by a Schaefer's fixed point argument (see Theorem 4 in Section 9.2.2 in [Eva10]) to get the result up to a fixed time $T > 0$. Consider the map $\mathcal{J}$ given by, for $t \in [0, T]$,

$$\mathcal{J}_t(\tilde{Y}, Y_0, X) = -\int_0^t P_{t-s} \alpha f_\varepsilon : e^{\alpha(g_\varepsilon * X_s)} : e^{\alpha(g_\varepsilon * \tilde{Y}_s)} e^{\alpha P_s(g_\varepsilon * Y_0)} \, ds.$$

We need to show that $\mathcal{J}(\cdot, Y_0, X) : A \to \bar{A}$ is a continuous and bounded map, where $A$ and $\bar{A}$ are the convex closed subsets, for $\kappa > 0$ small enough, $\delta > 0$, $\theta \in (\kappa \vee \delta, 1)$, $\ell > 0$ large enough, consisting of non-negative functions, and such that

$$A \subset C^{\theta-\kappa}([0, T], C_{-\ell+\kappa}^{2-2\theta-\delta-\kappa}(\mathbb{R}^2)), \qquad \bar{A} \subset C^{\theta}([0, T], C_{-\ell}^{2-2\theta-\delta}(\mathbb{R}^2)).$$

Exploiting the compact embedding

$$C^{\theta-\kappa}([0, T], C_{-\ell+\kappa}^{2-2\theta-\delta-\kappa}(\mathbb{R}^2)) \hookrightarrow C^{\theta}([0, T], C_{-\ell}^{2-2\theta-\delta}(\mathbb{R}^2)),$$

given by Besov embedding and Corollary 3 in [Sim87], we can then proceed applying Schaefer's fixed point theorem to get existence for every compact subset of $\mathbb{R}_+$ of the form $[0, \tau]$, for some $\tau > 0$. We have

$$\begin{aligned}
&\|\mathcal{J}_s(\tilde{Y}, Y_0, X)\|_{C^\theta([0,T], C_{-\ell}^{2-2\theta-\delta}(\mathbb{R}^2))} \\
&\lesssim \|\alpha f_\varepsilon : e^{\alpha(g_\varepsilon * X_s)} : e^{\alpha(g_\varepsilon * \tilde{Y}_s)} e^{\alpha P_s(g_\varepsilon * Y_0)}\|_{L^\infty([0,T], L^\infty(\mathrm{supp}(f_\varepsilon)))} \\
&\lesssim \|\alpha f_\varepsilon : e^{\alpha(g_\varepsilon * X_s)} : e^{\alpha P_s(g_\varepsilon * Y_0)}\|_{L^\infty([0,T], L^\infty(\mathrm{supp}(f_\varepsilon)))} < +\infty,
\end{aligned} \tag{E.2}$$

by the regularization property of the convolution with $g_\varepsilon$.

We have to prove continuity. For every $\tilde{Y}, \tilde{Y}' \in A$, we have

$$\|\mathscr{J}(\tilde{Y}, Y_0, X) - \mathscr{J}(\tilde{Y}', Y_0, X)\|_{C^\theta([0,T], C^{2-2\theta-\delta}_{-t}(\mathbb{R}^2))}$$
$$\lesssim \|\alpha f_\varepsilon : e^{\alpha(g_\varepsilon * X_s)} : \left(e^{\alpha(g_\varepsilon * \tilde{Y}_s)} - e^{\alpha(g_\varepsilon * \tilde{Y}'_s)}\right) e^{\alpha P_s(g_\varepsilon * Y_0)}\|_{L^\infty([0,T], L^\infty(\mathrm{supp}(f_\varepsilon)))}$$
$$\lesssim \|e^{\alpha(g_\varepsilon * \tilde{Y}_s)} - e^{\alpha(g_\varepsilon * \tilde{Y}'_s)}\|_{L^\infty([0,T], L^\infty(\mathrm{supp}(f_\varepsilon)))}$$
$$\lesssim \|g_\varepsilon * \tilde{Y}_s - g_\varepsilon * \tilde{Y}'_s\|_{L^\infty([0,T], L^\infty(\mathrm{supp}(f_\varepsilon)))}$$
$$\lesssim \|g_\varepsilon\|_{L^1(\mathbb{R}^2)} \|\tilde{Y}_s - \tilde{Y}'_s\|_{L^\infty([0,T], L^\infty(\mathrm{supp}(f_\varepsilon)))}$$
$$\lesssim \|\tilde{Y}_s - \tilde{Y}'_s\|_{C^{\theta-\kappa}([0,T], C^{2-2\theta-\delta-\kappa}_{-t+\kappa}(\mathbb{R}^2))}.$$

We are left to show uniqueness. Take two solutions $Y$ and $Y'$ to equation (3.4). Notice that their difference is given by $Y - Y' = \tilde{Y} - \tilde{Y}'$ and it satisfies

$$(\partial_t - \Delta + m^2)(\tilde{Y}_t - \tilde{Y}'_t) = -\alpha f_\varepsilon : e^{\alpha(g_\varepsilon * X)} : \left(e^{\alpha(g_\varepsilon * Y_t)} - e^{\alpha(g_\varepsilon * Y'_t)}\right).$$

Introduce a positive function $h: \mathbb{R} \to \mathbb{R}$ such that $h(t) \to 0$ when $t \to -\infty$, and with $|h'(t)| \le Ch(t)$, for some constant $C > 0$. Multiplying the previous expression by $h(t)(\tilde{Y}_t - \tilde{Y}'_t)$ and integrating with respect to time and space, we have

$$\int_0^T \int_{\mathbb{R}^2} h(t)(\tilde{Y}_t - \tilde{Y}'_t)(\partial_t - \Delta + m^2)(\tilde{Y}_t - \tilde{Y}'_t) \, \mathrm{d}z \, \mathrm{d}t$$
$$= -\alpha \int_0^T \int_{\mathbb{R}^2} h(t)(\tilde{Y}_t - \tilde{Y}'_t) f_\varepsilon : e^{\alpha(g_\varepsilon * X)} : \left(e^{\alpha(g_\varepsilon * Y_t)} - e^{\alpha(g_\varepsilon * Y'_t)}\right) \mathrm{d}z \, \mathrm{d}t.$$

For the left-hand side, we get by integration by parts

$$\int_0^T \int_{\mathbb{R}^2} h(t)(\tilde{Y}_t - \tilde{Y}'_t)(\partial_t - \Delta + m^2)(\tilde{Y}_t - \tilde{Y}'_t) \, \mathrm{d}z \, \mathrm{d}t$$
$$= h(T)\|\tilde{Y}_T - \tilde{Y}'_T\|_{L^2}^2 - \int_0^T h'(t)\|\tilde{Y}_t - \tilde{Y}'_t\|_{L^2}^2 \, \mathrm{d}t + \int_0^T h(t)(\|\tilde{Y}_t - \tilde{Y}'_t\|_{H^1}^2 + m^2 \|\tilde{Y}_t - \tilde{Y}'_t\|_{L^2}^2) \, \mathrm{d}t$$
$$\ge h(T)\|\tilde{Y}_T - \tilde{Y}'_T\|_{L^2}^2 - \int_0^T |h'(t)| \|\tilde{Y}_t - \tilde{Y}'_t\|_{L^2}^2 \, \mathrm{d}t + \int_0^T h(t)(\|\tilde{Y}_t - \tilde{Y}'_t\|_{H^1}^2 + m^2 \|\tilde{Y}_t - \tilde{Y}'_t\|_{L^2}^2) \, \mathrm{d}t$$
$$\ge h(T)\|\tilde{Y}_T - \tilde{Y}'_T\|_{L^2}^2 + \int_0^T h(t)(\|\tilde{Y}_t - \tilde{Y}'_t\|_{H^1}^2 + (m^2 - C)\|\tilde{Y}_t - \tilde{Y}'_t\|_{L^2}^2) \, \mathrm{d}t.$$

Notice that the last line is positive. Therefore,

$$h(T)\|\tilde{Y}_T - \tilde{Y}'_T\|_{L^2}^2 + \int_0^T h(t)(\|\tilde{Y}_t - \tilde{Y}'_t\|_{H^1}^2 + (m^2 - C)\|\tilde{Y}_t - \tilde{Y}'_t\|_{L^2}^2) \, \mathrm{d}t$$
$$\lesssim -\alpha \int_0^T \int_{\mathbb{R}^2} h(t)(\tilde{Y}_t - \tilde{Y}'_t) : e^{\alpha(g_\varepsilon * X)} : \left(e^{\alpha(g_\varepsilon * Y_t)} - e^{\alpha(g_\varepsilon * Y'_t)}\right) \mathrm{d}z \, \mathrm{d}t \le 0,$$

which yields uniqueness. □

**Remark E.2.** The growth with respect to $T$ of the norm of the solution $Y$ to equation (3.4) is polynomial. Indeed, recalling that in the proof of Proposition E.1 we have

$$\tilde{Y}^\varepsilon_t = -\int_0^t P_{t-s} \alpha f_\varepsilon : e^{\alpha(g_\varepsilon * X_s)} : e^{\alpha(g_\varepsilon * \tilde{Y}_s)} e^{\alpha P_s(g_\varepsilon * Y_0)} \mathrm{d}s,$$

we get that the time growth of $\tilde{Y}_t$ is determined by the growth with respect to $s$ of the term

$$\alpha f_\varepsilon : e^{\alpha(g_\varepsilon * X_s)} : e^{\alpha(g_\varepsilon * \tilde{Y}_s)} e^{\alpha P_s(g_\varepsilon * Y_0)}.$$

Now, $\exp(\alpha e^{-(-\Delta + m^2)s}(g_\varepsilon * Y_0))$ does not increase in time, while $\exp(\alpha(g_\varepsilon * \tilde{Y}_s))$ is bounded since $\tilde{Y}_s \leq 0$, and hence the growth with respect to $s$ is determined only by the remaining term, that is the exponential $f_\varepsilon : \exp(\alpha(g_\varepsilon * X_s)) :$.

In the next result, we prove continuity of the solutions to equations (3.3)–(3.4) with respect to the initial data.

**Lemma E.3.** *For every $\varepsilon > 0$, the solution $(X, Y)$ to equations (3.3)–(3.4) are continuous with respect to $X_0$ in $\hat{B}_X$ and with respect to $Y_0$ in $B_Y \cup B_{\exp}^{r,\ell}$, respectively.*

**Proof.** As far as $X$ is concerned, continuity follows from the linearity of its representation (E.1). Let us focus on $Y$. Consider a sequence $(Y_0^n)_{n \in \mathbb{N}}$ converging to some limit $Y_0$ in $B_Y \cup B_{\exp}^{r,\ell}$. Then $(Y_0^n)_{n \in \mathbb{N}}$ is bounded in $B_Y \cup B_{\exp}^{r,\ell}$, and, by the regularization properties of $g_\varepsilon$, we have that $\exp(\alpha P_t Y_0^n)$ is uniformly bounded with respect to $z \in \mathbb{R}^2$, $t$ and $n$ on the support of $g_\varepsilon$. From inequality (E.2), it follows that the solution $\tilde{Y}_t(Y_0^{n_k})$ is bounded in the space $\bar{A} \subset C^\theta([0, \tau], C_{-\ell}^{2-2\theta-\delta}(\mathbb{R}^2))$ and pre-compact in $A \subset C^{\theta-\kappa}([0, \tau], C_{-\ell+\kappa}^{2-2\theta-\delta-\kappa}(\mathbb{R}^2))$. Thus, there exists a subsequence $n_k$ such that $\tilde{Y}_t(Y_0^{n_k})$ converges to some $\tilde{Z}$ in $A$.

On the other hand, since $\tilde{Y}_t(Y_0^{n_k})$ solves the fixed-point equation

$$\tilde{Y}_t(Y_0^{n_k}) = \mathscr{F}_t(\tilde{Y}_t(Y_0^{n_k}), Y_0^{n_k}, X),$$

and since $\mathscr{F}$ is continuous with respect to all its variables (see the proof of Proposition E.1), we have

$$\tilde{Z} = \mathscr{F}_t(\tilde{Z}, Y_0, X).$$

By uniqueness of solutions to the fixed point equation, we have that $\tilde{Y}_t(Y_0) = \tilde{Z}$. Finally, since the solution to equation (3.4) with initial data $Y_0^{n_k}$ is given by

$$Y_t(Y_0^{n_k}) = \tilde{Y}_t(Y_0^{n_k}) + P_t Y_0^{n_k},$$

by the continuity of the heat kernel and of $\tilde{Y}_t$, we get the thesis. □

## E.2 Derivatives of the flow

In this section, we denote by $X$ and $Y$ the solutions to the system of equations (3.3)–(3.4) and we study their derivatives with respect to the initial data $X_0$ and $Y_0$.

### E.2.1 Existence and equations

We now show point *i.* in Proposition 3.8.

We have that $X$ is differentiable and that its derivatives $\nabla_{X_0} X_t$, $\nabla^2_{X_0} X_t$, and $\nabla_{Y_0} X_t$ solve equations (3.5), (3.6), and (3.7), respectively. This is because we have the following explicit representation (cf. equation (E.1)) of the solution

$$X_t = P_t X_0 + \int_0^t P_{t-s} \xi_s \, ds,$$

which is linear with respect to the initial data. We immediately get that $\nabla^2_{X_0} X \equiv 0$ and $\nabla_{Y_0} X \equiv 0$.

We now focus on the derivatives $\nabla_{X_0} Y$, $\nabla^2_{X_0} Y$, and $\nabla_{Y_0} Y$ of $Y$, and show that they exist and satisfy equations (3.8), (3.9), and (3.10), respectively. Furthermore, they are all continuous functions with respect to $X_0$ and $Y_0$.

**Proposition E.4.** *For every $\varepsilon > 0$, we have that that the derivatives $\nabla_{Y_0} Y$, $\nabla_{X_0} Y$, and $\nabla^2_{X_0} Y$ of the solution $Y$ to equation (3.4) exist and satisfy equations (3.8), (3.9), and (3.10), respectively.*

**Proof.** We only give the proof for $\nabla_{Y_0} Y$, the other cases follows in a similar way. Consider the approximating equation (3.4), that is

$$(\partial_t - \Delta + m^2) Y = -\mathcal{G}_\varepsilon(X, Y) = -\alpha f_\varepsilon : e^{\alpha(g_\varepsilon * X)}: e^{\alpha(g_\varepsilon * Y)}, \qquad Y(0) = Y_0, \tag{E.3}$$

with $Y_0 \leq 0$. We denote the solution to equation (E.3) as $Y_t(Z)$, in order to stress the initial condition $Z \in B_Y \cup B_{\exp}^{r,\ell}$. We have the integral representation of the solution $Y$ to (E.3):

$$Y_t(Z) = P_t Z - \int_0^t P_{t-s} \mathcal{G}_\varepsilon(X_s, Y_s(Z)) \, ds.$$

In order to compute the derivative with respect to the initial data, we need to perturb it. Therefore, taking $\lambda > 0$ and $h \in B_Y \cup B_{\exp}^{r,\ell}$, and considering the difference, we have

$$Y_t(Y_0 + \lambda h) - Y_t(Y_0) = P_t h \lambda - \int_0^t P_{t-s}(\mathcal{G}_\varepsilon(X_s, Y_s(Y_0 + \lambda h)) - \mathcal{G}_\varepsilon(X_s, Y_s(Y_0))) \, ds.$$

But

$$\int_0^t P_{t-s}(\mathcal{G}_\varepsilon(X_s, Y_s(Y_0 + \lambda h)) - \mathcal{G}_\varepsilon(X_s, Y_s(Y_0))) \, ds$$
$$= \int_0^t \int_0^1 P_{t-s} D_Y \mathcal{G}_\varepsilon(X_s, \varsigma Y_s(Y_0 + \lambda h) + (1-\varsigma) Y_s(Y_0))(g_\varepsilon * (Y_s(Y_0 + \lambda h) - Y_s(Y_0))) \, d\varsigma \, ds,$$

since $\mathcal{G}_\varepsilon$ is differentiable in the direction $g_\varepsilon * (Y_s(Y_0 + \lambda h) - Y_s(Y_0))$.

Define

$$H(s, \lambda, h) = \int_0^1 D_Y \mathcal{G}_\varepsilon(X_s, \varsigma Y_s(Y_0 + \lambda h) + (1-\varsigma) Y_s(Y_0)) \, d\varsigma,$$

so that we can write

$$Y_t(Y_0 + \lambda h) - Y_t(Y_0) = P_t h \lambda - \int_0^t P_{t-s} H(s, \lambda, h) \cdot (g_\varepsilon * (Y_s(Y_0 + \lambda h) - Y_s(Y_0))) \, ds.$$

Let us write down an equation for

$$D_{\lambda, h} Y_t(Y_0) = \frac{Y_t(Y_0 + \lambda h) - Y_t(Y_0)}{\lambda}.$$

Notice that
$$D_{\lambda,h}Y_t(Y_0) = P_t h - \int_0^t P_{t-s} H(s, \lambda, h) \cdot (g_\varepsilon * D_{\lambda,h} Y_t(Y_0))\, ds,$$

and therefore
$$\partial_t (D_{\lambda,h} Y_t(Y_0)) = -(-\Delta + m^2) D_{\lambda,h} Y_t(Y_0) - H(t, \lambda, h)(g_\varepsilon * D_{\lambda,h} Y_t(Y_0)), \qquad D_{\lambda,h} Y_0(Y_0) = h. \tag{E.4}$$

Now consider, as in the proof of Proposition E.1,
$$\tilde{Y}_t = Y_t - P_t Y_0.$$

Then, we have
$$D_{\lambda,h} \tilde{Y}_t(Y_0) = \frac{\tilde{Y}_t(Y_0 + \lambda h) - \tilde{Y}_t(Y_0)}{\lambda} = D_{\lambda,h} Y_t(Y_0) - P_t h,$$

and such a difference satisfies the following equation
$$(\partial_t - \Delta + m^2) D_{\lambda,h} \tilde{Y}_t(Y_0) = -H(t, \lambda, h)(g_\varepsilon * D_{\lambda,h} \tilde{Y}_t(Y_0) + P_t(g_\varepsilon * h)), \quad D_{\lambda,h} Y_0(Y_0) = 0. \tag{E.5}$$

By Theorem C.1, we have then the following bound, for some constant $K > 0$ and every $0 < \sigma < 1$:
$$\|\tilde{g}_\varepsilon * D_{\lambda,h} \tilde{Y}_t(Y_0)\|_{L^2}^2 + \int_0^t (\|\tilde{g}_\varepsilon * D_{\lambda,h} \tilde{Y}_s(Y_0)\|_{H^1}^2 + (m^2 - K\sigma)\|\tilde{g}_\varepsilon * D_{\lambda,h} \tilde{Y}_s(Y_0)\|_{L^2}^2)\, ds$$
$$\lesssim \int_0^t \|H(t, \lambda, h) P_t(g_\varepsilon * h)\|_{L^2}^2\, ds. \tag{E.6}$$

Consider
$$\hat{Y}_t^\lambda(Y_0) = Y_t(Y_0 + \lambda h) - P_t h \lambda = P_t Y_0 - \int_0^t P_{t-s} \mathcal{G}_\varepsilon(X_s, Y_s(Y_0 + \lambda h))\, ds.$$

Notice that $\hat{Y}_t^\lambda(Y_0)$ is negative and moreover it solves
$$\hat{Y}_t^\lambda(Y_0) = P_t Y_0 - \alpha \int_0^t P_{t-s} f_\varepsilon \colon\! e^{\alpha(g_\varepsilon * X_s)} \colon\! e^{\alpha(g_\varepsilon * \hat{Y}_s^\lambda(Y_0))} e^{\lambda \alpha(g_\varepsilon * P_s h)}\, ds.$$

Recall that $g_\varepsilon$ is such that $g_\varepsilon * h \in L^\infty_{\mathrm{loc}}$ and $f_\varepsilon$ is compactly supported. Therefore,
$$\exp(\alpha(1-\varsigma)(g_\varepsilon * Y_s)(Y_0)), \exp(\alpha g_\varepsilon * \hat{Y}_s^\lambda(Y_0)) \in L^\infty, \quad \text{since the exponents are negative,}$$
$$\mathbb{I}_{\mathrm{supp}(f_\varepsilon)} \exp(\lambda \alpha(g_\varepsilon * P_s h)) \in L^\infty, \qquad\qquad \text{since the exponent is negative on } \mathrm{supp}(f_\varepsilon).$$

We then have the uniform estimate on $H$ given by, for any $s > 0$ small enough and $1 < p < +\infty$,
$$\|H(s, \lambda, h)\|_{B_{p,p,-\ell}^{-s}} \leq \|f_\varepsilon \colon\! e^{\alpha(g_\varepsilon * X_s)} \colon\!\|_{B_{p,p,-\ell}^{-s}} \|e^{\alpha\varsigma(g_\varepsilon * \hat{Y}_s^\lambda(Y_0))} e^{\lambda\varsigma\alpha(g_\varepsilon * P_s h)} e^{\alpha(1-\varsigma)(g_\varepsilon * Y_s)(Y_0)}\|_{L^\infty(\mathrm{supp}(f_\varepsilon))}$$
$$\leq \|f_\varepsilon \colon\! e^{\alpha(g_\varepsilon * X_s)} \colon\!\|_{B_{p,p,-\ell}^{-s}} \|e^{\lambda\varsigma\alpha(g_\varepsilon * P_s h)}\|_{L^\infty(\mathrm{supp}(f_\varepsilon))}$$
$$\leq \|f_\varepsilon \colon\! e^{\alpha(g_\varepsilon * X_s)} \colon\!\|_{B_{p,p,-\ell}^{-s}} \|e^{\alpha|g_\varepsilon * P_s h|}\|_{L^\infty(\mathrm{supp}(f_\varepsilon))}, \tag{E.7}$$

which is uniform both in $\lambda$ and $\xi$.

This gives existence of a limit for $D_{\lambda,h} Y_t(Y_0)$ as $\lambda \to 0$.

We are left to show that the limit satisfies equation (3.8). First, we have to prove that the following limit holds

$$\lim_{\lambda \to 0} H(t, \lambda, h) = \alpha^2 f_\varepsilon : e^{\alpha(g_\varepsilon * X_t)} : e^{\alpha(g_\varepsilon * Y_t(Y_0))},$$

in some suitable space.

By Lemma E.3, we have that $Y_t(Y_0 + \lambda h) \to Y_t(Y_0)$ as $\lambda \to 0$ in $B_Y \cup B_{\exp}^{r,\ell}$, and hence, thanks to the regularization properties of $g_\varepsilon$, $g_\varepsilon * Y_t(Y_0 + \lambda h)$ converges to $g_\varepsilon * Y_t(Y_0)$ uniformly on the support of $f_\varepsilon$. Thus, we have the weak convergence

$$H(t, \lambda, h) \to \alpha^2 f_\varepsilon : e^{\alpha(g_\varepsilon * X_t)} : e^{\alpha(g_\varepsilon * Y_t(Y_0))}, \qquad \text{as } \lambda \to 0, \tag{E.8}$$

and, by the uniform estimates (E.7) and the compact embedding $B_{p,p,-\ell}^{-s} \hookrightarrow B_{p,p,0}^{-s'}$, with $s' > s$ (see Proposition A.2), we get the strong convergence in the space $B_{p,p,0}^{-s'}$.

We then have a weak limit for $D_{\lambda,h} Y_t(Y_0)$ as $\lambda \to 0$. Thus, proceeding as in the proof of Lemma E.3, thanks to the a priori estimate (E.6), to the uniform bound (E.7), and to the convergence (E.8), we get that the limit of $D_{\lambda,h} Y_t(Y_0)$ as $\lambda \to 0$ is a solution to equation (3.8). □

### E.2.2 Properties of the flow derivatives

We prove here some bounds on $\nabla_{X_0} Y_t(Y_0)$, $\nabla_{Y_0} Y_t(Y_0)$, and on the trace of $\nabla_{X_0}^2 Y_t(Y_0)$. Let us recall that the sets $\tilde{B}_Y$ and $\tilde{B}_X$ are defined as in (3.11).

**Proposition E.5.** *For every* $\delta \in (0, 1)$, $\theta \in (0, 1 - \delta)$, $\ell, \ell' \geq 1$ *and* $h \in \tilde{B}_Y$, *we have the estimate*

$$\|\nabla_{Y_0} Y_t(Y_0)(h)\|_{C_{\ell'}^\theta(\mathbb{R}_+, C_{-\ell}^{2-2\theta-2\delta}(\mathbb{R}^2)) \oplus L^\infty(\mathbb{R}_+, \tilde{B}_Y)}$$
$$\lesssim_{g_\varepsilon} \mathfrak{P}_2(\|f_\varepsilon : e^{\alpha(g_\varepsilon * X_t)} : e^{\alpha P_t(g_\varepsilon * h)}\|_{L_{\ell'}^\infty(\mathbb{R}_+, L_{-\ell}^\infty(\mathbb{R}^2))}, \|h\|_{\tilde{B}_Y}),$$

*where* $\mathfrak{P}_2$ *is a second degree polynomial.*

**Proof.** We use the same notation as in the proof of Proposition E.4. Recall that

$$D_{\lambda,h} Y_t(Y_0) = D_{\lambda,h} \tilde{Y}_t(Y_0) + P_t h,$$

and then, since $P_t h$ is uniformly bounded in $\tilde{B}_Y$, in order to prove the result, it suffices to give an estimate on $D_{\lambda,h} \tilde{Y}_t(Y_0)$ in the space $C_{\ell'}^\theta(\mathbb{R}_+, C_{-\ell}^{2-2\theta-2\delta}(\mathbb{R}^2))$.

Let us consider a time-weight $\rho_{\ell'}(t)$ and a space-weight $\rho_{-\ell}(x)$ defines as

$$\rho_l^k(z) = (1 + k|z|^2)^{-l/2}, \qquad z \in \mathbb{R}_+, \mathbb{R}^2.$$

Then, by Theorem C.2 we have, for some $0 < \sigma < 1$,

$$\rho_{2\ell'}(t) \|\rho_{-\ell}(\tilde{g}_\varepsilon * D_{\lambda,h} \tilde{Y}_t(Y_0))\|_{L^2}^2$$
$$+ \int_0^t \rho_{2\ell'}(s)(\|\rho_{-\ell}(\tilde{g}_\varepsilon * D_{\lambda,h} \tilde{Y}_s(Y_0))\|_{H^1}^2 + (m^2 - \sigma C)\|\rho_{-\ell}(\tilde{g}_\varepsilon * D_{\lambda,h} \tilde{Y}_s(Y_0))\|_{L^2}^2) \, ds$$
$$\lesssim C_\delta \int_0^t \rho_{2\ell'}(s) \int \rho_{-2\ell}(x)(\tilde{g}_\varepsilon * (H(s, \lambda, h) e^{-(-\Delta + m^2)s}(g_\varepsilon * h)))^2 \, dx \, ds.$$

Moreover, applying Corollary C.3 yields

$$\|H(s,\lambda,h)(\tilde{g}_\varepsilon * D_{\lambda,h}\tilde{Y}_s(Y_0) + P_s(\tilde{g}_\varepsilon * h))\|_{L^q_\ell(\mathbb{R}_+, B^0_{p,p,-\ell}(\mathbb{R}^2))}$$
$$\leq \mathfrak{P}_2(\|H\|_{L^\infty_\ell(\mathbb{R}_+, L^\infty_{-\ell}(\mathbb{R}^2))}, \|\tilde{g}_\varepsilon * h\|_{L^\infty_{-\ell}(\mathbb{R}^2)}).$$

Together with the previous estimate, this yields uniform bounds on the norm $\|D_{\lambda,h}\tilde{Y}_t(Y_0)\|_{B^\beta_{q,q,\ell}(\mathbb{R}_+, B^{2-2\beta-\delta}_{p,p,-\ell}(\mathbb{R}^2))}$, and choosing $\beta$ and $p$ accordingly we then deduce that we have $D_{\lambda,h}\tilde{Y}_t(Y_0) \in C^\delta_\ell(\mathbb{R}_+, C^{2-\delta}_{-\ell}(\mathbb{R}^2))$, uniformly in $\lambda$. Now, letting $\lambda \to 0$ yields the result. □

**Proposition E.6.** *For every $\delta \in (0,1)$, $\theta \in (0, 1-\delta)$, $\ell, \ell' \geq 1$ and $h \in \tilde{B}_X$, we have the estimate*

$$\|\nabla_{X_0} Y_t(Y_0)(h)\|_{C^\theta_\ell(\mathbb{R}_+, C^{2-2\theta-2\delta}_{-\ell}(\mathbb{R}^2))} \lesssim_{g_\varepsilon} P_2(\|f_\varepsilon : e^{\alpha(g_\varepsilon * X_t)} : e^{\alpha P_t(g_\varepsilon * h)}\|_{L^\infty_\ell(\mathbb{R}_+, L^\infty_{-\ell}(\mathbb{R}^2))}),$$

*where $P_2$ is a second degree polynomial.*

**Proof.** The proof is similar to the one of Proposition E.5, the only difference is that here the initial data of $\nabla_{X_0} Y_t$ is zero and therefore we do not need to subtract it before doing the estimates in Theorem C.2 and in Corollary C.3. □

We now deal with the trace term appearing in the definition of the operator $\mathscr{L}$.

**Proposition E.7.** *For every $\ell, \kappa \geq 0$, there exist $\beta, \delta > 0$ such that*

$$\|\nabla^2_{X_0} Y_t(Y_0)\|_{L(H^{-\kappa}_\ell, H^\kappa_{-\ell})} \lesssim \left(\int_{\mathbb{R}^2} \alpha^2 f_\varepsilon(z') : e^{\alpha(g_\varepsilon * X_t)(z')} : e^{\alpha(g_\varepsilon * Y_t)(z')} e^{(\delta - m^2)t}(1+|z'|^\beta)\,\mathrm{d}z'\right)^2.$$

*It follows that, whenever $\ell > 1$ and $\kappa > 1$,*

$$\mathrm{tr}(|\nabla^2_{X_0} Y_t(Y_0)\rho_{-\ell}|) \lesssim \left(\int_{\mathbb{R}^2} \alpha^2 f_\varepsilon(z') : e^{\alpha(g_\varepsilon * X_t)(z')} : e^{\alpha(g_\varepsilon * Y_t)(z')} e^{(\delta - m^2)t}(1+|z'|^\beta)\,\mathrm{d}z'\right)^2.$$

**Proof.** We suppose that all the computations involving equations (3.9) and (3.10) make sense since a rigorous proof of this fact can be given in a similar way as in the proof of Proposition E.5.

Let $\delta_.$ be the Dirac delta distribution, and consider the map, for $X_0 \in B_X$ and $t \in \mathbb{R}_+$,

$$T_{X_0, t} : H^{-\kappa}_\ell(\mathbb{R}^2) \to H^\kappa_{-\ell}(\mathbb{R}^2), \quad h \mapsto \nabla^2_{X_0} Y_t(Y_0)(h, \delta_.).$$

Hereafter, we drop the dependence on $Y_0$ in the derivatives of $Y$. Let us then evaluate equations (3.9) and (3.10) at $\delta_z(\cdot)$, and $(h, \delta_z(\cdot))$, respectively, to get

$$(\partial_t - (\Delta - m^2))\nabla_{X_0} Y_t(\delta_z)(z')$$
$$= -\alpha^2 f_\varepsilon(z') : e^{\alpha g_\varepsilon * X_t(z')} : e^{\alpha g_\varepsilon * Y_t(z')}(P_t(g_\varepsilon * \delta_z)(z') + (g_\varepsilon * \nabla_{X_0} Y_t(\delta_z))(z'))$$
$$= -\alpha^2 f_\varepsilon(z') : e^{\alpha g_\varepsilon * X_t(z')} : e^{\alpha g_\varepsilon * Y_t(z')}(P_t g_\varepsilon(z - z') + (g_\varepsilon * \nabla_{X_0} Y_t(\delta_z))(z')), \qquad (\text{E.9})$$

and

$$\begin{aligned}&(\partial_t - (\Delta - m^2))\nabla^2_{X_0} Y_t(h, \delta_z)(z') \\ &= -\alpha^2 f_\varepsilon : e^{\alpha(g_\varepsilon * X_t)(z')} : e^{\alpha(g_\varepsilon * Y_t)(z')} (g_\varepsilon * \nabla^2_{X_0} Y_t(h, \delta_z))(z') \\ &\quad -\alpha^3 f_\varepsilon : e^{\alpha(g_\varepsilon * X_t)(z')} : e^{\alpha(g_\varepsilon * Y_t)(z')} (g_\varepsilon * P_t h)(z')((g_\varepsilon * P_t \delta_z)(z') + (g_\varepsilon * \nabla_{X_0} Y_t(\delta_z))(z')) \\ &\quad -\alpha^3 f_\varepsilon : e^{\alpha(g_\varepsilon * X_t)(z')} : e^{\alpha(g_\varepsilon * Y_t)(z')} (g_\varepsilon * \nabla_{X_0} Y_t(h))(z') ((g_\varepsilon * P_t \delta_z)(z') + (g_\varepsilon * \nabla_{X_0} Y_t(\delta_z))(z')). \end{aligned} \quad (E.10)$$

Using similar methods as in the proofs of Proposition E.4 and Proposition E.5, it is possible to prove that $\nabla_{X_0} Y_t(\delta_z)$ and $\nabla^2_{X_0} Y_t(h, \delta_z)$ are differentiable infinitely many times with respect to $z$ and their derivatives with respect to $z$ for any multi-index $\beta$ satisfy

$$\partial^\beta_z (g_\varepsilon * \nabla_{X_0} Y_t(\delta_z))(z') = \int g_\varepsilon(z' - y) \partial^\beta_z \nabla_{X_0} Y_t(\delta_z)(y) \, dy.$$

Furthermore, by equations (E.9) and (E.10), we get

$$\begin{aligned}(\partial_t - (\Delta - m^2))\partial^\beta_z \nabla_{X_0} Y_t(\delta_z)(z') &= -\alpha^2 f_\varepsilon(z') : e^{\alpha(g_\varepsilon * X_t)(z')} : e^{\alpha(g_\varepsilon * Y_t)(z')} P_t \partial^\beta_z g_\varepsilon(z - z') \\ &\quad -\alpha^2 f_\varepsilon(z') : e^{\alpha(g_\varepsilon * X_t)(z')} : e^{\alpha(g_\varepsilon * Y_t)(z')} (g_\varepsilon * \partial^\beta_z \nabla_{X_0} Y_t(\delta_z))(z'), \end{aligned} \quad (E.11)$$

and

$$\begin{aligned}&(\partial_t - (\Delta - m^2))\partial^\beta_z \nabla^2_{X_0} Y_t(h, \delta_z)(z') \\ &= -\alpha^2 f_\varepsilon(z') : e^{\alpha(g_\varepsilon * X_t)(z')} : e^{\alpha(g_\varepsilon * Y_t)(z')} (\partial^\beta_. g_\varepsilon * \nabla^2_{X_0} Y_t(h, \delta_z))(z') \\ &\quad -\alpha^3 f_\varepsilon(z') : e^{\alpha(g_\varepsilon * X_t)(z')} : e^{\alpha(g_\varepsilon * Y_t)(z')} (g_\varepsilon * P_t h)(z')((\partial^\beta_. g_\varepsilon * P_t \delta_z)(z') + (g_\varepsilon * \partial^\beta_z \nabla_{X_0} Y_t(\delta_z))(z')) \\ &\quad -\alpha^3 f_\varepsilon(z') : e^{\alpha(g_\varepsilon * X_t)(z')} : e^{\alpha(g_\varepsilon * Y_t)(z')} (g_\varepsilon * \nabla_{X_0} Y_t(h))(z')((\partial^\beta_. g_\varepsilon * P_t \delta_z)(z') + (g_\varepsilon * \partial^\beta_z \nabla_{X_0} Y_t(\delta_z))(z')). \end{aligned} \quad (E.12)$$

Exploiting Theorem C.2 applied to equation (E.11), we get

$$\begin{aligned}&\|\tilde{g}_\varepsilon * \partial^\beta_z \nabla_{X_0} Y_t(\delta_z)\|^2_{L^2} + \int_0^t \|\tilde{g}_\varepsilon * \partial^\beta_z \nabla_{X_0} Y_\tau(\delta_z)\|^2_{H^1} d\tau \\ &\lesssim \int \alpha^2 f_\varepsilon(z') : e^{\alpha(g_\varepsilon * X_t)(z')} : e^{\alpha(g_\varepsilon * Y_t)(z')} |P_t \partial^\beta_z g_\varepsilon(z - z')|^2 dz' =: F_\beta(z). \end{aligned}$$

By a bootstrap argument, i.e. applying Corollary C.3, we can conclude that

$$\|\partial^\beta_z \nabla_{X_0} Y_t(\delta_z)\|^2_{B^s_{2,2}} \lesssim \rho_{-\ell'}(t) F_\beta(z) q_s(F_\beta(z)),$$

where $q_s$ is an $s$-dependent polynomial, because of the relation of the heat kernel with the support of $f_\varepsilon$.

By a similar method as above, we get the following estimate concerning the second derivative

$$\|\partial^\beta_z \nabla^2_{X_0} Y_t(h, \delta_z)\|_{B^s_{2,2}} \lesssim \|h\|_{H^{-\kappa}_{\ell}} F_\beta(z) \tilde{q}_s(F_\beta(z)),$$

where $\tilde{q}_s$ is another $s$-dependent polynomial. Taking $s$ large enough, we have

$$\|\partial^\beta_z \nabla^2_{X_0} Y_t(h, \delta_z)\|_{L^\infty} \lesssim \|h\|_{H^{-\kappa}_\ell} F_\beta(z) \tilde{q}_s(F_\beta(z)).$$

This proves that the map is linear with respect to $h$ as a map from $L^2$, moreover, if $s \in \mathbb{N}$ and fixing $z'$ and $h$, the norm with respect to $z$ is given by

$$\begin{aligned}
\|\nabla^2_{X_0} Y_t(h, \delta_\cdot)(z')\|^2_{B^s_{\ell,2,-\ell}(\mathrm{d}z)} &= \sum_{|\beta| \leq s} \int \rho^2_{-\ell}(z) \, (\partial^\beta_z \nabla^2_{X_0} Y_t(h, \delta_z)(z'))^2 \, \mathrm{d}z \\
&\lesssim \sum_{|\beta| \leq s} \int \rho^2_{-\ell}(z) \|\partial^\beta_z \nabla^2_{X_0} Y_t(h, \delta_z)(z')\|^2_{L^\infty(\mathrm{d}z')} \, \mathrm{d}z \\
&\lesssim \|h\|_{H^{-\kappa}_\ell} \sum_{|\beta| \leq s} \int \rho^2_{-\ell}(z) \, |F_\beta(z) \tilde{q}_s(F_\beta(z))|^2 \, \mathrm{d}z.
\end{aligned} \qquad (E.13)$$

We are left to show that the last integral is finite. Recall

$$F_\beta(z) = \int \alpha^2 f_\varepsilon(z') :e^{\alpha(g_\varepsilon * X_t)(z')}: e^{\alpha(g_\varepsilon * Y_t)(z')} |P_t \partial^\beta_z g_\varepsilon(z - z')|^2 \, \mathrm{d}z'.$$

Since we have to show that the integral is finite for some polynomial, we multiply (and divide) the heat kernel by some weight $\rho_{-2\ell}$, to get

$$F_\beta(z) \rho_{-2\ell}(z) = \int \alpha^2 f_\varepsilon(z') :e^{\alpha(g_\varepsilon * X_t)(z')}: e^{\alpha(g_\varepsilon * Y_t)(z')} |\rho_{-\ell}(z) P_t \partial^\beta_z g_\varepsilon(z - z')|^2 \, \mathrm{d}z'.$$

By inequality (6.2) in Section 6 of [Tri06], we have, for some $\kappa > 0$,

$$\rho_{-\ell}(z) \lesssim \rho_{-\ell}(z - \tilde{z})(1 + |\tilde{z}|)^\kappa.$$

Therefore,

$$\begin{aligned}
|\rho_{-\ell}(z) P_t \partial^\beta_z g_\varepsilon(z - z')| &= \left| \int_{\mathbb{R}^2} \frac{1}{2\pi t} e^{-\frac{|\tilde{z}|^2}{2t} - m^2 t} \partial^\beta_z g_\varepsilon(z - z' - \tilde{z}) \rho_{-\ell}(z) \, \mathrm{d}\tilde{z} \right| \\
&\leq \int_{\mathbb{R}^2} \frac{1}{2\pi t} e^{-\frac{|\tilde{z}|^2}{2t} - m^2 t} |\partial^\beta_z g_\varepsilon(z - z' - \tilde{z})| \rho_{-\ell}(z) \, \mathrm{d}\tilde{z} \\
&\leq \int_{\mathbb{R}^2} \frac{1}{2\pi t} e^{-\frac{|\tilde{z}|^2}{2t} - m^2 t} |\partial^\beta_z g_\varepsilon(z - z' - \tilde{z})| \rho_{-\ell}(z - z' - \tilde{z}) (1 + |z' + \tilde{z}|)^\kappa \, \mathrm{d}\tilde{z}.
\end{aligned}$$

By the compact support of $g_\varepsilon$, we have

$$\begin{aligned}
|\rho_{-\ell}(z) P_t \partial^\beta_z g_\varepsilon(z - z')| &\leq C_{\ell, g, \varepsilon, \beta} \int_{\mathbb{R}^2} \frac{1}{2\pi t} e^{-\frac{|\tilde{z}|^2}{2t} - m^2 t} (1 + |z' + \tilde{z}|)^\kappa \, \mathrm{d}\tilde{z} \\
&= C_{\ell, g, \varepsilon, \beta} e^{-m^2 t} \mathbb{E}_W[(1 + |z' + W_t|)^\kappa],
\end{aligned}$$

where $W$ is some two-dimensional Brownian motion. Thus, if we take $\kappa$ large enough, we get

$$\begin{aligned}
|\rho_{-\ell}(z) P_t \partial^\beta_z g_\varepsilon(z - z')| &\leq C_{\kappa, \ell, g, \varepsilon, \beta} e^{-m^2 t} (1 + |z'|^\kappa + \mathbb{E}[|W_t|^\kappa]) \\
&\leq C_{\kappa, \ell, g, \varepsilon, \beta, \delta} e^{(\delta - m^2) t} (1 + |z'|^\beta).
\end{aligned}$$

Then,

$$F_\beta(z) \rho_{-2\ell}(z) \leq C_{\kappa, \ell, g, \varepsilon, \beta, \delta} \int \alpha^2 f_\varepsilon(z') :e^{\alpha(g_\varepsilon * X_t)(z')}: e^{\alpha(g_\varepsilon * Y_t)(z')} e^{(\delta - m^2) t} (1 + |z'|^\kappa) \, \mathrm{d}z',$$

which is finite, since $f_\varepsilon$ is compactly supported. □